\documentclass[11pt]{article}
\usepackage[all]{xy}
\usepackage{epsfig}

\usepackage{amsfonts,amssymb,amsthm,amsmath,amscd}
\usepackage[mathscr]{eucal}
\usepackage{mathrsfs}
\usepackage{hyperref}
\usepackage{url}
\usepackage{appendix}

      \makeatletter

\usepackage{titlesec}

\titleformat{\subsection}[runin]
{\normalfont\large\bfseries}{\thesubsection}{1em}{}

\titleformat{\subsubsection}[runin]
{\normalfont\large\bfseries}{\thesubsubsection}{1em}{}

\theoremstyle{plain}
\newtheorem{thm}[subsubsection]{Theorem}
\newtheorem{thm*}{Theorem}
\newtheorem{cor}[subsubsection]{Corollary}
\newtheorem{lem}[subsubsection]{Lemma}
\newtheorem{prop}[subsubsection]{Proposition}
\newtheorem{conj}[subsubsection]{Conjecture}

\theoremstyle{definition}
\newtheorem{defn}[subsubsection]{Definition}

\theoremstyle{remark}
\newtheorem{rem}[subsubsection]{Remark}

\numberwithin{equation}{subsubsection}

\newcommand{\N}{\mathbb N}
\newcommand{\Z}{\mathbb Z}
\newcommand{\Q}{\mathbb Q}
\newcommand{\R}{\mathbb R}
\newcommand{\C}{\mathbb C}
\newcommand{\A}{\mathbb A}

\newcommand{\F}{\mathbb F}

\newcommand{\Fpb}{\bar{\mathbb{F}}_{p}}
\newcommand{\Zp}{{\mathbb Z}_p}
\newcommand{\Qp}{{\mathbb Q}_p}
\newcommand{\Ql}{{\mathbb Q}_l}
\newcommand{\Qv}{{\mathbb Q}_v}
\newcommand{\Qb}{\overline{\mathbb Q}}
\newcommand{\Qpb}{\overline{\mathbb{Q}}_p}
\newcommand{\Qlb}{\overline{\mathbb{Q}}_l}
\newcommand{\Qvb}{\overline{\mathbb{Q}}_v}
\newcommand{\et}{\text{\'et}} 
\newcommand{\nr}{\mathrm{ur}} 
\newcommand{\Gal}{\mathrm{Gal}} 

\newcommand{\Hom}{\mathrm{Hom}}
\newcommand{\Nm}{\mathrm{Nm}} 

\newcommand{\Int }{\mathrm{Int}} 
\newcommand{\Cent}{\mathrm{Cent}} 
\newcommand{\isom}{\stackrel{\sim}{\rightarrow}} 
\newcommand{\ra}{\rightarrow}
\newcommand{\hra}{\hookrightarrow}
\newcommand{\cO}{\mathcal{O}} 
\newcommand{\Gm}{\mathbb{G}_{\mathrm{m}}}

\newcommand{\dS}{\mathbb{S}} 
\newcommand{\Sh}{\mathrm{Sh}} 
\newcommand{\sS}{\mathscr{S}} 
\newcommand{\DD}{\mathbb{D}} 
\newcommand{\cris}{\mathrm{cris}} 
\newcommand{\Betti}{\mathrm{B}} 
\newcommand{\dR}{\mathrm{dR}} 

\begin{document}

\title{Non-emptiness of Newton strata of Shimura varieties of Hodge type}
\author{Dong Uk Lee}
\date{}

\maketitle

\begin{abstract}
For a Shimura variety of Hodge type with hyperspecial level at a prime $p$, the Newton stratification on its special fiber at $p$ is a stratification defined in terms of the isomorphism class of the Dieudonne module of parameterized abelian varieties endowed with a certain fixed set of Frobenius-invariant crystalline cycles (``$F$-isocrystal with $G_{\Qp}$-structure"). There has been a conjectural group-theoretic description of the F-isocrystals that are expected to show up in the special fiber. We confirm this conjecture by two different methods. More precisely, for any $F$-isocrystal with $G_{\Qp}$-structure that is expected to appear (in a precise sense), first we construct a special point which has good reduction and whose reduction has associated $F$-isocrystal equal to given one. Secondly, we produce a Kottwtiz triple (with trivial Kottwitz invariant) with the $F$-isocrystal component being the given one. According to a recent result of Kisin which establishes the Langlands-Rapoport conjecture, such Kottwitz triple arises from a point in the reduction.
\end{abstract}

\tableofcontents

\section{Introduction}

Fix a prime $p>0$, and for $g\geq1$, let $\mathcal{A}_g$ be the moduli space of principally polarized abelian varieties of dimension $g$ in characteristic $p$. Then, the Newton stratification on $\mathcal{A}_g$ is the stratification such that each stratum consists of points $x=(A_x,\lambda_x)$ whose associated $p$-divisible group $(A_x[p^{\infty}],\lambda_x[p^{\infty}])$ with quasi-polarization is quasi-isogenous to a fixed one. The Dieudonn\'e module $\DD(X)$ of a $p$-divisible group $X$ over a perfect field $k$ of characteristic $p>0$ provides an equivalence of categories between the isogeny category of (quasi-polarized) $p$-divisible groups and the category of (quasi-polarized) $F$-isocrystals. An \textit{$F$-isocrystal} (or simply \textit{isocrystal}) over a perfect field $k$ of $\mathrm{char} k=p$ is a finite-dimensional $L(k)$-vector space $M$ with a $\sigma$-linear bijective operator $\Phi:M\ra M$, where $L(k):=\mathrm{Frac}(W(k))$ and $\sigma$ is its Frobenius automophism. According to the Dieudonn\'e-Manin classification, an isocrystal over an algebraically closed field is determined by the slope sequence of $\Phi$. The latter combinatorial datum in turn can be faithfullly represented by a lower-convex (piecewise-linear) polygon with integral break points lying in the first quadrant of $\Z^2$, which will be called \textit{Newton polygon} hereafter. In other words, the Newton polygon of $\DD(A_x[p^{\infty}])$ determines the quasi-isogeny class of $A_x[p^{\infty}]$.
The existence of a quasi-polarization forces the Newton polygon of $\DD(A_x[p^{\infty}])$ to be ``symmetric". Then, a natural question arises which symmetric Newton polygon can be the Newton polygon of a point of $\mathcal{A}_g$. In fact, the Newton polygons arising from points of $\mathcal{A}_g$ in the just described way meet two more restrictions. First, when the initial point is located at the origin, the end point is always $(2g,g)\in \Z^2$. Secondly, it ``lies above" the \textit{ordinary} Newton polygon, which by definition is the Newton polygon having two slopes $(0,1)$ with same multiplicities $g$.
And, it has long been known (after it was conjectured by Manin) that \emph{any} Newton polygon subject to these two restraints is the Newton polygon of a point on $\mathcal{A}_g$; a proof can be found e.g. in \cite[Cor.7.8]{Oort13}, which uses the Honda-Tate theory, thus gives examples defined over finite fields.

The main result of this article is a generalization of this fact to more general moduli spaces of abelian varieties, namely to Shimura varieties of Hodge type. 

Let $(G,X)$ be a Shimura datum of Hodge type. This means that there exists an embedding of Shimura data $(G,X)\hra (\mathrm{GSp}(W,\psi),\mathfrak{H}_g^{\pm})$ into a Siegel Shimura datum, in the sense that there exists an embedding $\rho_{W}:G\hra \mathrm{GSp}(W,\psi)$ of $\Q$-groups which sends each morphism in $X$ to a member of $\mathfrak{H}_g^{\pm}$; we fix such an embedding. For a compact open subgroup $K\subset G(\A_f)$ (which will be tacitly assumed sufficiently small), the associated (canonical model over the reflex field $E(G,X)$ of the) Shimura variety $\Sh_K(G,X)$ parameterizes polarized abelian varieties endowed with a fixed set of (absolute) Hodge cycles and a $K$-level structure, and carries a universal family of abelian varieties $A\ra S_K:=\Sh_K(G,X)$ equipped with a similar set of (absolute Hodge) tensors on $H^1_{\dR}(A/S_K)$. To study its reduction modulo $p$, we fix a prime $\wp$ of $E:=E(G,X)$ and let $\cO:=\cO_{(\wp)}$ be the localization at $\wp$ of the ring of integers of $E$ with residue field $\kappa(\wp)$. To obtain an integral model over $\cO$ with good reduction, we assume that $G\otimes_{\Q}\Qp$ is unramified, which is the same as that there exists a reductive group scheme $G_{\Zp}$ over $\Zp$ with generic fibre $G\otimes\Qp$. We choose one such model $G_{\Zp}$ and set $K_p:=G_{\Zp}(\Zp)$ (such compact open subgroups of $G(\Qp)$ are called \textit{hyperspecial}). For $K$, we take $K=K_p\times K^p$ for a (sufficiently small) compact open subgroup $K^p\subset G(\A_f^p)$. Then, by Vasiu (\cite{Vasiu99}) and Kisin (\cite{Kisin10}), it is known that there exists a smooth integral model $\sS_K:=\sS_K(G,X)$ over $\cO$ with generic fibre $S_K$, which is furthermore uniquely characterized by a ``Neron-extension property", and has a universal abelian scheme $\mathcal{A}$ over it by construction  (see $\S$3.1 for details). To define the $F$-isocrystal attached to a point of the reduction $\sS_K\otimes_{\cO}\kappa(\wp)$, we fix (for simplicity) an algebraically closed field $k$ of finite transcendence degree over $\Fpb$. Let $z$ be a point of $\sS_K\otimes_{\cO}\kappa(\wp)$ defined over $k$. By smoothness of $\sS_K$ and the moduli interpretation of $S_K$, the Dieudonne module $\DD(\mathcal{A}_z)$ is supplied with a set of Frobnius-invariant tensors $\{t_{\alpha,z}\}_{\alpha\in\mathcal{J}}$, and there exists an $L(k)$-isomorphism between the dual space $W^{\vee}\otimes L(k)$ and the Dieudonne module $\DD(\mathcal{A}_z)$ which matches the $G$-invariant tensors on $W^{\vee}$ and the tensors $\{t_{\alpha,z}\}_{\alpha\in\mathcal{J}}$ on $\DD(\mathcal{A}_z)$, which is thus canonically determined up to the action of $G_{L(k)}$ on $W^{\vee}_{L(k)}$ (cf. $\S$3.2). Choosing such an isomorphism, we transport the Frobenius operator $\Phi$ to $W^{\vee}_{L(k)}$ and get an element $b\in G(L(k))$ by $\Phi=\rho_{W^{\vee}}(b)(\mathrm{id}_{W^{\vee}}\otimes\sigma)$, where $\rho_{W^{\vee}}:G\hra \mathrm{GL}(W^{\vee})$ is the contragredient representation of the symplectic representation $\rho_W$ fixed above. Recall that two elements $g_1$, $g_2$ of $G(L(k))$ are called \textit{$\sigma$-conjugate} if there exists $h\in G(L(k))$ with $g_2=hg_1\sigma(h)^{-1}$. Then, the $\sigma$-conjugacy class $\overline{b}$ of $b$ is independent of the choice of an isomorsphim $W^{\vee}_{L(k)}\ra \DD(A_x)$ just explained, and the isomorphism class of the $F$-isocrystal $(\DD(\mathcal{A}_z),\Phi)$ equipped with the Frobenius-invariant tensors $\{t_{\alpha,z}\}_{\alpha\in\mathcal{J}}$, called \textit{$F$-isocrystal with $G_{\Qp}$-structure}, corresponds to a unique $\sigma$-conjugacy class of elements of $G(L(k))$. If $B(G_{\Qp})$ denotes the set of the $\sigma$-conjugacy classes of elements of $G(L(k))$, there exists a subset $B(G_{\Qp},X)$ of $B(G_{\Qp})$ which is expected to be the set of isocrystals with $G_{\Qp}$-structure coming from points on $\sS_K\otimes_{\cO}\kappa(\wp)$. We remark that this is a subset of $B(G_{\Qp})$ defined by certain two conditions that are analogues of those discussed above in the Siegel case, and it is easy to see that for any point of $\sS_K\otimes_{\cO}\kappa(\wp)$, its associated $F$-isocrystal with $G_{\Qp}$-structure belongs to this subset $B(G_{\Qp},X)$. A detailed presentation of the discussion in this paragraph is given in Section 3.

For $\overline{b}\in B(G_{\Qp},X)$, let $\sS_{\overline{b}}$ be the subset of points of $\sS_K\otimes_{\cO}\kappa(\wp)$ whose associated $F$-isocrystal with $G_{\Qp}$-structure is $\overline{b}$. The main result of this paper is then that for every $\overline{b}\in B(G_{\Qp},X)$, $\sS_{\overline{b}}$ is non-empty. We present two proofs. The main idea of the first proof is to find a special point with required reduction property. We recall that a point $[h,g_f\cdot K]\in \Sh_K(G,X)(\C)=G(\Q)\backslash X\times G(\A_f)/K$ is called \textit{special}  if $h\in X$ factors through $T_{\R}$ for some maximal $\Q$-torus $T$ of $G$. Special points are known to be defined over $\Qb\subset\C$ and correspond to abelian varieties with complex multiplication. More precisely, in the first proof, we construct a special point which extends to a point $\sS_K(\cO_{\iota_p})$, where $\cO_{\iota_p}$ is the valuation ring of $\Qb$ defined by chosen embedding $\iota_p:\Qb\hra\Qpb$ (in which case, we will say that it has \textit{good reduction at $\iota_p$}) and such that the isocrystal of its reduction in $\sS_K(\Fpb)$ is the given isocrystal. 

\begin{thm} [Thm. \ref{thm:existence_of_special_point_with_good_reduction_and_prescribed_given_isocrystal}]

Let $(G,X)$ be a Shimura datum and $p$ a rational prime. Assume that $G_{\Qp}$ is unramified. Fix a hyperspecial subgroup $K_p$ of $G(\Qp)$ and a sufficiently small compact open subgroup $K^p\subset G(\A_f^p)$; put $K=K_p\times K^p$. Choose a prime $\wp$ of $E(G,X)$ above $p$, and an embedding $\iota_p:\Qb\hra\Qpb$ inducing $\wp$. 

Then, for every $\overline{b}\in B(G_{\Qp},X)$, there exists a special sub-Shimura datum $(T,h\in \Hom(\dS, T_{\R})\cap X)$ such that $T_{\Qp}$ is unramified and the (unique) hyperspecial subgroup of $T(\Qp)$ is contained in $K_p$, and that the $\sigma$-conjugacy class of $\mu_h^{-1}(p)\in G(\Qp^{\nr})$ equals $\overline{b}$.

In particular, if $(G,X)$ is of Hodge type, for every $g_f\in G(\A_f)$, the special point $[h,g_f]_{K}\in \Sh_{K}(G,X)(\Qb)$ has good reduction at $\iota_p$ and the $F$-isocrystal of its reduction in $\sS_{K}(\Fpb)$ is $\overline{b}$.
\end{thm}

The non-emptiness of Newton strata has been conjectured by Fargues \cite[Conjecture 3.1.1]{Fargues04} and Rapoport \cite[Conjecture 7.1]{Rapoport05}. And, some partial cases confirming this conjecture were known. In this introduction, we only mention previous works in two directions.
For PEL-type Shimura varieties (which form a subclass of Hodge-type Shimura varieties), this conjecture (i.e. non-emptiness of \emph{all} Newton strata in $B(G_{\Qp},X)$) was proved by C.-F. Yu \cite{Yu05} in the Lie-type $C$ cases, then by Viehmann-Wedhorn \cite[Thm.1.6]{VW13} in general cases, and by Kret \cite{Kret13} for some simple groups of Lie type $A$ or $C$.
For a general Hodge-type Shimura variety, to the best of the author's knowledge, there are two partial results. First, Wortmann \cite{Wortmann13} showed non-emptiness of \textit{$\mu$-ordinary locus}; the $\mu$-ordinary locus is a generalization of the (usual) ordinary locus in $\mathcal{A}_g$. And, for \emph{projective} Shimura varieties of Hodge type, Koskivirta \cite{Koskivirta14} proved the conjecture (non-emptiness of \emph{all} Newton strata). 

While all these previous results on the non-emptiness of Newton strata, except for that of Kret, are obtained by characteristic-$p$ geometric methods, which are often of limited applicability for Hodge-type Shimura varieties, we draw reader's attention to that the first claim of our theorem above is a purely group-theoretic statement on Shimura \emph{data} without any reference to geometry of Shimura variety. When the Shimura datum is of Hodge type, this group-theoretic statement carries a geometric meaning of the second claim (by Lemma \ref{lem:isocrystal_of_special_pt}). In this geometric interpretation of the group-theoretic statement of the first claim, the only geometric property that we use about the canonical integral model $\sS_K$ is the Neron-extension property, which is in the end a consequence of smoothness of $\sS_K$. 

We also point out that our method of proof of the theorem, more precisely that of constructing a special point with certain prescribed properties, allows us to obtain some finer results. First, as already apparent, we prove more than just non-emptiness of Newton strata, namely we show existence of a special point whose reduction lying in any given Newton stratum. As another example, we also prove a generalization of a result of Wedhorn \cite[(1.6.3)]{Wedhorn99} (the next result in the PEL-type cases of Lie type $A$ or $C$) (by argument different from the original one of Wedhorn, cf. Remark \ref{rem:mu-ordinary_locus}).

\begin{cor} (Cor. \ref{cor:non-emptiness_criterion_of_ordinary_locus})
Let $(G,X)$ be a Shimura datum of Hodge type. Suppose that $G_{\Q_p}$ is unramified and choose a hyperspecial subgroup $K_p$ of $G(\Qp)$. Let $\wp$ be a prime of $E(G,X)$ above $p$. Then the reduction $\sS_{K_p}(G,X)\times\kappa(\wp)$ has non-empty ordinary locus if and only if $\wp$ has absolute height one (i.e. $E(G,X)_{\wp}=\Q_p$).
\end{cor}

Our second proof of the non-emptiness of Newton strata is by construction of a Kottwtiz triple. A Kottwitz triple is a triple of elements in the product 
\begin{equation*}
(\gamma_0;\gamma=(\gamma_l)_{l\neq p},\delta)\in G(\Q)\times G(\A_f^p)\times G(\Qp^{\nr}),
\end{equation*}
satisfying certain compatibility conditions among themselves and some other conditions. It has been expected that one can attach such triple to every point of $\sS_K$ over a finite field, such that $\gamma_l (l\neq p)$ is the (relative) Frobenius acting on the $l$-adic etale cohomology and $\delta$ is the (absolute) Frobenius on the crystalline cohomology of the abelian variety corresponding to the point. This was established by Kisin in the course of his proof of the Langlands-Rapoport conjecture (\cite{Kisin13}). In fact, this Langlands-Rapoport conjecture implies that every Kottwitz triple meeting suitable conditions (especially, vanishing of the Kottwitz invariant) also comes from a point over a finite field. Therefore, in view of this, our second proof is deduced from the following result:

\begin{thm} [Thm. \ref{thm:existence_of_Kottwitz_triple_for_any_b_in_B(G,X)}, Cor. \ref{for:second_proof_of_non-emptiness}]
Keep the assumptions and notations from the above theorem, and assume further that the derived group $G^{\mathrm{der}}$ of $G$ is simply connected.

Then, for every $\overline{b}\in B(G_{\Qp},X)$, there exists a Kottwitz triple $(\gamma_0;\gamma,\delta)\in G(\Q\times \A_f^p\times \Qp^{\nr})$ with $\delta\in\overline{b}$ and having trivial Kottwitz invariant.
In particular, there exists a point $x\in Sh_{K_p}(G,X)(\Fpb)$ whose associated $F$-isocrystal is $\overline{b}$, where $\Fpb$ is an algebraic closure of $\kappa(\wp)$.
\end{thm}

We refer readers to $\S$6 for the definitions of Kottwitz triple and Kottwitz invariant, and also to Theorem \ref{thm:existence_of_Kottwitz_triple_for_any_b_in_B(G,X)} for a more precise statement of the result in this theorem.

Our two proofs of the non-emptiness of Newton strata are not unrelated. Both are based on the fact that for every $\overline{b}\in B(G_{\Qp},X)$, the condition $(\dagger)$ in Theorem \ref{thm:existence_of_special_SD_with_given_F-isocrystal} holds. We verify this in Proposition \ref{prop:existence_of_Levi_and_Hodge_cocharacter} in two ways.
In fact, this was proved in \cite{LR87}, Lemma 5.11. However, we also present our own proof which was found before the proof of Langalnds-Rapoport came to our attention. Our proof is similar to theirs, but also differs in details. 

Finally, we remark that since our methods of proof of the above results are purely group-theoretic (making heavy use of Galois cohomology of algebraic groups), they are very likely to apply to more general Shimura varieties, particularly to Shimura varieties of abelian-type and/or with non-hyperspecial level. We did not pursue such ideas in this article, though.

This article is organized as follows. In Section 2, we give a brief account of the theory of Newton map for reductive groups over $p$-adic fields and provide a collection of some basic facts. In Section 3, we review the construction of the canonical integral model $\sS_{K_p}(G,X)$ of a Hodge-type Shimura variety $\Sh_{K_p}(G,X)$ with hyperspecial level, due to Vasiu and Kisin, and the Newton stratification on its (good) reduction $\sS_{K_p}(G,X)\otimes\kappa(\wp)$. Here, we give a criterion when a special point on $\Sh_K(G,X)$ has good reduction at given embedding $\iota_p:\Qb\hra\Qpb$, and identifies the $F$-isocrystal of its reduction in terms of the defining special Shimura datum. 
Section 4 is devoted to the first proof of the non-emptiness of Newton strata, by constructing a special point with good reduction at $\iota_p$ and with prescribed $F$-isocrystal. 
Our construction is carried out, assuming that given isocrystal $\overline{b}\in B(G_{\Qp})$ satisfies the condition $(\dagger)$ in Theorem \ref{thm:existence_of_special_SD_with_given_F-isocrystal}. Then, we verify that the condition $(\dagger)$ holds for \emph{every} isocrystal $\overline{b}\in B(G_{\Qp})$. 
We first present the proof of Langlands-Rapoport \cite[Lemma5.11]{LR87} of this, and then our own argument. Combining these results, we get our first proof of non-emptiness of Newton strata. As a corollary, we also obtain a criterion of when the good reduction of a Shimura variety with hyperspecial level has an ordinary point. In the final section, Section 5, we give our second proof of the non-emptiness of Newton strata, by constructing a Kottwtiz triple with prescribed isocrystal component. We begin by a review of the definitions of Kottwitz triple and Kottwitz invariant. Our construction of such Kottwitz triple (with trivial Kottwitz invariant) also uses the aforementioned fact that every isocrystal $\overline{b}\in B(G_{\Qp},X)$ satisfies the condition $(\dagger)$.

\

\textbf{Acknowledgements}
This work was supported by IBS-R003-D1. Also, parts of this research were conducted during the author's stay at PMI (Pohang Mathematics Institute), Postech, and ASARC (Algebraic Structures and Applications Research Center), Kaist, with support from the National Research Foundation of Korea (NRF) funded by the Korea government (MSIP) (No. 2011-0030749 and No. 2012-0000798).

\

\textbf{Notations}

Throughout this paper, $\Qb$ denotes the algebraic closure of $\Q$ inside $\C$ (so $\Qb$ has a previliged embedding into $\C$). 

For a (connected) reductive group $G$ over a field, we let $G^{\mathrm{sc}}$ be the universal covering of its derived group $G^{\mathrm{der}}$, and for a (general linear algebraic) group $G$, $Z(G)$ and $G^{\mathrm{ad}}$ denote its center and the adjoint group $G/Z(G)$, respectively.



\section{$F$-isocrystals and Newton maps for reductive groups over $p$-adic fields}

In this section, we give a brief review of the Newton map which is defined for connected reductive groups over $p$-adic fields. For more details on the material discussed in this section, see, e.g.  
\cite{Kottwitz85}, \cite{RR96}, \cite{Rapoport05} and references therein.

\subsection{Preliminaries} 
Let $G$ be a connected reductive group over an arbitrary field $F$ of characteristic zero.

\subsubsection{} 
Any Borel pair $(T,B)$ of $G_{\overline{F}}$ (i.e. $T$ is a maximal torus of $G_{\overline{F}}$ and $B$ is a Borel subgroup of $G_{\overline{F}}$ containing $T$) determines a based root datum 
\begin{equation*}\mathcal{BR}(G,T,B)=(X^{\ast}(T),R^{\ast}(T),X_{\ast}(T),R_{\ast}(T),\Delta(T,B)),\end{equation*}
where $\mathcal{R}(G,T)=(X^{\ast}(T),R^{\ast}(T),X_{\ast}(T),R_{\ast}(T))$ is the root
datum of $(G,T)$ (the usual notation) and $\Delta(T,B)$ is the root base of $R^{\ast}(T)$ determined by $B$. A based root datum of $G$ has a canonical Galois action
$\tau:\Gal(\overline{F}/F)\rightarrow\mathrm{Aut}(\mathcal{BR}(G,T,B))$
\cite[15.5]{Springer98} such that for $x\in X^{\ast}(T),\
\sigma\in\Gal(\overline{F}/F)$,
\begin{equation*}\tau(\sigma)x=g\sigma(x)g^{-1},\end{equation*}
where $g\in G(\overline{F})$ satisfies $g\sigma(T,B)g^{-1}=(T,B)$; this action can be also understood as a natural action of $\Gal(\overline{F}/F)$ on the orbit space $X_{\ast}(T)/\Omega(G,T)$, where $\Omega(G,T)=N_G(T)/Z_G(T)$ is the Weyl group of $G$, i.e. $N_G(T)$ (resp. $Z_G(T)$) is the normalizer (resp. centralizer) of $T$ in $G$.

\subsubsection{} 
Let $\overline{C}(T,B)\subset X_{\ast}(T)_{\Q}$ be the closed Weyl chamber corresponding to a root base $\Delta(T,B)$ of $R^{\ast}(T)$:
\begin{equation*}\overline{C}(T,B)=\{ x\in X_{\ast}(T)_{\Q}\ |\ \langle\alpha,x\rangle\geq0,\ \forall \alpha\in\Delta(T,B)\}.\end{equation*}
This is a fundamental domain for the action of the Weyl group on $X_{\ast}(T)_{\Q}$ and is stable under the action of $\Gal(\overline{F}/F)$ just defined. Note that when $T$ is defined over $F$, we also have a ``naive" Galois action on $X_{\ast}(T)$ defined in terms of the $F$-structure of $T$ only.
In this article, we call the former action on $\overline{C}(T,B)$ \emph{canonical} action, to distinguish it from the naive action on $X_{\ast}(T)$. When $G$ is \emph{quasi-split}, i.e. when there is a Borel subgroup $B$ of $G$ defined over $F$, the naive action on $X_{\ast}(T)$ restricts to the canonical action on $\overline{C}(T,B)$. 
For any other Borel pair $(T',B')$, there is $g\in G(\overline{F})$ such that $(T',B')=g(T,B)g^{-1}$ and the resulting isomorphism of based root data $\mathcal{BR}(G,T,B)\isom \mathcal{BR}(G,T',B')$ is independent of the choice of $g$ and is $\Gal(\overline{F}/F)$-equivariant. The \emph{based root datum} $\mathcal{BR}(G)=(X^{\ast},R^{\ast},X_{\ast},R_{\ast},\Delta)$ of $G$ is the projective limit of all such based root data $\mathcal{BR}(G,T,B)$ (where $(T,B)$ runs through the set of Borel pairs of $G_{\overline{F}}$). This is equipped with a canonical Galois action $\Gal(\overline{F}/F)\rightarrow\mathrm{Aut}(\mathcal{BR}(G))$.

\subsubsection{} 
For an algebraically closed extension $\overline{K}$
of $F$, let $\mathcal{C}_G(\overline{K})$ denote the set of
$G(\overline{K})$-conjugacy classes of (algebraic group)
homomorphisms $\mathbb{G}_{m,\overline{K}}\rightarrow
G_{\overline{K}}$. Then, for any maximal torus $T$ of $G_{\overline{F}}$, we have a canonical identification
\begin{equation*}\mathcal{C}_G(\overline{K}) \cong X_{\ast}(T)/\Omega.\end{equation*}
In particular, we have
$\mathcal{C}_G(\overline{K})=\mathcal{C}_G(\overline{F})$ for any
algebraically closed extension $\overline{K}$ of $F$.

\subsection{Newton maps} 
Next, we review the definition of Newton map. 
Let $k$ be an algebraically closed field of characteristic $p$, and let $W(k)$ (resp. $K_0=L(k)$) be its Witt vector ring (resp. the fraction field of $W(k)$). Fix an algebraic closure $\overline{K_0}$ of $K_0$. In this subsection, let $F$ be a finite extension of $\Q_p$ in $\overline{K_0}$, and let $L$ be the composite of $K_0$ and $F$ in $\overline{K_0}$; for our application to Shimura varieties, the most interesting case is when $k=\Fpb$,  $F=\Qp$ so that $L=L(\Fpb)$. Denote by $\sigma$ the Frobenius automorphism on $K_0$ (and its extension to an element of $\mathrm{Aut}(L/F)$). Let $G$ be a connected reductive group $G$ over $F$. Let $\Gal_F:=\Gal(\overline{F}/F)$.

\subsubsection{} 
Let $B(G)$ be the set of $\sigma$-conjugacy classes in $G(L)$:
\begin{equation*}B(G)=G(L)/\sim,\end{equation*}
where two elements $x,y$ of $G(L)$ are $\sigma$-conjugated, denoted
by $x\sim y$, if $x=gy\sigma(g)^{-1}$ for some $g\in G(L)$. An element of $B(G)$ is called ($F$-)\textit{isocrystal with $G$-structure}.

\subsubsection{} 
Let $\mathbb{D}$ be the pro-algebraic torus with character group $\Q$. We define $\mathcal{N}(G)$ to be the set of $\sigma$-invariants in the set of conjugacy classes of homomorphisms $\mathbb{D}_{L}\rightarrow G_{L}$:
\begin{equation*}\mathcal{N}(G)=(\mathrm{Int}\ G(L)\ \backslash\ \Hom_{L}(\mathbb{D},G))^{\langle\sigma\rangle},\end{equation*}
where $\langle\sigma\rangle$ denotes the infinite cyclic group generated by $\sigma$ (and is endowed with the discrete topology). For example, when $G=T$ is a torus, we have $\mathcal{N}(T)=X_{\ast}(T)^{\Gal_F}\otimes\Q$. 

If $T\subset G$ is a maximal $F$-torus with (absolute) Weyl group $\Omega$, there is a natural identification
\begin{equation*}\mathcal{N}(G)=(X_{\ast}(T)_{\Q}/\Omega)^{\Gal_F}.\end{equation*}
In particular, if $G$ is quasi-split and $T$ is a maximal $F$-torus contained in a Borel subgroup $B$ defined over $F$, $\mathcal{N}(G)$ can be identified with the subset of $\Gal_F$-invariants in $\overline{C}(T,B)$. 

\begin{thm} \cite[Theorem 1.8]{RR96} \label{thm:Thm.1.8-RR96}
For $b\in G(L)$, there exists a unique element $\nu=\nu_b\in\Hom_L(\mathbb{D},G)$ for which there are an integer $s>0$, an element $c\in G(L)$ and a uniformizing element $\pi$ of $F$ such that:
\begin{itemize}\addtolength{\itemsep}{-4pt}
\item[(i)] $s\nu\in\Hom_L(\Gm,G)$.
\item[(ii)] $\Int (c)\circ s\nu$ is defined over the fixed field of $\sigma^s$ in $L$.
\item[(iii)] $c\cdot b\cdot \sigma(b)\cdots\sigma^{s-1}(b)\cdot \sigma^{s}(c)^{-1}=c\cdot 
(s\nu)(\pi)\cdot c^{-1}$.
\end{itemize}
The map $b\mapsto \nu_b$ has the following properties.
\begin{itemize}\addtolength{\itemsep}{-4pt}
\item[(a)] $\nu_{\sigma(b)}=\sigma(\nu_b)$.
\item[(b)] $gb\sigma(g)^{-1}\mapsto \Int (g)\circ \nu,\ g\in G(L)$.
\item[(c)] $\nu_b=\Int (b)\circ\sigma(\nu_b)$.
\end{itemize}
\end{thm}
By the properties (b) and (c), we obtain a well-defined map $\nu_G:B(G)\rightarrow \mathcal{N}(G)$. 
The induced functorial map
$\overline{\nu}$, defined on the category of connected reductive
groups, is called the \textit{Newton map}:
\begin{equation*}\overline{\nu}:B(\cdot)\rightarrow \mathcal{N}(\cdot)\ ;\
\overline{\nu}_{G}(\overline{b})=\overline{\nu}_{b},\quad
b\in\overline{b}.\end{equation*} Here, $\overline{b}$ (resp.
$\overline{\nu}_{b}$) is the $\sigma$-conjugacy class of $b$ (resp.
the conjugacy class of $\nu_b$). When $G=\mathrm{GL}_h$, $B(G)$ classifies the isomorphism classes of $\sigma-L$-spaces of height $h$, and the Newton map sends a $\sigma-L$-space of height $h$ to its Newton polygon which is represented in $\mathcal{N}(G)$ by the corresponding slope homomorphism (cf. Example 1.10 of \cite{RR96}).

\subsubsection{} \label{subsubsec:algebraic_fundamental_group}
Recall the canonical action of $\Gal_F$ on a based root datum $\mathcal{BR}(G,T,B)$ (cf. Subsection 2.1). The induced action of $\Gal_F$ on
\begin{equation*}\pi_1(G,T):=X_{\ast}(T)/\Sigma_{\alpha\in\Delta(G,T)}\Z\alpha^{\vee}\end{equation*}
is independent of the choice of $B$, where $\Phi(G,T)$ denotes the set of roots of $T$ and for $\alpha\in\Phi(G,T)$, $\alpha^{\vee}$ is the coroot of $\alpha$. Moreover, for another maximal torus $T'$ of $G$, we have a canonical isomorphism $\pi_1(G,T)=\pi_1(G,T')$ of $\Gal_F$-modules. The algebraic fundamental group $\pi_1(G)$ of $G$ is defined to be this common $\Gal_F$-module; The functor is an exact functor from the category of connected reductive groups over $F$ to the category of finitely generated discrete $\Gamma_F$-modules, cf. \cite[1.13]{RR96}. For a torus $T=G$, we have
\begin{equation*}\pi_1(T)=X_{\ast}(T).\end{equation*}
The functor $\pi_1(\cdot)^{\Gal_F}\otimes\Q$ defines an exact functor from the category of connected reductive groups over $F$ to the category of finite-dimensional $\Q$-vector spaces. There is a canonical isomorphism 
\begin{equation*}X^{\ast}(Z(\widehat{G})^{\Gal_F})=\pi_1(G)_{\Gal_F},\end{equation*}
where $Z(\widehat{G})$ denotes the center of the Langlands dual group of $G$ and 
\begin{equation*}
\pi_1(G)_{\Gal_F}=\pi_1(G)/\langle\tau x-x\ |\ \tau\in\Gal_F,\ x\in\pi_1(G)\rangle
\end{equation*} 
is the group of coinvariants of $\pi_1(G)$ for the canonical action of $\Gal_F$.

\begin{thm} [\cite{RR96}, Theorem 1.15] \label{thm;RR96-Thm.1.15} 
(i) There exists a unique natural transformation
\footnote{In \cite{RR96}, Theorem 1.15, this map is denoted $\gamma$, but here we changed it to $\kappa$, a notation which seems to be more popular in recent literatures.}
\begin{equation*}\kappa:B(\cdot)\rightarrow \pi_1(\cdot)_{\Gal_F}\end{equation*}
of set-valued functors on the category of connected reductive groups
over $F$ such that on $\Gm$ we have
\begin{equation*}\kappa_{\Gm}(\overline{b})=\mathrm{ord}_L(b)\  (b\in L^{\times}),\end{equation*}
where the valuation on $L$ is normalized by
$\mathrm{ord}_L(\pi_L)=1$ for a uniformizer $\pi_L$, and for tori it
is an isomorphism. 

(ii) There is a
unique natural transformation of functors on the category of
connected reductive groups
\begin{equation*}\delta:\mathcal{N}(\cdot)\rightarrow \pi_1(\cdot)^{\Gal_F}\otimes\Q\end{equation*}
such that  for a torus this is the natural identification
\begin{equation*}\mathcal{N}(T)=X_{\ast}(T)^{\Gal_F}\otimes\Q=\pi_1(T)^{\Gal_F}\otimes\Q.\end{equation*}

(iii) The following diagram is functorial (in the sense of pointed sets)
\begin{equation*}
\xymatrix{B(G) \ar[r]^{\overline{\nu}_G} \ar[d]_{\kappa_G} & \mathcal{N}(G) \ar[d]^{\delta_G} \\
\pi_1(G)_{\Gal_F}\ar[r]^(.5){\alpha_G} & \pi_1(G)^{\Gal_F}\otimes\Q.}
\end{equation*}
Here, the map $\alpha_{G}$ in the bottom is given by
\begin{equation*}\overline{\mu}\mapsto |\Gal_F\cdot\mu|^{-1}\sum_{\mu'\in \Gal_F\cdot\mu}\mu'.\end{equation*}
\end{thm}

For a detailed discussion of this theorem, we refer to \cite{Kottwitz85}.


\section{Newton stratification on good reduction of Shimura variety}

In this section, we give an account of the construction of the Newton stratification on the reduction of the integral canonical model of a Hodge-type Shimura variety. For the theory of the integral canonical model of Shimura variety of Hodge type, we follow \cite{Vasiu99}, \cite{Kisin10} (see also \cite{Vasiu08}, \cite{Kisin09} for a survey of the main results).

\subsubsection{} \label{sssec:Hodge_cocharacter:mu_h}
Let $(G,X)$ be a Shimura datum. For a morphism $h:\dS\rightarrow G_{\R}$ in $X$, the associated \emph{Hodge cocharacter}
\begin{equation*}\mu_{h}:\mathbb{G}_{m\C}\rightarrow G_{\C}\end{equation*}
is the composite of $h_{\C}:\dS_{\C}\rightarrow G_{\C}$ and the cocharacter of $\dS_{\C}\cong\prod_{\C\hra\C}\mathbb{G}_{m\C}$ corresponding to the identity embedding $\C\hookrightarrow\C$. Let $c(G,X)$ denote the $G(\C)$-conjugacy class of cocharacters of $G_{\C}$ containing $\mu_h^{-1}$ (\emph{not} $\mu_h$), for an arbitrary $h\in X$. For a maximal torus $T$ of $G_{\Qb}$, we can consider $c(G,X)$ as an element of $X_{\ast}(T)/\Omega$. Alternatively, when we fix a based root datum $\mathcal{BR}(G,T,B)$, $c(G,X)$ has a unique representative in the associated closed Weyl chamber $\overline{C}(T,B)$. We will also identify the conjugacy class $c(G,X)$ with this representative: 
\begin{equation*} \label{eq:representataive_Hodge_cocharacter}
c(G,X)\in \overline{C}(T,B).
\end{equation*}
The \emph{reflex field} $E(G,X)$ of a Shimura datum $(G,X)$ is the field of definition of $c(G,X)\in\mathcal{C}_G(\Qb)$, i.e. the fixed field of the stabilizer of $c(G,X)$ in $\Gal(\Qb/\Q)$; so a reflex field, which is a finite extension of $\Q$, is always a subfield of $\C$, cf. \cite[Lemme1.2.4]{Deligne77}. When $T$ is a torus, the reflex field $E(T,\{h\})$ is just the smallest subfield of $\Qb\subset\C$ over which the single
morphism $\mu_h$ is defined.

\subsection{Integral canonical model} Let $(G,X)$ be a Shimura datum of Hodge type provided with an embedding $\rho:(G,X)\rightarrow(\mathrm{GSp}(W,\psi),\mathfrak{H}^{\pm})$. Consider a compact open subgroup $K$ of $G(\A_f)$ of the form $K=K_pK^p$, where $K_p$ (resp.$K^p$) is a compact open subgroup of $G(\Q_p)$ (resp. of $G(\A_f^p)$). 

\begin{lem} \label{lem:existence_of_H}
Suppose $p>2$. If $K_p$ is a hyperspecial subgroup (i.e. $K_p=G_{\Z_p}(\Z_p)$ inside $G(\Qp)$ for a reductive group scheme $G_{\Z_p}$ over $\Z_p$ with generic fiber $G_{\Qp}:=G\otimes_{\Q}\Qp$), 
for every sufficiently small $K^p$, there exist a lattice $W_{\Z}$ of $W$, compact open subgroups $H_p\subset \mathrm{GSp}(W,\psi)(\Q_p)$, $H^p\subset \mathrm{GSp}(W,\psi)(\A_f^p)$, such that $H:=H_pH^p$ leaves $W_{\widehat{\Z}}:=W_{\Z}\otimes\widehat{\Z}$ stable and that $\rho$ induces an embedding (of weakly canonical models of Shimura varieties)
\begin{equation*}\Sh_K(G,X)\hookrightarrow \Sh_H(\mathrm{GSp}(W_{\Z},\psi),\mathfrak{H}^{\pm})\otimes_{\Q}E(G,X).\end{equation*}
\end{lem}
\textsc{Proof.} This is shown in \cite[(2.3)]{Kisin10}. We give a brief summary of it (also see the introduction of \cite{Vasiu08}). According to \cite[Lem.(2.3.1)]{Kisin10}, there exists a lattice $W_{\Z}$ of $W$ such that $\rho_{\Q_p}:G_{\Qp}\hra \mathrm{GL}(W_{\Qp})$ is induced by a closed embedding $\rho_{\Zp}: G_{\Z_p}\hookrightarrow \mathrm{GL}(W_{\Z_p})$, where $W_{\Z_p}:=W_{\Z}\otimes_{\Z}\Z_p$. Fix such a $W_{\Z}$. Then, for $K^p$ small enough, $K^p$ leaves $W_{\Z}\otimes_{\Z}\widehat{\Z}^{(p)}$ as well. Finally, when we take $H_p\subset \mathrm{GSp}(W,\psi)(\Q_p)$ to be the stabilizer of $W_{\Z_p}$, the existence of $H^p$ with the required property is established in (2.1.2) of \cite{Kisin10} (cf. \cite[Lem. 3.3]{Noot96}). $\qed$


\subsubsection{} For sufficiently small compact open subgroup $H$ of $\mathrm{GSp}(\A_f)$, let $\mathcal{A}_{g,d,H}$ denote the fine moduli scheme over $\Z_{(p)}$ representing the functor which, with any $\Z_{(p)}$-scheme $T$, associates the isomorphism classes of an abelian scheme equipped with a polarization of degree $d$ and a level $H$-structure (cf. \cite{Kottwitz85}). Choose a prime $\wp$ of $E(G,X)$ above $p$ and let $\mathcal{O}_{(\wp)}$ be the localization at $\wp$ of the ring of integers $\mathcal{O}_{E(G,X)}$ of $E(G,X)$. 
After modifying $W_{\Z}$ by a scalar multiple, we may assume that $\psi$ is $\Z$-valued on $W_{\Z}$. Let $W_{\Z}^{\vee}$ be the dual lattice of $W_{\Z}$ (with respect to $\psi$). If we take a compact open subgroup $H$ of $\mathrm{GSp}(\A_f)$ as in Lemma \ref{lem:existence_of_H} and $d:=[W_{\Z}^{\vee}:W_{\Z}]$, there is a natural embedding of $\mathcal{O}_{(\wp)}$-schemes
\begin{equation*}\Sh_{K}(G,X)\hookrightarrow \mathcal{A}_{g,d,H}\otimes_{\Z_{(p)}}\mathcal{O}_{(\wp)}.\end{equation*}
We denote by $\sS_{K}(G,X)$ the
normalization of the closure of $\Sh_{K}(G,X)$ in 
$\mathcal{A}_{g,d,H}\otimes_{\Z_{(p)}}\mathcal{O}_{(\wp)}$. 
Set
\begin{equation*}\sS_{K_p}(G,X)=\varprojlim_{K^p}\sS_{K_pK^p}(G,X),\end{equation*}
where $K^p$ runs over sufficiently small compact open subgroups of $G(\A_f^p)$.


The following theorem was conjectured by Langlands, Milne  (cf. \cite{Milne92}, Conjecture 4.25), and was established by Vasiu (\cite{Vasiu99}, \cite{Vasiu07}, \cite{Vasiu08}) and Kisin (\cite{Kisin09}, \cite{Kisin10}).

\begin{thm} \label{thm:integral_canonical_model}
Retain previous assumptions (particularly, $p>2$ and $K_p$ is hyperspecial). 

(1) The scheme $\sS_{K_p}(G,X)$ is an inverse limit of smooth $\mathcal{O}_{(\wp)}$-schemes with finite \'etale transition maps, whose restriction to $E(G,X)$ can be $G(\A_f^p)$-equivariantly identified with $\Sh_{K_p}(G,X)$. 

(2) The scheme $\sS_{K_p}(G,X)$ has the following extension property: for any regular, formally smooth $\mathcal{O}_{(\wp)}$-scheme $S$, every $E(G,X)$-morphism $S_{E(G,X)}\rightarrow \Sh_{K_p}(G,X)$ extends uniquely over $\mathcal{O}_{(\wp)}$.

\end{thm}

In particular, $\sS_{K_p}(G,X)$ itself is regular and formally smooth over $\mathcal{O}_{(\wp)}$. We call the $\mathcal{O}_{(\wp)}$-(pro)schemes $\sS_{K}(G,X)$ and $\sS_{K_p}(G,X)$ \emph{integral canonical model} of $\Sh_{K}(G,X)$ and $\Sh_{K_p}(G,X)$, respectively.

\begin{rem} \label{rem:Neron_extension_property}
This definition of the extension property is due to Kisin (\cite{Kisin10}, Thm. (2.3.8)). It differs slightly from that of Vasiu, who uses, as test schemes, healthy regular schemes over $\mathcal{O}_{(\wp)}$ (\cite[3.2]{Vasiu99}) instead of regular, formally smooth $\mathcal{O}_{(\wp)}$-schemes. But, when $G_{\Q_p}$ is unramified, every regular, formally smooth scheme over $\mathcal{O}_{(\wp)}$ is also healthy regular over $\mathcal{O}_{(\wp)}$, because $\mathcal{O}_{(\wp)}$ is unramified over $\Z_{(p)}$ (\cite{Milne94}, Cor. 4.7 (a)) and over any d.v.r unramified over $\Z_{(p)}$, every regular, formally smooth scheme is healthy regular (if $p>2$), according to a lemma of Faltings (\cite[3.6]{Moonen98}). Consequently, this difference (and related other minor differences) will not affect our arguments involving the integral canonical models, since our use of the integral canonical models will be only restricted to its extension property  for discrete valuation rings unramified over $\mathcal{O}_{(\wp)}$.
\end{rem}

\subsection{Newton stratification} 
We keep the notations and assumptions from the previous subsection. In particular, we are given a hyperspecial subgroup $K_p$ of $G(\Qp)$, and accordingly a reductive group scheme $G_{\Zp}$ over $\Zp$ with generic fibre $G_{\Qp}$ and such that $G_{\Zp}(\Zp)=K_p$. Also, choose a lattice $W_{\Z}$ of $W$ and a closed embedding 
\begin{equation*}
\Sh_K(G,X)\hra \Sh_H(\mathrm{GSp}(W_{\Z},\psi),\mathfrak{H}^{\pm})\otimes_{\Q}E(G,X)
\end{equation*}
as in Lemma \ref{lem:existence_of_H}. Further, we fix an embedding $\iota_p:E(G,X)\hra\Qpb$ inducing $\wp$ on $E(G,X)$. Let $\sS_{K}\otimes\kappa(\wp)$ and $\sS_{K_p}\otimes\kappa(\wp)$ be the reductions of $\sS_{K}$ and $\sS_{K_p}$ at $\wp$, respectively. In this subsection, we define the Newton stratification on these schemes, following \cite{Rapoport05}, \cite{Vasiu08}, \cite{Wortmann13}.

\subsubsection{}
For a vector space or a free module over a d.v.r $V$, let $\mathcal{T}(V)$ be the tensor space attached to $V$, i.e. $\mathcal{T}(V)=\oplus_{s,t\in\Z_{\geq0}}V^{\otimes s}\otimes (V^{\vee})^{\otimes t}$. Let $\{s_{\alpha}\}_{\alpha\in \mathcal{J}}$ be the set of tensors on $W_{\Q}$ fixed under the extended action of $G_{\Q}\hra \mathrm{GL}(W)$ on $\mathcal{T}(W)$. 

Let $(\pi:A\rightarrow \Sh_{K},\lambda_A)$ be the pull back to $\Sh_{K}=\Sh_{K}(G,X)$ of the universal abelian scheme with a polarization of degree $d$ over $\mathcal{A}_{g,d,H}$, where as usual we assume $K^p$ to be sufficiently small for this to make sense. We have the local system $\mathcal{W}:=R^1(\pi^{an})_{\ast}\Q$ of $\Q$-vector spaces and the (analytic) vector bundle $\mathbb{W}^{an}_{\C}:=R^1(\pi^{an})_{\ast}\Omega^{\bullet}_{A_{\C}/\Sh_{K}}$ with Gauss-Manin connection (which is an integrable connection with regular singularities); by Deligne's theorem, the latter is the analytification of a unique algebraic vector bundle $\mathbb{W}:=R^1\pi_{\ast}\Omega^{\bullet}_{A/\Sh_{K}}$ with integrable connection over $\Sh_{K}$. 
Each tensor $s_{\alpha}$ defines a global section $s_{\alpha,\Betti}$ of the local system $\mathcal{T}(\mathcal{W})$ and also a global section $s_{\alpha,\dR}$ of the vector bundle $\mathcal{T}(\mathbb{W})$ (\cite[(2.2)]{Kisin10}). 

Now, let $F\supset E(G,X)$ be an extension which can be embedded in $\C$; we fix an embedding $\sigma_{\infty}:\overline{F}\hra\C$ which extends the given embedding $E(G,X)\hra\C$. 
For $x\in \Sh_K(G,X)(F)$, let $A_x$ be the corresponding abelian variety over $F$, and let $H^m_{\dR}(A_x/F)$ and $H^m_{\et}(A_x\otimes_{F}{\overline{F}},\Qp)$ denote respectively the de Rham cohomology of $A_x/F$ and the \'etale cohomology of $A_x\otimes_{F}\overline{F}$.
For each tensor $s_{\alpha}\in \mathcal{T}(W)$, let $s_{\alpha,\Betti,\sigma_{\infty}(x)}\in \mathcal{T}(H^1_{\Betti}(\sigma_{\infty}(A_x),\Q))$ denote the fibre of $s_{\alpha,\Betti}$ at $\sigma_{\infty}(x)\in \Sh_K(G,X)(\C)$; by the moduli interpretation of the complex points of $\Sh_{K}(G,X)$ (\cite[Prop.3.9]{Milne94}), there is an isomorphism $W^{\vee}\isom H^1_{\Betti}(A_{\sigma_{\infty}(x)},\Q)$, uniquely determined up to action of $G(\Q)$, where $G$ acts on $W^{\vee}$ by the contragredient representation, and $s_{\alpha,\Betti,\sigma_{\infty}(x)}$ is the image of $s_{\alpha}$ by any such isomorphism (so, does not depend on the choice of such isomorphism). For each prime $l$, let $s_{\alpha,l,\sigma_{\infty}(x)}$ be the tensor in $\mathcal{T}(H^1_{\et}(A_x\otimes_{F}{\overline{F}},\Ql))$ which is the preimage of $s_{\alpha,\Betti,\sigma_{\infty}(x)}$ under the canonical isomorphism
\begin{equation} \label{eq:l-adic=Betti}
H^m_{\et}(A_x\otimes_{F}{\overline{F}},\Ql)\isom H^m_{\et}(\sigma_{\infty}(A_x\otimes_{F}{\overline{F}}),\Ql)=H^m_{\Betti}(\sigma_{\infty}(A_x),\Q)\otimes\Ql,
\end{equation}
where the first isomorphism $\sigma_{\infty}$ is the proper base change isomorphism in \'etale cohomology.  
Also, it is known  (\cite[Cor.2.8]{DMOS82}) that the image of $s_{\alpha,\Betti,\sigma_{\infty}(x)}$ under the comparison isomorphism
\begin{equation*}H^m_{\Betti}(\sigma_{\infty}(A_x),\Q)\otimes\C= H^m_{\dR}(A_x\otimes_{F,\sigma_{\infty}}\C)\end{equation*} lies in $\mathcal{T}(H^1_{\dR}(A_x))\otimes_{F,\sigma_{\infty}}\sigma_{\infty}(\overline{F})$, i.e. equals $\sigma_{\infty}(s_{\alpha,\dR,\sigma_{\infty}(x)})$, for some $s_{\alpha,\dR,\sigma_{\infty}(x)}\in \mathcal{T}(H^1_{\dR}(A_x\otimes_F\overline{F}))$.
The element $(s_{\alpha,\dR,\sigma_{\infty}(x)},s_{\alpha,\et,\sigma_{\infty}(x)}:=(s_{\alpha,l,\sigma_{\infty}(x)})_l)$ of the tensor space $\mathcal{T}(H^1_{\A}(A_x),\Q(1))$ is an absolute Hodge cycle on $A_x$, where $H^1_{\A}(A_x):=H^1_{\dR}(A_x)\times \prod H^1_{\et}(A_x,\Ql)$ (cf. \cite{DMOS82}), from which it follows (\cite[(2.2.1)]{Kisin10}) that $s_{\alpha,\et,\sigma_{\infty}(x)}$ is fixed under the canonical Galois action of $\Gal(\overline{F}/F)$, hence in fact, $s_{\alpha,\dR,\sigma_{\infty}(x)}$ belongs to the subspace $\mathcal{T}(H^1_{\dR}(A_x/F))$.

\subsubsection{} 
We show how every geometric point $z\in\sS_{K}\otimes\Fpb(k)$ ($k=\overline{k}$, $\mathrm{char} k=p$) gives rise to an $F$-isocrystal over $k$ with $G_{\Q_p}$-structure; here, we assume that $L(k)$ can embed in $\C$, for example, it suffices if $k$ has finite transcendence degree over $\Fpb$. The $p$-divisible group $(A_z[p^{\infty}],\lambda_{A_z}[p^{\infty}])$ with quasi-polarization gives rise to an $F$-isocrystal over $k$ equipped with a non-degenerate alternating pairing. 
More precisely, the crystalline cohomology $M=H^1_{\cris}(A_{z}/W(k))$ of $A_{z}$ with quasi-polarization $\lambda_{z}$ is a quasi-polarized $F$-crystal $(M,\phi,\langle\ ,\ \rangle)$: $M$ is a free $W(k)$-module of rank equal to $2\dim A_z$, $\phi\in\mathrm{End}_{\Z_p}(M)$ is a $\sigma$-linear endomorphism such that $pM\subset\phi(M)$, and $\langle\ ,\ \rangle:M\times M\rightarrow W(k)$ is an alternating form with the property that $\langle\phi(v_1),\phi(v_2)\rangle=p\langle v_1,v_2\rangle^{\sigma}$ for $v_1,v_2\in M$.

To show that the $F$-isocrystal $(M\otimes L(k),\phi,\psi)$ has $G_{\Q_p}$-structure, let $x:\mathrm{Spec}(W(k))\rightarrow \sS_{K}$ be a lift of $z$, which exists since $\sS_{K}(G,X)$ is smooth over $\cO_{(\wp)}$ and $\cO_{(\wp)}$ is unramified over $\Z_{(p)}$ (\cite[Cor.4.7]{Milne94}), and set $(A_{x},\lambda_{x}):=x^{\ast}(A,\lambda_A)$; here we identified $E_{\wp}=L(\F_q)$ for the residue field $\F_q$ of $\cO_{(\wp)}$. 
Let $\{t_{\beta,\dR,x},t_{\beta,\et,x}\}_{\beta\in\mathcal{J}_x}\in \mathcal{T}(H^1_{\A}(A_x),\Q(1))$ be the set of absolute Hodge cycles on $A_x$ (see  \cite{DMOS82} for the notations appearing here), and for each $\beta\in\mathcal{J}_x$, let $t_{\beta,z}\in \mathcal{T}(M\otimes_{W(k)}L(k))$ denote the image of $t_{\beta,\dR,x}$ under the canonical isomorphism $H^1_{\dR}(A_{x}/L(k))=H^1_{\cris}(A_{z}/W(k))\otimes L(k)$. Then, one can identify $\mathcal{J}$ as a subset of $\mathcal{J}_x$ (see below for a proof) and the triple 
\begin{equation*}(M\otimes F,\phi,(t_{\beta,z})_{\beta\in\mathcal{J}}) \end{equation*}
depends only on $z$, not on the lift $x$%
%
: this is basically because for any two lifts $x,x'$ of $z$, the composite map $H^1_{\dR}(A_{x}/L(k))=H^1_{\cris}(A_{z}/W(k))\otimes L(k)=H^1_{\dR}(A_{x'}/L(k))$ is given by parallel transport with respect to the Gauss-Manin connection (\cite[(2.9)]{BO83}). 
Also, the cycle $t_{\beta,z}$ is \textit{crystalline}, namely $t_{\beta,z}$ belongs to the $F^{0}$-filtration of $\mathcal{T}(M\otimes L(k))$ and is fixed under the Frobenius $\phi$ (originally, this was proved when $x$ is defined over a number field, by Wintenberger and Blasius \cite[$\S$5]{Blasius94}, as a consequence of the fact that the Hodge cycles on an abelian variety over a number field are \textit{de Rham} and the compatibility of the crystalline and the de Rham comparison isomorphisms \cite[5.1(5)]{Blasius94}. Then, Vasiu \cite[Sec.8]{Vasiu08} generalized it to arbitrary fields). 
Moreover, there exists an $L(k)$-isomorphism
\begin{equation} \label{eq:comparison_over_L(k)}
W^{\vee}\otimes_{\Q} L(k)\isom M\otimes_{W(k)}L(k)
\end{equation}
which maps $s_{\alpha}$, for every $\alpha\in\mathcal{J}$, to $t_{\beta,z}$ for some $\beta=\beta(\alpha)\in\mathcal{J}_x$. 
Indeed, as before we choose an embedding $\sigma_{\infty}:\overline{L(k)}\hra\C$ which extends the given embedding $E(G,X)\hra\C$. Then, for each $\beta\in\mathcal{J}_x$, 
as $t_{\beta,\dR,x}$ is the de Rham component of an absolute Hodge cycle on $A_x$, 
there exists a tensor $v_{\beta}\in \mathcal{T}(W)$ such that $\sigma_{\infty}(t_{\beta,\dR,x})$ is the image of $v_{\beta}$ under the comparison isomorphism 
\begin{equation} \label{eq:W_dual=Betti=deRham}
W^{\vee}\otimes\C\isom H^1_{\Betti}(\sigma_{\infty}(A_x),\Q)\otimes\C\isom H^1_{\dR}(A_x/L(k))\otimes_{L(k),\sigma_{\infty}}\C,
\end{equation}
where as noted above the first isomorphism exists and is unique up to the action of $G(\Q)$ on $W^{\vee}$.
Because for each $\alpha\in\mathcal{J}$, $\sigma_{\infty}(s_{\alpha,\dR,\sigma_{\infty}(x)})$ is the image of $s_{\alpha}$ of the same isomorphism and also is the de Rham component of an absolute Hodge cycle on $A_x$, for every $\alpha\in\mathcal{J}$, $s_{\alpha}$ equals $v_{\beta}$ for some $\beta=\beta(\alpha)\in\mathcal{J}_x$ (which gives an injection $\mathcal{J}\hra \mathcal{J}_x$ mentioned above), i.e. 
the set $\{v_{\beta}\}_{\beta\in\mathcal{J}_x}$ contains the set $\{s_{\alpha}\}_{\alpha\in\mathcal{J}}$. Then, the functor defined on $L(k)$-schemes 
\begin{equation*}R\mapsto \Hom_{R}((W^{\vee}_{L(k)}\otimes R,\{v_{\beta}\otimes1\}),(H^1_{\dR}(A_{x}/L(k))\otimes R,\{t_{\beta,\dR,x}\otimes1\}))\end{equation*} 
of isomorphisms between $W^{\vee}\otimes_{\Q}L(k)$ and $H^1_{\dR}(A_{x}/L(k))$ taking $v_{\beta}$ to $t_{\beta,\dR,x}$ for every $\beta\in\mathcal{J}(\subset \mathcal{J}_x)$ is represented by a scheme which is non-empty as it has a $\C$-valued point, and thus is a torsor over $L(k)$ under $G_{L(k)}$. Since $H^1(L(k),G_{L(k)})=\{0\}$ (Steinberg's theorem), this torsor is trivial, namely has an $L(k)$-valued point.

Therefore, when one chooses an isomorphism as in (\ref{eq:comparison_over_L(k)}) to transport the $\sigma$-linear map $\phi$ to $W^{\vee}$, there exists an element $b\in G(L(k))$ with \begin{equation*}\phi=\rho_{W^{\vee}}(b)(\mathrm{id}_{W^{\vee}}\otimes\sigma),\end{equation*}
where $\rho_{W^{\vee}}:G\hra \mathrm{GL}(W^{\vee})$ denotes the contragredient representation: indeed, $\rho_{W^{\vee}}(b)(\mathrm{id}_{W^{\vee}}\otimes\sigma)$ fixes each $v_{\beta}\otimes1\in \mathcal{T}(W^{\vee}\otimes L(k))$. Although $b$ depends on the choice of an isomorphism (\ref{eq:comparison_over_L(k)}), its $\sigma$-conjugacy class $\overline{b}\in B(G_{\Q_p})$ is independent of such choice. 
This shows that the $F$-isocrystal over $k$ attached to $z$ is an $F$-isocrystal with $G_{\Q_p}$-structure.
For arbitrary point $z$ of $\sS_K(G,X)\otimes\kappa(\wp)$, we define the $F$-isocrystal attached to $z$ to be the $F$-isocrystal of the geometric point in $\sS_K(G,X)(k)$ induced from $z$ for \textit{any} algebraically closed field $k$ containing the residue field $\kappa(z)$ of $z$: the resulting $F$-isocrystal does not depend on the choice of $k'$ (cf. \cite[$\S$8]{VW13}).
Therefore, we obtain a (set-theoretic) map
\begin{equation*}\Theta:\sS_K(G,X)\otimes\kappa(\wp)\rightarrow B(G).\end{equation*}

We note that the same argument also establishes that if $x\in \Sh_K(G,X)(F)$ is a special point corresponding to a special Shimura datum $(T,h\in\Hom(\dS, T_{\R})\cap X)$, the associated isocrystal has a representative in $T(L)$. Indeed, by the moduli interpretation, there is an isomorphism of polarized $\Q$-Hodge structures $W^{\vee}\isom H_1(\sigma_{\infty}(A_x),\Q)$, where on the left side the Hodge structure is given by $h$, so $T\subset\mathrm{GL}(W)$ is the stabilizer of the Hodge tensors $\{v_{\beta}\}_{\beta\in\mathcal{J}_x}\subset \mathcal{T}(W)$, and according to the above discussion, the associated element $b\in G(L(k))$ fixes all of the tensors $\{v_{\beta}\}_{\beta\in\mathcal{J}_x}$.

Finally, we remark that although our definition of the $F$-isocrystal attached to a point of $\sS_K(G,X)(k)$ uses \emph{cohomology} spaces, one can equally work with homology spaces, as adopted by other people (such as \cite{VW13}). This does not alter the definition of the map $\Theta$.

\subsubsection{} We recall that for a torus $T$ over a non-archimedean local field $F$, 
there exists a unique maximal compact subgroup of $T(F)$ (which is also open). In fact, $T$ has a canonical integral form $\mathcal{T}$ over $\cO_F$ such that for every finite extension $F'$ of $F$, the maximal compact subgroup of $T(F')$ is $\mathcal{T}(\cO_{F'})$. If $E$ is a finite Galois extension of $F$ splitting $T$, the maximal compact subgroup of $T(E)=\Hom(X^{\ast}(T),F^{\times})(\simeq (E^{\times})^{\dim T})$ is $\mathcal{T}(\cO_E)=\Hom(X^{\ast}(T),\cO_E^{\times})(\simeq (\cO_E^{\times})^{\dim T})$; furthermore we have $\mathcal{T}(\cO_F)=\Hom_{\Gal(E/F)}(X^{\ast}(T),\cO_E^{\times})$. If $T$ is unramified over $F$, $\mathcal{T}(\cO_F)$ is the unique hyperspecial subgroup of $T(F)$ 
(cf.  Thm. 1 and the discussion after Prop. 2 in $\S$10.3 of \cite{Voskresenskii98}, \cite{Tits79}).  We will also write simply $T(\cO_F)$ for $\mathcal{T}(\cO_F)$.

\begin{lem} \label{lem:isocrystal_of_special_pt}
Choose an embedding $\iota_p:\Qb\hra \Qpb$ inducing given prime $\wp$ of $E\subset\Qb$. Let $K_p$ be a hyperspecial subgroup of $G(\Qp)$ and $K^p$ a compact open subgroup of $G(\A_f^p)$. Let $x=[h,k_f]\in \Sh_K(G,X)(\Qb)$ be a special point, where $h\in \Hom(\dS, T_{\R})\cap X$ for a maximal $\Q$-torus $T$ of $G$ and $k_f\in G(\A_f)$. Suppose that $T$ is unramified over $\Qp$ and the (unique) hyperspecial subgroup $T(\Z_p)$ of $T(\Qp)$ is contained in $K_p$. Then, 

(1) $x\in \Sh_K(G,X)(\Qpb)$ is defined over an unramified extension $F$ of $E(G,X)_{\wp}$, thus extends over $\mathcal{O}_F$, and

(2) the $F$-isocrystal with $G_{\Qp}$-structure of its reduction $z\in\Sh_K(G,X)(\Fpb)$ is represented by $\mu_h^{-1}(p)\in T(\Qp^{\nr})$. In particular, the associated Newton point $\nu_{T_{\Qp}}(\mu_h(p))\in \mathcal{N}(T_{\Qp})=X_{\ast}(T)^{\Gamma(p)}_{\Q}$ is
\begin{equation*}-\frac{1}{r}\Nm_{L(\F_q)/\Qp}(\mu_h),\end{equation*}
where $L(\F_q)=\mathrm{Frac}(W(\F_q))$ is any splitting field of $T_{\Qp}$ ($q=p^r$).
\end{lem}

A special point of $\Sh_K$ will be said to have \textit{good reduction at $\iota_p$} if it
extends to an $\cO_{\iota_p}$-valued point of $\sS_K$ for the valuation ring $\cO_{\iota_p}$ of $\Qb$ defined by $\iota_p$. 

\textsc{Proof.} (1) The statement on the field of definition is an easy consequence of the reciprocity law chararcterizing the canonical model of Shimura varieties (\cite[$\S$3]{Deligne71}) and the local class field theory.
For any compact open subgroup $K'$ contained in $T(\A_f)\cap K$, there is a natural map $\Sh_{K'}(T,\{h\})\rightarrow \Sh_K(G,X)\otimes_{E(G,X)}E(T,\{h\})$ of weakly canonical models of Shimura varieties (\cite[(3.13)]{Deligne71}).
Hence, it suffices to show that for the compact open subgroup $K'=T(\Z_p)\times K'^p$ for any $K'^p\subset T(\A^p)\cap K^p$, the connected components of the finite scheme $\Sh_{K'}(T,\{h\})$ over $E(T,\{h\})$ are defined over an (abelian) extension $E'$ of $E:=E(T,\{h\})$ such that the prime $\mathfrak p'$ of $E'$ induced by $\iota_p$ is unramified over $\wp$ (\cite[3.15]{Deligne71}). The action of $\Gal(\Qb/E)$ (which factors through the abelianization $\Gal(E^{\text{ab}}/E)$) on the group $\Sh_{K'}(T,\{h\})(\Qb)=T(\Q)\backslash T(\A_f)/K'$ is given by (the adelic points of) the reciprocity map (\cite[3.9]{Deligne71})
\begin{equation} \label{eq:reciprocity_map}
r_{T,\mu_h}:\mathrm{Res}_{E/\Q}(\mathbb{G}_{mE})\stackrel{\mathrm{Res}_{E/\Q}(\mu_h)}{\longrightarrow}\mathrm{Res}_{E/\Q}(T_E)\stackrel{N_{E/\Q}}{\longrightarrow} T
\end{equation}
via class field theory: for $\rho\in \Gal(E_{\wp}^{\text{ab}}/E_{\wp})=E_{\wp}^{\times}$ and $[t]_{K'}\in  \Sh_{K'}(T,\{h\})(\overline{E_{\wp}})=T(\Q)\backslash T(\A_f)/K'$, one has 
\begin{equation*}\rho[t]=[t\cdot r_{T,\mu}(\rho)^{-1}],\end{equation*}
where one uses the convention on the class field theory that under the identification $\Gal(E_{\wp}^{\text{ab}}/E_{\wp})=E_{\wp}^{\times}$, 
the \textit{geometric Frobenius} corresponds to a uniformizer of $E_{\wp}$.
But, the image in $T(\Qp)$ of the unit group $(\mathcal{O}_E\otimes\Z_p)^{\times}$ is obviously contained in the maximal compact open subgroup $T(\Z_p)$. Since $\Gal(E_{\wp}^{\mathrm{ab}}/E_{\wp}^{\nr})=(\mathcal{O}_E)_{\wp}^{\times}\subset (\mathcal{O}_E\otimes\Z_p)^{\times}$, this proves the claim.  
Then, by the Neron extension property (Thm. \ref{thm:integral_canonical_model}, (2) and Remark \ref{rem:Neron_extension_property}) of the integral canonical model, $x$ extends to an $\mathcal{O}_F$-valued point of $\sS_K$.

(2) Set $L:=L(\Fpb)$; by (1), we may regard $L$ as an extension of $F$. 
The statement on the Newton point follows from the first statement and Theorem \ref{thm;RR96-Thm.1.15}. The claim that the $F$-isocrystal of $z$  is represented by $\mu_h^{-1}(p)\in T(L)$ is Lemma 13.1 of \cite{Kottwitz92}; as was observed above, we already know that it is represented by an element of $T(L)$. More precisely, this lemma is a corollary of Lemma 12.1 of loc.cit., which is a statement purely for a torus over a non-archimedean local field and a cocharacter of it (in particular, in Lemma 12.1 it is not required that the torus is defined by a PEL-datum), and
the argument in \cite{Kottwitz92} of deducing Lemma 13.1 from Lemma 12.1 still works in our situation of a general Hodge-type. To see that in detail, we reproduce the Kottwitz's argument here, pointing out necessary modifications explicitly along the way. 

We fix an embedding $\sigma_{\infty}:\overline{L}\hra\C$ which extends the given embedding $E(G,X)\hra\C$. The moduli interpretation of $\Sh_K(G,X)(\C)$ provides an isomorphism $W^{\vee}\isom H^1_{\Betti}(\sigma_{\infty}(A_x),\Q)$, determined uniquely up to action of $G$; we fix one. Accordingly, we get an isomorphism $W^{\vee}_{\Qp}\isom H^1_{\et}((A_x)_{\Qpb},\Qp)$ by the isomorphism (\ref{eq:l-adic=Betti}), via which we transport the action of $T_{\Qp}$ on $W^{\vee}_{\Qp}$ (i.e. the contragredient representation) to $H^1_{\et}((A_x)_{\Qpb},\Qp)$. Then, the Galois representation $\Gal(\Qpb/F)\ra \mathrm{GL}(H_1((A_x)_{\Qpb},\Qp))$ factors through $T_{\Qp}$ and the restriction to $\Gal(\Qpb/\Qp^{\nr})=\Gal(\Qpb/F^{\nr})$ of the resulting representation equals the homomorphism
\begin{equation} \label{eq:rho_T_mu}
\rho_{T_{\Qp},\mu_h}:\Gal(\Qpb/\Qp^{\nr})\ra T_{\Qp}
\end{equation}
attached to $(T_{\Qp},\mu_h)$ by the procedure of $\S$12 of \cite{Kottwitz92} (just as in the PEL-type case); in particular, it is a crystalline representation. Indeed, as we have seen from \cite[(2.2.1)]{Kisin10}, the absolute Hodge cycles on $A_x$ are invariant under $\Gal(\overline{F}/F)$, which proves the first claim. The second claim holds because we can reduce our case to the Siegel case by means of the symplectic embedding $G\hra \mathrm{GSp}(W,\psi)$.

Next, we need three $\Qp$-linear Tannakian categories: 
$\mathscr{REP}(T_{\Qp})$ is the neutral $\Qp$-linear Tannakian category of finite-dimensional representations of the torus $T_{\Qp}$, $\mathscr{CR}(F)$ is the neutral $\Qp$-linear Tannakian category of finite-dimensional crystalline representations of $\Gal(\Qpb/F^{\nr})$, and $\mathscr{CRYS}$ denotes the $\Qp$-linear Tannakian category of $F$-isocrystals $(U,\Phi)$, where $U$ is a finite dimensional vector space over $L$ and $\Phi$ is a $\sigma$-linear bijection $\Phi:U\ra U$.
Then, to any cocharacter $\mu\in X_{\ast}(T_{\Qp})$, there is attached a (composite) tensor functor
\begin{equation} \label{eq:defn_of_Xi}
\Xi_{\mu}:\mathscr{REP}(T_{\Qp})\stackrel{\rho_{\mu_h}^{\ast}}{\ra} \mathscr{CR}(F)\stackrel{\DD}{\ra} \mathscr{CRYS}.
\end{equation}
Here, the first functor $\rho_{\mu_h}^{\ast}$ is dual to the homomorphism $\rho_{\mu_h}:=\rho_{T_{\Qp},\mu_h}:\Gal(\Qpb/\Qp^{\nr})\ra T_{\Qp}$ (\ref{eq:rho_T_mu})  and the second functor $\DD$ is the Dieudonn\'e functor in the Fontaine-Messing theory: for a finite unramified extension $E(\subset \Qp^{\nr})$ of $F$ and 
a crystalline representation $V$ of $\Gal(\Qpb/E)$, 
\begin{equation*}\DD(V)=(B_{\cris}\otimes_{\Qp}V)^{\Gal(\Qpb/E)},\end{equation*}
where $B_{\cris}$ is the crystalline period ring of Fontaine-Messing (\cite{FM87}). 
On the other hand, if for a field $K$, $\mathscr{VEC}_K$ denotes the category of finite-dimensional $K$-vector spaces, from the fibre functor $\mathscr{CRYS}\ra \mathscr{VEC}_{L}: (U,\Phi)\mapsto U$, we get a non-standard fibre functor $\omega_1:\mathscr{REP}(T_{\Qp})\ra \mathscr{VEC}_{L}$.
By Steinberg's theorem ($H^1(L,T_{L})=\{0\}$), this fibre functor is isomorphic to the standard fibre functor $\omega_0:\mathscr{REP}(T_{\Qp})\ra \mathscr{VEC}_{L}:(\rho,V)\mapsto V\otimes_{\Qp} L$.
We fix an isormophsm $f:\omega_1\isom\omega_0$ of fibre functors.

Now, it follows from the above discussion (on the homomorphism $\rho_{T_{\Qp},\mu_h}$) that for the natural representation $W_{\Qp}$ of $T_{\Qp}$, the associated crystalline representation $\rho_{u_h}^{\ast}(W_{\Qp})$ is the natural representation of $\Gal(\Qpb/\Qp^{\nr})$ on $H_1((A_x)_{\Qpb},\Qp)$, hence $\Xi_{\mu_h}(W_{\Qp})$ is canonically isomorphic to the \emph{dual} of the $F$-isocrystal $H^1_{\cris}(A_z/L)=H^1_{\cris}(A_z/W(\Fpb))\otimes L$  (cf. proof of \cite[Lemma13.1]{Kottwitz92}); here, by definition, the $F$-isocrystal dual to $H^1_{\cris}(A_z/L)$ is the linear dual $H_1^{\cris}:=\Hom_L(H^1_{\cris}(A_z/L),L)$ equipped with the Frobenius operator $\Phi_1(f)(v):=p^{-1}\cdot{}^{\sigma}f(Vv)$ for $f\in H_1^{\cris}$ and $v\in H^1_{\cris}(A_z/L)$ ($V$ being the Verschiebung). 
Then, Lemma 12.1 of \cite{Kottwitz92} says that if one transports the Frobenius operator $\Phi_1$ on $\Xi_{\mu_h}(W_{\Qp})=H_1^{\cris}$ to $W\otimes L$
via the isomorphism
\begin{equation*}f(W_{\Qp}):\omega_1(W_{\Qp})=H_1^{\cris}\isom \omega_0(W_{\Qp})=W_{\Qp}\otimes L,\end{equation*} 
and writes $(H_1^{\cris},\Phi_1)\simeq(W\otimes L,\rho_{W}(b)(\mathrm{id}_{W}\otimes\sigma))$, then one has an equality \begin{equation*}\overline{b}=\overline{\mu_h^{-1}(p)}\end{equation*}
in $B(T_{\Qp})$. It also holds that 
\begin{equation*}H^1_{\cris}(A_x/L)\simeq(W^{\vee}\otimes L,\rho_{W^{\vee}}(b)(\mathrm{id}_{W^{\vee}}\otimes\sigma)).\end{equation*} 
Therefore, to prove the statement (2), we only need to show that for any identification $f:\omega_1\isom\omega_0$ of fibre functors, $f(W_{\Qp}):H_1^{\cris}\isom W_L$ induces an isomorphism as in (\ref{eq:comparison_over_L(k)}), namely $f(W_{\Qp})^{\vee}: W^{\vee}_L\isom H^1_{\cris}(A_z/L)$ sends $s_{\alpha}$, for every $\alpha\in\mathcal{J}$, to $t_{\beta,z}$ for some $\beta=\beta(\alpha)\in\mathcal{J}_x$.  We remark that in the PEL-type setting, the tensors involved are morphisms in tensor categories, so this was obvious. In a general Hodge-type case, the proof of this fact requires more sophisticated argument and also a non-trivial fact (\cite{Blasius94}) that the Hodge cycles on abelian varieties over number fields are de Rham.

In more detail, we consider each tensor $s_{\alpha}$ as a morphism $s_{\alpha}:\mathbf{1}\ra W_{\Qp}^{\otimes}$ in $\mathscr{REP}(T_{\Qp})$, where $\mathbf{1}$ is the trivial representation $\Qp$ (which is an identity object of the tensor category $\mathscr{REP}(T_{\Qp})$) and $W_{\Qp}^{\otimes}$ denotes some specific object of $\mathscr{REP}(T_{\Qp})$. Apply the functor $\Xi_{\mu}$, we get a morphism $t_{\alpha}':\Xi_{\mu}(\mathbf{1})\ra \Xi_{\mu}(W_{\Qp}^{\otimes})$ in $\mathscr{CRYS}$. But, as $\Xi_{\mu}(\mathbf{1})$ is the trivial $F$-isocrystal $(L,\sigma)$, again this is equivalently regarded as a crystalline tensor on $\Xi_{\mu}(W_{\Qp})=H_1^{\cris}$. Moreover, as $f:\omega_1\ra \omega_0$ is an isomorphism of fibre functors, we have a commutative diagram of $L$-vector spaces
\begin{equation*}\xymatrix{ \omega_1(\mathbf{1}) \ar[r]^{t_{\alpha}'} \ar[d]_{f(\mathbf{1})}^{\simeq} & \omega_1(W_{\Qp}^{\otimes})\ar[d]_{\simeq}^{f(W_{\Qp}^{\otimes})} \\
\omega_0(\mathbf{1}) \ar[r]^{s_{\alpha}} & \omega_0(W_{\Qp}^{\otimes}). }\end{equation*} 
As the isomorphism $f(\mathbf{1})$ is the multiplication by an element of $L$, by modifying $f$, we may assume that $f(\mathbf{1})$ is the identity on $\mathbf{1}=\Qp$. This implies that under $f(W_{\Qp})^{-1}$, the tensor $s_{\alpha}$ on $W_{\Qp}=\omega_0(W_{\Qp})$ maps to a crystalline tensor $t_{\alpha}'$ on $H_1^{\cris}=\omega_1(W_{\Qp})$. So, it remains to check that the tensor $t_{\alpha}'$ is the de Rham component of an absolute Hodge cycle on $A_x$. But, the image of the tensor $s_{\alpha}:\mathbf{1}\ra W_{\Qp}^{\otimes}$ under $\rho_{\mu_h}^{\ast}$ is just the image of it under our chosen identification $W^{\vee}_{\Qp}\isom H^1_{\et}((A_x)_{\Qpb},\Qp)$, hence equals the \'etale component $s_{\alpha,p,\sigma_{\infty}(x)}:\mathbf{1}\ra H_1^{\et}((A_x)_{\Qpb},\Qp)^{\otimes}$ of a Hodge cycle $s_{\alpha,\Betti,\sigma_{\infty}(x)}$ on $\sigma_{\infty}(A_x)$ by the moduli interpretation of $\Sh_K(G,X)(\C)$. In turn, the image of $s_{\alpha,p,\sigma_{\infty}(x)}$ under the second functor $\DD$ is $s_{\alpha,\dR,\sigma_{\infty}(x)}:\mathbf{1}\ra H_1^{\dR}(A_x/L)^{\otimes}:=\Hom_L(H^1_{\dR}(A_x/L),L)^{\otimes}$, because the etale and the de Rham component of an absolute Hodge cycle on an abelian variety over a number field match under the $p$-adic comparison isomorphism from which the functor $\DD$ is induced (\cite[Thm.0.3]{Blasius94}). So, we have just proved that $t_{\alpha}'=\Xi_{\mu_h}(s_{\alpha})=s_{\alpha,\dR,\sigma_{\infty}(x)}$, and consequenty the statement (2).
$\qed$

\subsubsection{} \label{overline{mu}(G,X)} 
We fix an embedding $\iota_p:\Qb\hookrightarrow\Qpb$ for an algebraic closure $\Qpb$ of $\Q_p$, thereby an inclusion 
\begin{equation*}\Gamma(p):=\Gal(\Qpb/\Q_p)\hookrightarrow
\Gamma:=\Gal(\Qb/\Q).\end{equation*}
This allows us to consider $c(G,X)$ as an element of
$\mathcal{C}_{G}(\Qpb)$. As such, the field of
definition of $c(G,X)$ is equal to the $\wp$-adic completion
$E(G,X)_{\wp}$ of $E(G,X)$, where $\wp$ is the prime of $E(G,X)$
induced by $\iota_p$. 

Let $\mathcal{BR}(G_{\Q_p})=(X^{\ast},R^{\ast},X_{\ast},R_{\ast},\Delta)$ be the based root datum of $G_{\Q_p}$ and $\overline{C}\subset (X_{\ast})_{\Q}$ the closed Weyl chamber associated with the root base $\Delta$. Denoting also by $c(G,X)$ its (unique) representative in $\overline{C}$, we set 
\begin{align*} 
\overline{\mu}(G,X)&:=
|\Gamma(p)\cdot c(G,X)|^{-1}
\sum_{\mu'\in\Gamma(p)\cdot c(G,X)}\mu'\qquad\in\overline{C}\\
&=[\Gamma(p):\Gamma(p)_c]^{-1} \sum_{\gamma\in
\Gamma(p)/\Gamma(p)_c}\gamma(c(G,X)),
\end{align*}
Here, the action of $\Gamma(p)$ on $\overline{C}$ is the canonical action and $\Gamma(p)_c:=\Gal(\Qpb/E(G,X)_{\wp})$, i.e. the stabilizer of $c(G,X)$ in $\Gamma(p)$.

\subsubsection{} \label{subsubsec:mu_natural}
As $X_{\ast}(T)=X^{\ast}(\widehat{T})$ for the dual torus $\widehat{T}$ of $T$, regarded as a character on $\widehat{T}$, we can restrict $c(G,X)\in\overline{C}(T,B)$ to the subgroup $Z(\widehat{G})^{\Gamma(p)}$ of $\widehat{T}$, obtaining an element 
\begin{equation*}\mu^{\natural}\in X^{\ast}(Z(\widehat{G})^{\Gamma(p)}).\end{equation*}
This depends only on the pair $(G,X)$ (more precisely, on the pair $(G_{\Q_p},c(G,X)\in \mathcal{C}_G(\Qpb))$). 
Alternatively, under the identification $X^{\ast}(Z(\widehat{G})^{\Gamma(p)})=\pi_1(G)_{\Gamma(p)}$, the element $\mu^{\natural}$ equals the image of $\mu\in X_{\ast}(T)$ under the canonical map $X_{\ast}(T)\rightarrow \pi_1(G)_{\Gamma(p)}=(X_{\ast}(T)/\Sigma_{\alpha\in\Phi(G,T)}\Z\alpha^{\vee})_{\Gamma(p)}$, where $\mu$ is any cocharacter of $T_{\Qpb}$ that lies in the conjugacy class $c(G,X)$ of cocharacters of $G_{\Qpb}$.

\subsubsection{} \label{subsubsec:B(G,X)} We define a finite subset $B(G_{\Q_p},X)$ of $B(G_{\Q_p})$ (following Kottwitz \cite[Sec.6]{Kottwitz97} and Rapoport \cite[Sec.4]{Rapoport05}): 
\begin{equation*}B(G_{\Q_p},X):=\left\{\ \overline{b}\in B(G_{\Q_p})\ |\quad \kappa_{G_{\Q_p}}(\overline{b})=\mu^{\natural},\quad \overline{\nu}_{G_{\Q_p}}(\overline{b})\preceq \overline{\mu}(G,X)\ \right\},\end{equation*}
where $\kappa_{G_{\Q_p}}$ is as defined in Theorem \ref{thm;RR96-Thm.1.15} and $\preceq$ is the natural partial order on the closed Weyl chamber $\overline{C}(T,B)$ defined by that $\nu\preceq \nu'$ if $\nu'-\nu$ is a nonnegative linear combination with \emph{real} coefficients of simple coroots in $R_{\ast}(T)$ (\cite{RR96}, Lemma 2.2). 

One knows (\cite[4.13]{Kottwitz97}) that the map 
$(\overline{\nu},\kappa):B(G_{\Q_p})\rightarrow \mathcal{N}(G_{\Q_p})\times X^{\ast}(\widehat{Z}(G)^{\Gamma(p)})$ is injective, hence $B(G_{\Q_p},X)$ can be identified with a subset of $\mathcal{N}(G_{\Q_p})$.

\begin{prop} 
Suppose $G$ is unramified over $\Q_p$ and $K=K_pK^p$ is hyperspecial at $p$. 

(1) The image of $\Theta:\sS_K(G,X)(\kappa(\wp))\rightarrow B(G)$ is contained in $B(G_{\Q_p},X)$.

(2) For $\overline{b}\in B(G_{\Q_p})$, the subset 
\begin{equation*}\mathcal{S}_{\overline{b}}:=\Theta^{-1}(\overline{b})\end{equation*}
is locally closed inside $\sS_{K}(G,X)\otimes\Fpb$.
\end{prop}

For (1), see e.g. \cite{RR96}. The second statement (2) is due to \cite[5.3.1]{Vasiu11}.

Endowed with reduced induced subscheme structure, the subvarieties $\mathcal{S}_{\overline{b}}$ of $\sS_{K}(G,X)\otimes\Fpb$ are called the \emph{Newton strata} of $\sS_K(G,X)$.

In \cite{Rapoport05}, Conjecture 7.1, Rapoport conjectures 
\begin{conj} \label{conj:non-emptiness_of_Newton_strata}
\begin{equation*}\mathrm{Im}(\Theta)=B(G_{\Q_p},X).\end{equation*}
\end{conj}

\begin{rem}
(1) As mentioned briefly in the introduction, this non-emptiness conjecture of Newton strata was proved for PEL-type Shimura varieties by C.-F. Yu \cite{Yu05} in Lie-type $C$ cases and by Viehmann-Wedhorn \cite{VW13} in general cases. Their methods of proof are both characteristic-$p$ geometric, which are not applicable for Hodge-type Shimura varieties.

(2) It is known (cf. \cite{RR96}, \cite{Chai00}) that with respect to the partial order $\preceq$ on $B(G_{\Qp},X)$, there exist a unique maximal element, called \textit{$\mu$-ordinary} element, and a unique minimal element, called \textit{basic} element. The $\mu$-ordinary element is nothing other than (the $\sigma$-conjugacy class in $B(G_{\Qp},X)$ with Newton point) $\overline{\mu}(G,X)$: it is clear from definition that $\overline{\mu}(G,X)\in B(G_{\Qp},X)$. Previously, non-emptiness of these two special strata have been known. For the $\mu$-ordinary locus, this was first proved by Wedhorn (\cite{Wedhorn99}) in the PEL-type cases, using an equi-characteristic deformation argument (which is thus yet unavailable for Hodge type Shimura varieties) and by Wortmann \cite{Wortmann13} for general Hodge-type cases, along the line of Viehmann-Wedhorn \cite{VW13} comparing the Newton stratification with  the Ekedahl-Oort stratification. 
As such, Wortmann's work makes an essential use of the Kisin's results on the integral $p$-adic comparison between the etale cohomology with $\Zp$-coefficients and the crystalline cohomology with integral coefficients, endowed with tensors (which is also established by \cite{Vasiu12}).
There was also a group-theoretic approach of Bultel \cite{Bultel01}. The $\mu$-ordinary locus is also known to be dense and open. In the PEL-type cases, non-emptiness of the basic locus was shown by Fargues \cite{Fargues04}.

(3) There is a conjectural formula for the (co)dimension of each Newton stratum, given by Chai \cite{Chai00}. The result of Kret \cite{Kret13}, alluded in the introduction, in fact proves this formula in the simple PEL-type cases of Lie type $A$ or $C$. His result depends on the resolution of Kisin \cite{Kisin13} of the Langlands-Rapoport conjecture and the stabilization of the twisted trace formula. We also mention the work of Scholze-Shin \cite[Cor.8.4]{SS13} of similar flavor.

(4) Originally, Conjecture \ref{conj:non-emptiness_of_Newton_strata} was formulated for general parahoric subgroup $K_p$. But, if $K_p$ is not hyperspecial, the $\mu$-ordinary locus may not be dense, as the example of Stamm \cite{Stamm97} shows (the Hilbert-Blumenthal surfaces with Iwahoric level structure at $p$). Still, our result (the first method of proof) proves that if there exists an integral model for which a suitable Neron-extension property similar to Theorem \ref{thm:integral_canonical_model} holds, then the Newton strata (parametrized by the same set $B(G_{\Qp},X)$) on the reduction of such model are all non-empty.

(5) We do not know yet that our definition of the isocrystal with $G_{\Qp}$-structure is an isocrystal with $G_{\Qp}$-structure on $\sS_K\otimes\kappa(\wp)$ in the sense of \cite[$\S$3]{RR96}. 
Due to this flaw, some basic properties on Newton stratification that are known for PEL-type Shimura varieties, such as the Grothendieck specialization theorem, are not yet established for Hodge type Shimura varieties.
For some known properties, we refer to \cite{Rapoport03}, \cite{Vasiu08}. 

\end{rem}


\section{Construction of special Shimura datum with prescribed $F$-isocrystal}

Suppose given a Shimura datum $(G,X)$ of Hodge type such that $G_{\Qp}$ is unramified. Fix a hyperspecial subgroup $K_p$ of $G(\Qp)$ and a prime $\wp$ of $E(G,X)$. We choose embeddings $\iota_{\infty}:\Qb\hra\C$ and $\iota_p:\Qb\hra\Qpb$ inducing $\wp$ on $E(G,X)$. Let $\sS_{K_p}(G,X)$ be the canonical integral model of the Shimura variety $\Sh_{K_p}(G,X)$. 
The goal of this section is to prove the fact that for any $F$-isocrystal with $G_{\Qp}$-structure $\overline{b}$ lying in $B(G_{\Qp},X)$, there exists a special Shimura sub-datum $(T,h\in \Hom(\dS,T_{\R})\cap X)$ such that $T_{\Qp}$ is unramified and the unique hyperspecial subgroup is contained in a $G(\Qp)$-conjugate of $K_p$, and that the $\sigma$-conjugacy class of $\mu_h^{-1}(p)\in G(\Qp^{\nr})$ equals $\overline{b}$.

Later, from such special Shimura datum, we will find a special Shimura sub-datum $(T',h')$ with unramified $T'_{\Qp}$ and $T'(\Zp)\subset K_p$, and further such that for any $g_f\in G(\A_f)$, the special point $x=[h',g_f]\in \sS_{K_p}(G,X)(\Qb)$ has good reduction at $\iota_p$ and the associated $F$-isocrystal with $G_{\Qp}$-structure of its reduction is $\overline{b}$.

\subsection{Construction of special Shimura datum with prescribed $F$-isocrystal satisfying condition $(\dagger)$} 
We will use the following notation in the remaining of this section.
Let $M$ be a reductive group over a field $F$ (of characteristic zero) and $T$ a maximal $F$-torus of $M$. For any $\mu\in X_{\ast}(T)$, we let $[\mu]_M$ denote its image under the natural map \begin{equation*}
X_{\ast}(T)\rightarrow \pi_1(M)_{\Gal(\overline{F}/F)}=(X_{\ast}(T)/\sum_{\alpha\in\Delta(T,M)}\Z\alpha^{\vee})_{\Gal(\overline{F}/F)}.\end{equation*} 
(Although the field $F$ does not enter the notation, in most cases, this will not cause confusion; $F$ will be the natural field over which $M$ is to be considered.)


Remember that we had fixed two embeddings $\Qb\hookrightarrow \C$, $\Qb\hookrightarrow\Qpb$, and using these we will regard every conjugacy class of cocharacters of $G_{\C}$ (for example, $c(G,X)$) as a conjugacy class of cocharacters of $G_{\Qpb}$. Suppose that

\begin{thm} \label{thm:existence_of_special_SD_with_given_F-isocrystal}
Let $(G,X)$ be a Shimura datum. Let $p$ be a rational prime such that $G_{\Q_p}$ is unramified, and $\overline{b'}\in B(G_{\Q_p})$. Suppose that 

$(\dagger)$ there exists a pair $(T',\mu')$, where $T'$ is an unramified maximal $\Q_p$-torus of $G_{\Qp}$, 
and $\mu'\in X_{\ast}(T')\cap c(G,X)$ such that the Newton point $\overline{\nu}_{G_{\Qp}}(\overline{b'})$ of $\overline{b'}\in B(G_{\Qp})$ equals the image of
\begin{equation*}\frac{1}{r}\Nm_{L_r/\Qp}(\mu')\ \in\ \mathcal{N}(T')=X_{\ast}(T')^{\Gamma(p)}_{\Q}\end{equation*}
under the canonical map $\mathcal{N}(T')\ra \mathcal{N}(G_{\Qp})$, where $T'$ splits over $L_r:=L(\F_{p^r})$.

Then, there exists a special Shimura datum $(T,h\in X\cap \Hom(\dS,T_{\R}))$ such that $(T_{\Qp},\mu_h)=\Int g_p (T',\mu')$ for some $g_p\in G(\Qp)$ and that $b'$ is $\sigma$-conjugate in $G(L)$ to $\mu_h^{-1}(p)\in G(\Qp^{\nr})$.
Moreover, for any finite prime $l\neq p$ of $\Q$, we can choose such $(T,h)$ with $T_{\Ql}$ being elliptic in $G_{\Ql}$
\end{thm}


We note that when $\overline{b'}\in B(G_{\Qp},X)$, the equality $\overline{\nu}_{G_{\Qp}}(\overline{b'})=\frac{1}{r}\Nm_{L_r/\Qp}(\mu')$ of elements in $\mathcal{N}(G_{\Qp})$ (i.e. the condition in $(\dagger)$) is equivalent to the equality 
\[ \overline{b'}=\overline{\mu'(p)}\] 
of elements in $B(G_{\Qp})$, because of
Lemma \ref{lem:isocrystal_of_special_pt} and the fact that the map 
\begin{equation*}
(\overline{\nu},\kappa):B(G_{\Q_p})\rightarrow \mathcal{N}(G_{\Q_p})\times \pi_1(G)_{\Gamma(p)} 
\end{equation*}
is injective (\cite[4.13]{Kottwitz97}).

\textsc{Proof.} Our proof consists of two steps.
Fix a rational prime $l$ different from $p$.

\textbf{Step 1.} First, we claim that if the assumption $(\dagger)$ holds, then it holds with $T'=(T_0)_{\Qp}$ for some maximal $\Q$-torus $T_0$ of $G$ which is elliptic at both $\R$ and $l$. More precisely, we find a maximal $\Q$-torus $T_0$ of $G$ such that $(T_0)_{\Qv}$ is elliptic in $G_{\Qv}$ for $v=\infty,l$ and that $(T_0)_{\Q_p}=\Int (g_p)(T')$ for some $g_p\in G(\Q_p)$. To that end, we consider the following three sets: 
\begin{align*}
X_{v}=\{x\in G^{\mathrm{sc}}(\Q_v) &\ |\ x \text{ is regular semisimple and }\\ &\quad   \Cent_G(x)\text{ is elliptic in }G_{\Q_v}\ \}\quad (v=\infty,l); \\
X_{p}=\{x\in G^{\mathrm{sc}}(\Q_p) &\ |\ x \text{ is regular semisimple and }\\ &\quad   \Cent_G(x)=\Int (g_p)(T')\text{ for some }g_p\in G(\Q_p)\ \},
\end{align*}
These sets are all non-empty: for $v=\infty$, this follows from the Deligne's condition on Shimura data (\cite[(2.1.1.2)]{Deligne77}), and for $v=l$, from \cite[Thm. 6.21]{PR94}. They are also open in the real, $l$-adic, and $p$-adic topology, respectively: ``for sufficiently close two regular semisimple elements, their centralizers are conjugated" (cf. \cite{Bultel01}, proof of Lemma 3.1). 
So, by weak approximation theorem, there exists an element $x\in G^{\mathrm{sc}}(\Q)$
which lies in all of them. Its centralizer $T_0:=\Cent_G(x)$ is then the desired maximal torus.

Clearly, for any $g\in G(\Q_p)$ (especially, $g\in G(\Q_p)$ with $(T_0)_{\Qp}=gT'g^{-1}$), the new pair 
$(T_1,\mu_1)=(gT'g^{-1}, g\mu' g^{-1})$ still satisfies the condition $(\dagger)$.
Hence, in the rest of our proof, we may and do assume that the condition $(\dagger)$ holds for a pair $(\mu',T')$ such that $T'=(T_0)_{\Qp}$ for some maximal $\Q$-torus $T_0$ of $G$ which is elliptic at $\R$ and at our chosen prime $l\neq p$.

\textbf{Step 2.} Next, for any maximal $\Q$-torus $T_0$ of $G$, elliptic at $\R$, and each cocharacter $\mu'\in X_{\ast}(T_0)$, we claim that there exists $u\in G(\Qb)$ such that 
\begin{itemize}
\item[(i)] $\Int u:(T_0)_{\Qb}\hra G_{\Qb}$ is defined over $\Q$ and the $\Q$-torus $T:=\Int u(T_0)$ is elliptic in $G$ over $\R$, 
\item[(ii)] $\Int u(\mu')$ equals $\mu_h^{-1}$ for some $h\in \Hom(\dS,T_{\R})\cap X$.
\end{itemize}
Moreover, if $(T_0)_{\Ql}$ is elliptic in $G_{\Ql}$ for some $l\neq p$, 
there exist such $u\in G(\Qb)$, $h\in X$ satisfying, in addition to these properties, that 
\begin{itemize}
\item[(iii)] there exists $y\in G(\Qp)$ such that $(T_{\Qp},\mu_h^{-1})=\Int y((T_0)_{\Qp},\mu')$.
\end{itemize}

Our proof is an adaptation of the argument of \cite{LR87}, Lemma 5.12. 
Let $T_0^{\mathrm{sc}}$ be the inverse image of $(T_0\cap G^{\mathrm{der}})^0$ under the canonical isogeny $G^{\mathrm{sc}}\ra G^{\mathrm{der}}$, and pick any morphism $h_0\in X$ factoring through $(T_0)_{\R}$, which exists since $(T_0)_{\R}$ is elliptic in $G_{\R}$. 
Choose $w\in N_{G}(T_0)(\C)$ such that $\mu'=w(\mu_{h_0}^{-1})$,%
\footnote{Recall our convention that $c(G,X)$ is the conjugacy class of cocharcters into $G_{\C}$ containing $\mu_h^{-1}$ for some $h\in X$.}
and consider the cocycle $\alpha^{\infty}\in Z^1(\Gal(\C/\R),G^{\mathrm{sc}}(\C))$ defined by 
\begin{equation*}\alpha^{\infty}_{\iota}=w\cdot\iota(w^{-1}).\end{equation*}
As a mater of fact, this has values in $T_0^{\mathrm{sc}}(\C)$. Indeed, by \cite[Prop. 2.2]{Shelstad79}, the automorphism $\Int(w^{-1})$ of $(T_0^{\mathrm{sc}})_{\C}$ is defined over $\R$, so $\Int(w^{-1})(\iota(t))=\iota(\Int(w^{-1})t)=\Int(\iota(w^{-1}))(\iota(t))$ for all $t\in T_0^{\mathrm{sc}}(\C)$, i.e. $\iota(w)w^{-1}\in \mathrm{Cent}_{G^{\mathrm{sc}}}(T_0^{\mathrm{sc}})(\C)=T_0^{\mathrm{sc}}(\C)$, and so is $w\iota(w^{-1})=\iota(\iota(w)w^{-1})$.
Then, by Lemma 7.16 of \cite{Langlands83}, there exists a global cocycle 
\begin{equation*}\alpha\in Z^1(\Q,T_0^{\mathrm{sc}})\end{equation*} mapping to $\alpha^{\infty}\in  H^1(\Q_{\infty},T_0^{\mathrm{sc}})$. When $(T_0)_{\Ql}$ is elliptic in $G_{\Ql}$ for some $l\neq p$, then the following lemma ensures the existence of a cocycle $\alpha\in Z^1(\Q,T_0^{\mathrm{sc}})$ mapping to $\alpha^{\infty}\in H^1(\Q_{\infty},T_0^{\mathrm{sc}})$ and further to zero in $H^1(\Qp,T_0^{\mathrm{sc}})$. For a finitely generated abelian group $A$, let $A_{\mathrm{tors}}$ be the subgroup of its torsion elements.  
\begin{lem} 
(1) For a maximal $\Q$-torus $T$ of $G$ which is elliptic at $l\neq p$, the natural map $(\pi_1(T^{\mathrm{sc}})_{\Gamma(l)})_{\mathrm{tors}} \ra (\pi_1(T^{\mathrm{sc}})_{\Gamma})_{\mathrm{tors}}$
is surjective.

(2) There exists $\alpha\in Z^1(\Q,T_0^{\mathrm{sc}})$ mapping to $\alpha^{\infty}\in H^1(\Q_{\infty},T_0^{\mathrm{sc}})$ and to zero in $H^1(\Qp,T_0^{\mathrm{sc}})$.
\end{lem}

\textsc{Proof of lemma.} (1) 
This is equal to the composite:
\begin{equation*}(\pi_1(T^{\mathrm{sc}})_{\Gamma(l)})_{\mathrm{tors}}\hra \pi_1(T^{\mathrm{sc}})_{\Gamma(l)}\twoheadrightarrow \pi_1(T^{\mathrm{sc}})_{\Gamma}\twoheadrightarrow (\pi_1(T^{\mathrm{sc}})_{\Gamma})_{\mathrm{tors}},\end{equation*}
where the last two maps are clearly surjective. Hence, it suffices to prove that $\pi_1(T^{\mathrm{sc}})_{\Gamma(l)}$ is a torsion group. 
But, since $T^{\mathrm{sc}}_{\Q_l}$ is anisotropic, $\widehat{T^{\mathrm{sc}}}^{\Gamma(l)}$ is a finite group, and so is 
\begin{equation*}
\pi_1(T^{\mathrm{sc}})_{\Gamma(l)}=X^{\ast}(\widehat{T^{\mathrm{sc}}}^{\Gamma(l)})=\Hom(\widehat{T^{\mathrm{sc}}}^{\Gamma(l)},\C^{\times}).
\end{equation*}

(2) For any place $v$ of $\Q$, there exists a canonical isomorphism (\cite[(3.3.1)]{Kottwitz84a}, \cite[1.14]{RR96})
\begin{equation*}H^1(F,T^{\mathrm{sc}})\isom \Hom(\pi_0(\widehat{T^{\mathrm{sc}}}^{\Gamma(v)}),\C^{\times})=((\widehat{T^{\mathrm{sc}}}^{\Gamma(v)})^D)_{\mathrm{tors}}\cong(\pi_1(T^{\mathrm{sc}})_{\Gamma(v)})_{\mathrm{tors}},\end{equation*}
and a short exact sequence (\cite[Prop.2.6]{Kottwitz86})
\begin{equation*}H^1(\Q,T^{\mathrm{sc}})\ra H^1(\Q,T^{\mathrm{sc}}(\overline{\A})):=\oplus_v H^1(\Qv,T^{\mathrm{sc}})
\stackrel{\theta}{\ra} \pi_0(\widehat{T^{\mathrm{sc}}}^{\Gamma})^D=(\pi_1(T^{\mathrm{sc}})_{\Gamma})_{\mathrm{tors}},\end{equation*}
Here, $\theta$ is the composite $\oplus_v H^1(\Qv,T^{\mathrm{sc}})\isom \oplus_v \pi_0(\widehat{T^{\mathrm{sc}}}^{\Gamma(v)})^D\ra \pi_0(\widehat{T^{\mathrm{sc}}}^{\Gamma})^D$, of which the second map is the one considered in (1) for $T=T_0$. 
Hence, by (1), there exists a class $\alpha^l\in H^1(\Ql,T^{\mathrm{sc}})$ such that $\theta(\alpha^{l})=-\theta(\alpha^{\infty})$, and thus the element $(\beta^v)_v\in H^1(\Q,T^{\mathrm{sc}}(\overline{\A}))$ defined by that $\beta^{v}=\alpha^{v}$ for $v=\infty,l$ and $\beta^{v}=0$ for $v\neq l,\infty$ maps to zero in $\pi_0(\widehat{T^{\mathrm{sc}}}^{\Gamma})^D$. By exactness of the sequence, we conclude existence of $\alpha\in Z^1(\Q,T^{\mathrm{sc}})$ whose cohomology class maps to $(\beta^v)_v$. $\qed$

Now, by changing $\alpha^{\infty}$ and $w$ further, we may assume that $\alpha^{\infty}$ \emph{is} the restriction of $\alpha$ to $\Gal(\C/\R)$. Then, since the restriction map $H^1(\Q,G^{\mathrm{sc}})\ra H^1(\Q_{\infty},G^{\mathrm{sc}})$ is injective (the Hasse principle), we see that $\alpha$ becomes trivial as a cocycle with values in $G^{\mathrm{sc}}(\Qb)$. We summarize this discussion in the diagram:
\begin{equation*}\xymatrix{
H^1(\Q_{\infty},T_0^{\mathrm{sc}})\ar[r] & H^1(\Q_{\infty},G^{\mathrm{sc}}) & \alpha^{\infty}_{\iota}=w\iota(w^{-1}) \ar@{|->}[r] & 0  & \\
H^1(\Q,T_0^{\mathrm{sc}})\ar[r] \ar[u] & H^1(\Q,G^{\mathrm{sc}}) \ar@{^{(}->}[u] & \exists\ \alpha \ar@{|->}[u]   \ar@{|->}[r] & \alpha'\ \ar@{|->}[u] & \Rightarrow \alpha'=0.} 
\end{equation*}
So there exists $u\in G^{\mathrm{sc}}(\Qb)$ such that for all $\rho\in\Gal(\Qb/\Q)$,
\begin{equation*}\alpha_{\rho}=u^{-1}\rho(u).\end{equation*}
Then, the homomorphism $\Int  u: (T_0)_{\Qb}\ra G_{\Qb}:t\ra t'=utu^{-1}$ is defined over $\Q$; in particular, $T:=u T_0 u^{-1}$ is a torus defined over $\Q$, and as the restriction of $\Int u$ to $Z(G)$ is the identity, $T _{\R}$ is also elliptic in $G_{\R}$. 
Moreover, since $u^{-1}\iota(u)=w\iota(w^{-1})$ for $\iota\in\Gal(\C/\R)$, one has that $uw\in G^{\mathrm{sc}}(\R)$ and \begin{equation*}\Int  u(\mu')=\Int (uw)\mu_{h_0}^{-1}=\mu_h^{-1}\end{equation*} for $h:=\Int  (uw)(h_0)\in X\cap \Hom(\dS,T)$. This establishes the claims (i), (ii). The proof of the claim (iii) is similar. As the restriction of $\alpha$ to $\Gal(\Qpb/\Qp)$ is trivial, there exists $x\in T(\Qpb)$ such that $x\rho(x^{-1})=\alpha_{\rho}=u^{-1}\rho(u)$ for all $\rho\in \Gal(\Qpb/\Qp)$, i.e. $ux\in G(\Qp)$. But, the homomorphism $\Int u:(T_0)_{\Qpb}\ra G_{\Qpb}$ also equals $\Int u=\Int ux$. This proves (iii).

To finish the proof, we observe that for the pair $((T_0)_{\Qp},\mu')$ produced at the end of Step 1, if $(T,h)$ is the special Shimura datum obtained from $(T_0,\mu')$ by the procedure of Step 2, by the property (iii),
$\overline{b'}$ equals the $\sigma$-conjugacy class of $\mu_h^{-1}(p)\in T(\Qp^{\nr})$. 
$\qed$


\subsection{Affine Deligne-Lusztig varieties and verification of the assumption $(\dagger)$}
In this subsection, we prove that in fact, the assumption $(\dagger)$ of Theorem \ref{thm:existence_of_special_SD_with_given_F-isocrystal} always holds (i.e. holds for every $\overline{b}\in B(G_{\Qp},X)$). For that, we use the fact (due to Witenberger \cite{Wintenberger05})) that for any $\overline{b}\in B(G,\{\mu\})$, the affine Deligne-Lusztig variety $X_{\mu}(b)$ is non-empty. From the latter fact, we give two ways of verifying the assumption $(\dagger)$.

\subsubsection{} 
Let $L$ be the completion of the maximal unramified extension of a non-archimedean local field $F$ in an algebraic closure $\overline{F}$ of $F$ and $\sigma\in\mathrm{Aut}(L/F)$ the relative Frobenius. Let $G$ be an unramified reductive group over $F$.%
%
For a conjugacy class $\{\mu\}$ of cocharacters of $G$ defined over $\overline{L}$, we have the finite set $B(G,\{\mu\})$ defined in (\ref{subsubsec:B(G,X)}). Let $S$ be a maximal $F$-split torus of $G$ and $T$ its centralizer (so $T$ is a maximal torus of $G$, since $G$ is quasi-split); corresponding to $S$, there exists an apartment $\mathbf{a}$ in the Bruhat-Tits building $\mathscr{B}(G,L)$ of the reductive group $G_L$ over $L$ (cf. \cite{Tits79}). Let $\widetilde{K}=\mathrm{Stab}(x)$ be the (hyper)special subgroup of $G(L)$ corresponding to a (hyper)special vertex $x$ of $\mathbf{a}$. Fix a uniformizer $\varpi$ of $F$. 
The \textit{affine Deligne-Lusztig variety} associated with $\mu\in X_{\ast}(T)/\Omega$ and $b\in G(L)$ is defined to be the set 
\begin{equation*}X_{\mu}(b)=\{ g\widetilde{K}\in G(L)/\widetilde{K}\ |\ g^{-1}b\sigma(g)\in \widetilde{K}\mu(\varpi) \widetilde{K}\}.\end{equation*}
When $b'\in G(L)$ is $\sigma$-conjugate to $b$, say $b'=h^{-1}b\sigma(h)$, the map $g\mapsto g'=h^{-1}g$ induces a bijection
\begin{equation} \label{eq:isom_of_affine_DL_varieties}
X_{\mu}(b)\isom X_{\mu}(b').
\end{equation}

Recall that a cocharacter $\mu'\in X_{\ast}(T')$ of a maximal torus $T'$ of $G_{\overline{F}}$ is said to be \emph{minuscule} if $\langle \alpha,\mu'\rangle\in\{-1,0,1\}$ for all simple roots $\alpha\in\Delta(G,T')$. It is well-known (\cite[1.2.5]{Deligne77}) that any Hodge cocharacter $\mu=\mu_h^{-1}\ (h\in X)$ is minuscule.

\begin{thm} \cite[Corollare 3]{Wintenberger05} \label{thm:non-emptiness_of_affine_Deligne-Lusztig_variety}
Suppose $\mu$ is minuscule. Then for any $\overline{b}\in B(G,\{\mu\})$, the set $X_{\mu}(b)$ is non-empty.
\end{thm}

\begin{lem} \label{eq:condition(dagger)}
For a connected reductive group $G$ over a $p$-adic field $F$ and a $G(\overline{F})$-conjugacy class $\mathcal{C}$ of cocharacters, the assumption $(\dagger)$ of Theorem \ref{thm:existence_of_special_SD_with_given_F-isocrystal} holds if there exists a quadruple \begin{equation*}(M',b',T',\mu'),\end{equation*} where $M'$ is a Levi subgroup of $G_{\Qp}$, $b'$ is a representative in $M'(L)$ of $\overline{b'}$ whose $\sigma$-conjugacy class in $B(M')$ is \emph{basic} (cf. \cite{Kottwitz85}), $T'$ is an unramified elliptic maximal $\Q_p$-torus of $M'$, and $\mu'\in X_{\ast}(T')\cap \mathcal{C}$ such that \begin{equation*}\kappa_{M'}(\overline{b'})=[\mu']_{M'}\end{equation*} (here, $\overline{b'}$ denoting the $\sigma$-conjugacy class of $b'$ in $B(M')$).

\end{lem}
\textsc{Proof.} It suffices to show that $b'\in M'(L)$ is  $\sigma$-conjugate (in $M'(L)$) to $\mu'(p)\in T'(L)$. 
Since $T'$ is an \emph{elliptic} maximal torus of $M'$, the image of $B(T')$ in $B(M')$ equals the subset $B(M')_{b}$ of basic elements (\cite{Kottwitz85}, Prop. 5.3), and there is a commutative diagram of natural maps with bijective vertical maps (\cite{Kottwitz85}, Prop. 5.6) 
\begin{equation} \label{eq:functoriality_diagram_of_kappa_T,M}
\xymatrix{
B(T') \ar[r] \ar[d]^{\cong}_{\kappa_{T'}} & B(M')_{\mathrm{basic}} \ar[d]^{\cong}_{\kappa_{M'}} \\
X_{\ast}(T')_{\Gamma(p)} \ar[r] & \pi_1(M')_{\Gamma(p)}.}
\end{equation}
Denoting again (by abuse of notation) by $\overline{\mu'(p)}$ the class of $\mu'(p)$ in $B(T')$, 
we have \begin{equation*}\kappa_{T'}(\overline{\mu'(p)})=[\mu']_{T'}\end{equation*} (\cite{RR96}, Thm. 1.15 (ii)).
So, since the elements $\overline{\mu'(p)}$, $\overline{b'}$ of $B(M')$ are both \emph{basic}, their equality can be checked after applying $\kappa_{M'}$. Then, the equality $\kappa_{M'}(\overline{b'})=\kappa_{M'}(\overline{\mu'(p)})$ is just the assumption $\kappa_{M'}(\overline{b'})=[\mu]_{M'}$. $\qed$

Recall that for a Levi subgroup $M$ of $G$, an element $b$ of $M(L)$ is said to be \textit{$G$-regular basic} for $M$ if the $\sigma$-conjugacy class of $b$ in $B(M)$ is basic (so that $\nu_b$ is defined over $\Qp$, cf. \cite[(4.4.3)]{Kottwitz85}) and $\Cent_{G_{\Qp}}(\nu_b)$ is equal to $M$ (cf. \cite[$\S$6]{Kottwitz85}). A Levi subgroup of $G$ is said to be \textit{standard} if it contains the centralizer $\Cent_{G}(S)$ of a maximal $\Qp$-split torus $S$ of $G$ ($\Cent_{G}(S)$ is a maximal torus since $G$ is quasi-split).

\begin{prop} \label{prop:existence_of_Levi_and_Hodge_cocharacter}
Let $G$ be an unramified connected reductive group over $\Qp$ and $\mathcal{C}$ be a $G(\Qpb)$-conjugacy class of  cocharacters, containing a \emph{minuscule} one $\mu_0:\mathbb{G}_{m\Qpb}\rightarrow G_{\Qpb}$. Let $G_{\Zp}$ be a reductive $\Zp$-group scheme with generic fibre $G_{\Zp}\otimes_{\Zp}\Qp=G$ and put $K_p:=G_{\Zp}(\Zp)$. Let $B(G,\mathcal{C})$ be the subset of $B(G)$ as defined in (\ref{subsubsec:B(G,X)}). 

(1) For any $\overline{b}\in B(G,\mathcal{C})$, there exists a pair $(T,\mu)$ satisfying the assumption $(\dagger)$ of Theorem \ref{thm:existence_of_special_SD_with_given_F-isocrystal} and furthermore such that the unique hyperspecial subgroup of $T(\Qp)$ is contained in $K_p$.

(2) Suppose that one of the following two conditions holds:

(i) $G^{\mathrm{der}}$ is simply connected; 

(ii) there exists a representation $\rho:G_{\Qpb}\ra\mathrm{GL}(V)$ over $\Qpb$ with finite kernel and such that every weight of $\rho\circ\mu_0$ is either $0$ or $1$; 

Then, for every triple
\begin{equation*}(M,b,T),\end{equation*}
where $M$ is a standard Levi $\Qp$-subgroup of $G_{\Qp}$ splitting over $L$, $b\in M(L)$ is a basic element (i.e. $\overline{b}\in B(M)_{basic}$), and $T$ is a maximal torus of $M$, there exists $\mu\in X_{\ast}(T)\cap \mathcal{C}$ such that 
\begin{equation*}\kappa_M(\overline{b})=[\mu]_{M}.\end{equation*}

(3) Under the assumption of (2), for any $\overline{b}\in B(G,\mathcal{C})$, there exists a quadruple $(M,b,T,\mu)$ as in (2), which further satisfies that $b\in M(L)$ is $G_{\Qp}$-regular basic, $T$ is an unramified elliptic maximal torus of $M$, and that the unique hyperspecial subgroup of $T(\Qp)$ is contained in a $G(\Qp)$-conjugate of $K_p$.
\end{prop}

Note that clearly the condition (ii) of (2) is satisfied when $(G,\mathcal{C})$ arises (via obvious base-change from $\Q$ to $\Qp$) from a Hodge-type Shimura datum. 

In the cases of the assumptions of (2), by Lemma \ref{eq:condition(dagger)}, 
(2) and (3) together give another proof of the fact that the condition $(\dagger)$ as well as the claim on the hyperspecial subgroup hold for every $\overline{b}\in B(G,\mathcal{C})$, with certain simplification
(i.e. simple proof of (\ref{eq:Lemma 2.3.3_Kottwitz84b})) in the Hodge-type case (ii). 

\textsc{Proof.}
(1) This is what is established in the proof of Lemma 5.11 of \cite{LR87}; a (slightly different) account of this proof is also given in \cite[Lem. (2.2.2)]{Kisin13}. From the start, Langlands-Rapoport considers only isocrystals $\overline{b}\in B(G_{\Qp})$ with non-empty $X_{\mu}(b)$. Also, in their argument, $T$ can be any (unramified) maximal $\Qp$-torus of $G_{\Qp}$. So in our case, we can simply take $T$ to be an unramified maximal $\Qp$-torus of $G_{\Qp}$ such that the unique hyperspecial subgroup of $T(\Qp)$ is contained in $K_p$ (which exists by Lemma \ref{lem:existence_of_aniostropic_and_unramified_torus}).
Here, for the convenience of readers, we give a sketch of the argument of \cite[Lem. 5.11]{LR87}; we also use the notations of loc.cit. 

Let $M$ be the centralizer of a maximal $\Qp$-split sub-torus of $T$. It is a Levi subgroup of a parabolic subgroup $Q$ of $G$ defined over $\Qp$. Note that $T$ is an elliptic maximal torus of $M$ and $M$ is unramified (since $T$ is assumed to be unramified). We may assume that $M(L)\cap G_{\Zp}(\cO_L)$ is a special subgroup of $M(L)$ (\cite{Tits79}). Then, choose an \emph{unramified} maximal $\Qp$-torus $T'$ of $M$ such that $T'$ contains a maximal $\Qp$-split torus whose corresponding apartment in the Bruhat-Tits building contains the hyperspecial point corresponding to $K_p=G_{\Zp}(\Zp)$ (so that the hyperspecial subgroup of $T'(\Qp)$ is contained in $K_p$); from the beginning, one can assume that $T$ is such a torus and take $T'=T$ (Lemma \ref{lem:existence_of_aniostropic_and_unramified_torus}). Choose representatives $b\in M(L)$ of $\overline{b}$ and $\mu'\in X_{\ast}(T')$ of $\mathcal{C}$. 
Then, the following facts are established step by step: 
\begin{itemize}
\item[(i)] There exist $m\in M(L)$ and $n'\in N(L)$, where $N$ is the unipotent radical of $Q$, such that 
\begin{equation*}m^{-1}b\sigma(m)n' \in K_L \mu'(p) K_L\quad (K_L:=G_{\Zp}(\cO_L)).\end{equation*}
\item[(ii)] Let $\mu''\in X_{\ast}(T')$ be defined, via the Cartan decomposition for $M$ (\cite[3.3.3]{Tits79}), by that
\begin{equation} \label{eq:defining_eqn_of_mu''}
m^{-1}b\sigma(m) \in (K_L\cap M(L))\mu''(p)(K_L\cap M(L)).
\end{equation}
Then, $\mu''$ is conjugate to $\mu'$.
\item[(iii)] Let $\nu_{\overline{b}}$ be a $\Qp$-rational representative of $\overline{\nu}_{\overline{b}}\in \mathcal{N}(G)=(\mathrm{Int} G(L) \backslash \Hom_{L}(\mathbb{D},G_L))^{\langle\sigma\rangle}$ that factors through $T$ (which exists by Prop. 6.2 of \cite{Kottwitz85}) and let $\mu$ be a cocharacter of $T$ that is conjugated to $\mu''$ in $M$. Then for any $L_r=L(\F_{p^r})$ splitting $T$,
we have an equality in $X_{\ast}(T)$
\begin{equation*}r\nu_{\overline{b}}=\Nm_{L_r/\Qp}(\mu)\end{equation*}
\end{itemize}
The statement (i) is a consequence of that $X_{\mu'}(b)\neq \emptyset$.

The claim (ii) is the technical heart of the proof. It is proved using \cite{Kottwitz84b}, Lemma 2.3.3 and 2.3.7: the latter lemma uses Satake transform, and the former lemma can be given a simple proof in Hodge-type situation (the case (ii) of the statement (2) of this proposition, cf. proof of (\ref{eq:Lemma 2.3.3_Kottwitz84b}) below).
We note that the equation (\ref{eq:defining_eqn_of_mu''}) implies
\begin{equation} \label{eq:lambda(b)=mu_M}
\lambda_M(m^{-1}b\sigma(m))=[\mu'']_M.
\end{equation}
Here, $\lambda_M:M(L)\ra X^{\ast}(Z(\widehat{M}))$ is the map defined in \cite[Lem.3.3]{Kottwitz84b} (strictly speaking, the map in loc.cit. is obtained from this map by taking Galois invariants, i.e. $H^0(\Qp,-)$), and $[-]_M$ denotes (by abuse of notation) the natural map $X_{\ast}(T')\ra \pi_1(M)$. Also, $\lambda_M$ induces the map $\kappa_M$ of Theorem \ref{thm:Thm.1.8-RR96} by taking $H^1(\Qp,-)$. 
%

The proof of (iii) goes as follows: clearly, the equality in (iii) is equivalent to the statement that for each $\Qp$-rational character $\chi$ of $T$, one has that
\begin{equation*}\langle\chi,\mu\rangle=\langle\chi,\nu_{\overline{b}}\rangle\end{equation*}
(this is the same as (5.n) on p.178 of \cite{LR87}: $\nu_p$ in (5.n) equals $r\nu_{\overline{b}}=[L_r:\Qp]\nu_{\overline{b}}$ in our notation, and the sign difference is due to that our sign convention of $c(G_{\Qp},X)$ is different from that of Langlands-Rapoport).
But because $T$ is \emph{elliptic} in $M$, each $\Qp$-rational character of $T$ can be also regarded as a character of $M$ (equivalently, of the maximal abelian quotient $M^{\mathrm{ab}}$ of $M$), and as $\langle\chi,\mu\rangle=\langle\chi,\mu''\rangle$ for any character $\chi$ of $M$, the claim (iii) is the same as that for every $\Qp$-rational character $\chi$ of $M$
\begin{equation*}\langle\chi,\mu''\rangle=\langle\chi,\nu_{\overline{b}}\rangle,\end{equation*}
and this is a direct consequence of the equation (\ref{eq:defining_eqn_of_mu''}), as shown in \cite[p.178-179]{LR87}. Alternatively, 
assuming that $T$ is an unramified and taking $T'=T$, $\mu=\mu''$, 
we can also proceed as follows: as $T$ is elliptic in $M$, $\nu_{\overline{b}}$ is central in $M$ (i.e. $m^{-1}b\sigma(m)\in M(L)$ is \emph{basic}), hence in view of the equation (\ref{eq:lambda(b)=mu_M}), $(M,m^{-1}b\sigma(m),T,\mu'')$ satisfies the conditions of Lemma \ref{eq:condition(dagger)}, which gives that the two $\Qp$-rational cocharacters $r\nu_{\overline{b}}$, $\Nm_{L_r/\Qp}(\mu'')$ of $T$ are conjugated in $M(L)$, so their projections to $M^{\mathrm{ab}}$ are equal. But, they factor through the maximal $\Qp$-split sub-torus $T_1$ of $T$, and as $T$ is elliptic in $M$, the natural map $T_1\ra M^{\mathrm{ab}}$ induces an injection $X_{\ast}(T_1)\hra X_{\ast}(M^{\mathrm{ab}})$. So, we must have that $r\nu_{\overline{b}}=\Nm_{L_r/\Qp}(\mu'')$.

(3) Let us first prove (3) assuming (2). Let $S$ be any maximal $\Qp$-split $\Qp$-torus of $G$. According to \cite[Prop.6.2]{Kottwitz85}, there exist a standard Levi subgroup $M$ of $G$ and $b\in M(L)$ with $b\in \overline{b}$ such that $b$ is basic $G$-regular for $M$. In fact, $M$ can be taken to be the centralizer of \emph{any} $\Qp$-rational representative $\nu_p\in \Hom_{\Qp}(\mathbb{D},G)$ of the Newton point $\overline{\nu}_{G}(\overline{b})\in \mathcal{N}(G)=(\mathrm{Int} G(L) \backslash \Hom_{L}(\mathbb{D},G_L))^{\langle\sigma\rangle}$ (proof of Prop. 6.2 of \cite{Kottwitz85}). 
Furthermore, there exists an unramified elliptic maximal $\Q_p$-torus $T$ of $M$ (Lemma \ref{lem:existence_of_aniostropic_and_unramified_torus}). So far these discussions concern $\Qp$-groups.

For $\Zp$-groups, we can choose a maximal $\Qp$-split $\Qp$-torus $S$ such that in the Bruhat-Tits building $\mathcal{B}(G,\Qp)$, the hyperspecial point $x$ fixed by $K_p$ is contained in the apartment associated with $S$. For such $S$, let $S_{\Zp}$ be the reductive $\Zp$-group scheme with generic fiber being $S$ and such that $S_{\Zp}(\Zp)$ is the (unique) hyperspecial subgroup of $S(\Qp)$.
We claim that
$S_{\Zp}$ and $T'_{\Zp}:=\mathrm{Cent}_{G_{\Zp}}(S_{\Zp})$ are smooth, closed subschemes of $G_{\Zp}$.
To see that $S_{\Zp}$ is a closed subscheme of $G_{\Zp}$, we first observe that $S_{\Zp}(\cO_L)$ maps to $G_{\Zp}(\cO_L)$ under the embedding $S\hra G$, where $\cO_L$ is the ring of integers of $L$. This follows from the two facts: first, the hyperspecial point $x$ of $\mathcal{B}(G,\Qp)$ is also a hyperspecial point of $\mathcal{B}(G,L)$ with the associated $\cO_L$-group scheme being $(G_{\Zp})_{\cO_L}$ (\cite[3.4, 3.8]{Tits79}), and secondly, for a maximal $L$-split torus $S_1$ of $G$ containing $S$ and defined over $\Qp$ (which exists since $G$ is quasi-split), $x\in \mathcal{B}(G,L)$ lies in the apartment corresponding to $S_1$ (\cite[1.10]{Tits79}). So, $x$ is fixed under the maximal bounded subgroup $S_1(\cO_L)$ of $S_1(L)$, and thus under $S_{\Zp}(\cO_L)$ as well. Then, according to \cite[3.1.2.1]{Vasiu99}, this implies that the $\Qp$-embedding $S\hra G$ extends to a $\Zp$-embedding $S_{\Zp}\hra G_{\Zp}$. 
The fact that
$T'_{\Zp}:=\mathrm{Cent}_{G_{\Zp}}(S_{\Zp})$ is a smooth closed group $\Zp$-subscheme of $G_{\Zp}$ is well-known (\cite[XI, 5.3]{SGA3}, \cite[2.2]{Conrad11}). Moreover, one knows that it represents the functor on $\Zp$-schemes:
\begin{equation*}R'\mapsto \{g\in G(R')\ |\ g(S_{\Zp})_{R'}g^{-1}=(S_{\Zp})_{R'}\},\end{equation*}
and its generic fibre $T':=\mathrm{Cent}_{G}(S)$ is a maximal torus of $G$ which splits over $L$, since $G_{\Qp}$ is quasi-split (\cite[1.10]{Tits79}).  
Next, there exists a $\Qp$-rational representative $\nu_p$ of $\overline{\nu}_G(b)\in \Hom_{\Qp}(\mathbb{D},G)^{\langle\sigma\rangle}$ which factors through $S$ (\cite[(1.1.3)]{Kottwitz84b}). Suppose that $N\nu_p\ (N\in \N)$ is a morphism from $\Gm$ to $S\subset T'$ (not just a quasi-cocharacter). Then, the image of $N\nu_p$ is a $\Qp$-split subtorus of $T'$, and since $T'$ extends over $\Zp$, this image extends uniquely to a $\Zp$-subtorus $U$ of $T'_{\Zp}$ (use the equivalence between the category of $\Qp$-tori and the category of free $\Z$-modules with $\Gamma(p)$-action, under which $\Zp$-tori correspond to unramified $\Gamma(p)$-modules). The morphism $U\ra T'_{\Zp}$ is a closed embedding (\cite[IX.2.5]{SGA3}). 

Now, as mentioned above, we may take $\mathrm{Cent}_{G}(N\nu)$ for $M$. It is the generic fibre of $M_{\Zp}:=\mathrm{Cent}_{G_{\Zp}}(U)$, which is a reductive group scheme over $\Zp$ (\cite[XI, 5.3]{SGA3}, \cite[2.2]{Conrad11}). Therefore, by Lemma \ref{lem:existence_of_aniostropic_and_unramified_torus}, there exists a maximal torus $T_{\Zp}$ defined over $\Zp$ of $M_{\Zp}$ whose generic fibre $T=T_{\Zp}\otimes\Qp$ is elliptic in $M=M_{\Zp}\otimes\Qp$ (and thus is unramified over $\Qp$). Finally, the claim follows since the unique hyperspecial subgroup of $T(\Qp)$ is $T_{\Zp}(\Zp)$ (which is also the unique maximal compact subgroup of $T(\Qp)$).

(2) Following \cite{Kottwitz03}, $\S$4.2, for arbitrary $\mu\in X_{\ast}(T)$, we put
\begin{equation*} \mathcal{P}_{\mu}:=\{\nu\in X_{\ast}(T)\ |\ [\nu]_G=[\mu]_G,\ \nu\in\mathrm{Conv}(\Omega_G\cdot\mu)\},
\end{equation*}
where $\mathrm{Conv}(\Omega_G\cdot\mu)$ denotes the convex hull in $X_{\ast}(T)_{\R}$ of the orbit $\Omega_G\cdot\mu$, $\Omega_G$ being the absolute Weyl group of $T$ in $G$. In the original definition of $\mathcal{P}_{\mu}$, instead of our $T$, Kottwitz used a maximal torus of $M$ containing a maximal $\Qp$-split torus of $G_{\Qp}$, but there is an obvious bijection between $\mathcal{P}_{\mu}$'s obtained using two different maximal tori of $M$. 
We also remark that when $G^{\mathrm{der}}$ is simply-connected, the first condition $[\nu]_G=[\mu]_G$ in the definition of $\mathcal{P}_{\mu}$ is implied by the second one (under this assumption, the second condition implies that $\mu-\nu$ lies in the coroot lattice $\sum_{\alpha\in\Delta}\Z\alpha^{\vee}=X_{\ast}(T\cap G^{\mathrm{der}})$; see below for a proof).

Let $\mathcal{P}_{\mu,M}$ denote the subset of $\pi_1(M)_{\Gamma(p)}$ (which was denoted by $Y_M$ in \cite{Kottwitz03}) obtained as the image of $\mathcal{P}_{\mu}$ under the canonical surjection $[\ \cdot\ ]_M:X_{\ast}(T)\rightarrow\pi_1(M)_{\Gamma(p)}$.
Then, Kottwitz shows (\cite[Thm. 4.3]{Kottwitz03}) that for every \emph{standard} Levi subgroup $M$ and any basic $b'\in M(L)$ and $\mu'\in X_{\ast}(T)$, we have that
\begin{equation*}\kappa_M(\overline{b'})\in \mathcal{P}_{\mu',M}\quad \text{ if and only if }\quad X_{\mu'}(b')\neq\emptyset.\end{equation*} 
On the other hand, the result of Wintenberger (Theorem \ref{thm:non-emptiness_of_affine_Deligne-Lusztig_variety}) tells us that for any \textit{minuscule} cocharacter $\mu'$ and any $b'\in G(L)$, 
\begin{equation*}X_{\mu'}(b')\neq\emptyset\quad \text{ if and only if }\quad  \overline{\nu}(\overline{b'})\in B(G,\{\mu'\}).\end{equation*} Hence, it follows from these together that for any $\overline{b'}\in B(G,\{\mu'\})$, $\kappa_M(\overline{b'})$ lies in $\mathcal{P}_{\mu',M}$. 

Finally, it remains to show that 
\begin{equation*}
\mathcal{P}_{\mu',M}(=[\mathcal{P}_{\mu'}]_M)=[\Omega_G\cdot\mu']_M,
\end{equation*}
In turn, since one has $[\nu]_G=[\mu']_G$ for any $\nu\in \Omega_G\cdot\mu'$, it suffices to show that
\begin{equation} \label{eq:Lemma 2.3.3_Kottwitz84b}
X_{\ast}(T)\cap \mathrm{Conv}(\Omega_G\cdot\mu')=\Omega_G\cdot\mu'.
\end{equation}
First, suppose that $G^{\mathrm{der}}$ is simply connected. Then, for any \emph{minuscule} cocharacter $\mu'\in X_{\ast}(T)$, one has 
\begin{equation*}X_{\ast}(T)\cap \mathrm{Conv}(\Omega_G\cdot\mu')\stackrel{(a)}{=}(\mu'+Q^{\vee})\cap \mathrm{Conv}(\Omega_G\cdot\mu')\stackrel{(b)}{=}\Omega_G\cdot\mu',\end{equation*}
where $Q^{\vee}=\sum_{\alpha\in \Delta} \Z\alpha^{\vee}$ is the coroot lattice.
The equality (\textit{a}) is due to the fact (\cite[Lemma 2.2]{RR96}) that when one fixes a set $\Delta$ of simple roots and $\mu'$ lies in the associated closed Weyl chamber $\overline{C}\subset X_{\ast}(T)_{\R}$, 
$\mu\in \mathrm{Conv}(\Omega_G\cdot\mu')$ implies that $\mu'-w\mu\in C^{\vee}:=\sum_{\alpha\in \Delta} \R_{\geq0}\alpha^{\vee}
\subset X_{\ast}(T^{\mathrm{der}})_{\R}$ for every $w\in\Omega$, and the fact that 
\begin{equation*} \label{eq:coroot_lattice=cocharacter_lattice}
Q^{\vee}=X_{\ast}(T^{\mathrm{der}})=X_{\ast}(T)\cap X_{\ast}(T^{\mathrm{der}})_{\R}\quad (\subset X_{\ast}(T)_{\R})
\end{equation*}
for simply connected $G^{\mathrm{der}}$. To see the second equality, we note that since \emph{$G^{\mathrm{der}}$ is simply connected}, $(T\cap G^{\mathrm{der}})^0=T\cap G^{\mathrm{der}}$, thus the submodule $X_{\ast}(T^{\mathrm{der}})$ of $X_{\ast}(T)$ is a direct summand (as follows, for example, by applying \cite[Cor. 2.3]{Kottwitz84b} to the exact sequence of diagonalizable groups $1\ra T\cap G^{\mathrm{der}}\ra T\ra G^{\mathrm{ab}}\ra 1$, to get an exact sequence of free $\Z$-modules).
\footnote{We also note the fact that when $G^{\mathrm{der}}$ is simply connected, for $x,x'\in X_{\ast}(T)$, if $x-x'\in X_{\ast}(T^{\mathrm{der}})_{\R}$, one has that $x'\geq x$ ($x'\succeq x$ in our notation) if and only if $x'-x\in\sum_{\alpha\in \Delta} \Z_{\geq0}\alpha^{\vee}$.}
The equality (\textit{b}) is \cite[Lemma 2.3.3]{Kottwitz84b}. 

Next, suppose that we are in the situation of case (ii). 
Fix a maximal torus $T'$ of $\mathrm{GL}(V)$ containing $\rho(T)$. Since the induced map $\rho_{\ast}:X_{\ast}(T)\ra X_{\ast}(T')$ is linear, we have the inclusion
\begin{equation*}\rho_{\ast}(X_{\ast}(T)\cap\mathrm{Conv}(\Omega_G\cdot\mu'))\subset X_{\ast}(T')\cap\mathrm{Conv}(\rho_{\ast}(\Omega_G\cdot\mu')).\end{equation*}
Also, since $\rho$ has finite kernel, $\rho_{\ast}$ is injective. 
Let $\{e_1,\cdots,e_n\}$ be the weights of the representation $T'\ra\mathrm{GL}(V)$; clearly, these form a basis $X^{\ast}(T)$, giving an identification $X^{\ast}(T')=\Z^n$, $X_{\ast}(T')=\Z^n$.  Likewise, the set $\{e_1\circ\rho,\cdots,e_ n\circ\rho\}$ consists of weights of $(T,\rho|_T)$, and for any $\mu\in X_{\ast}(T)$, each $e_i\circ\rho_{\ast}(\mu)\in\Z\ (i=1,\cdots,n)$ is also a weight of $\rho\circ\mu:\Gm\ra\mathrm{GL}(V)$. 
So, by assumption (which also implies that for any $\mu\in \Omega_G\cdot \mu$, every weight of $\rho\circ\mu$ is either $0$ or $1$), the set $\rho_{\ast}(\Omega_G\cdot\mu')\subset\Z^{n}$ is contained in the set of vertices of the hypercube $B_{n}:=\{(x_i)\in \R^{n}\ |\ 0\leq x_i\leq 1\}$; in particular, 
$\mathrm{Conv}(\rho_{\ast}(\Omega_G\cdot\mu'))\subset \R^{n}$ is contained in the hypercube $B_{n}$. Then, since $X_{\ast}(T')=\Z^{n}\subset\R^{n}$, the intersection $X_{\ast}(T')\cap\mathrm{Conv}(\rho_{\ast}(\Omega_G\cdot\mu'))$ must be equal to the vertex set of $\mathrm{Conv}(\rho_{\ast}(\Omega_G\cdot\mu'))$, and thus a subset of $\rho_{\ast}(\Omega_G\cdot\mu')$. Since $\rho_{\ast}$ is injective, this proves (\ref{eq:Lemma 2.3.3_Kottwitz84b}).
$\qed$

\begin{rem}
(1) The main ideas of the two proofs of the condition $(\dagger)$ are quite similar. But, the details are still different. For example, the Levi subgroup $M$ in the proof of (1) is in general different (albeit related) from the Levi subgroup appearing in (2) and (3): the latter is the centralizer of a $\Qp$-rational representative of the conjugacy class $\overline{\nu}_{\overline{b}}$, while the former is not necessarily so. 

(2) In fact, the claim of (2) holds more generally for any Levi subgroup \emph{splitting over $L$} (not just for standard Levi subgroups). 
Indeed, we only need to check that for such Levi subgroup,
\begin{equation*}\kappa_M(\overline{b'})\in \mathcal{P}_{\mu',M}\quad \text{ if }\quad X_{\mu'}(b')\neq\emptyset.\end{equation*}
Indeed, this direction of Theorem 4.3 of \cite{Kottwitz03} follows from Lemma 4.5, Lemma 4.6 there; this part does not use Lemma 4.4 of loc. cit., whose proof seems to require the group to be unramified, and Lemma 4.6 is Lemma 4.5 of \cite{KR03}. Both of these lemmas (and the mathematical facts used for their proofs, such as the Iwasawa (resp. Cartan) decompositions for $G_L$ (resp. for $M_L$)) are statements about either reductive groups defined over $F$ and splitting over $L$ (Lemma 4.5), or split reductive groups over $L$ (Lemma 4.6 of \cite{Kottwitz03}, i.e. Lemma 4.5 of \cite{KR03}); one also notes that the same argument in fact establishes the stronger fact that if $X_{\mu'}(b')\neq\emptyset$, the subset $\lambda_{M}(\overline{b'})$ of $\pi_1(M)=X^{\ast}(Z(\widehat{M}))$ equals the image of $\mathcal{P}_{\mu'}$ under the canonical surjection $X_{\ast}(T)\ra \pi_1(M)$ (in the introduction of \cite{Kottwitz03}, $\lambda_{M}$ was denoted by $w_M:G(L)\ra \pi_1(M)$). 

\end{rem}

\subsection{First proof of non-emptiness of Newton strata}

We now give our first proof of the main theorem.

\begin{thm} \label{thm:existence_of_special_point_with_good_reduction_and_prescribed_given_isocrystal}

Let $(G,X)$ be a Shimura datum. Let $p$ be a rational prime such that $G_{\Qp}$ is unramified, and fix a hyperspecial subgroup $K_p$ of $G(\Qp)$. Choose a prime $\wp$ of $E(G,X)$ above $p$ and an embedding $\iota_p:\Qb\hra\Qpb$ inducing $\wp$.

Then, for every $\overline{b}\in B(G_{\Qp},X)$, there exists a special sub-Shimura datum $(T,h:\dS\ra T_{\R})$ such that $T_{\Qp}$ is unramified and the unique hyperspecial subgroup of $T(\Qp)$ is contained in $K_p$, and that the $\sigma$-conjugacy class of $\mu_h^{-1}(p)\in G(\Qp^{\nr})$ equals $\overline{b}$.

In particular, if $(G,X)$ is of Hodge type, for every $g_f\in G(\A_f)$, the special point $[h,g_f]_{K_p}\in \Sh_{K_p}(G,X)(\Qb)$ has good reduction at $\iota_p$ and the $F$-isocrystal of its reduction in $\sS_{K_p}(\Fpb)$ is $\overline{b}$.
\end{thm}

Recall that a special point of $\Sh_K$ is said to have \textit{good reduction at $\iota_p$} if it
extends as a point in $\sS_K(\cO_{\iota_p})$ for the valuation ring $\cO_{\iota_p}$ of $\Qb$ defined by $\iota_p$. 

\textsc{Proof.} We can apply Theorem \ref{thm:existence_of_special_SD_with_given_F-isocrystal} and Proposition \ref{prop:existence_of_Levi_and_Hodge_cocharacter}, to obtain a special Shimura datum $(T_1,h_1:\dS\ra (T_1)_{\R})$ such that $(T_1)_{\Qp}$ is unramified and the unique hyperspecial subgroup of $T_1(\Qp)$ is contained in $\Int(g_p)(K_p)$ for some $g_p\in G(\Qp)$ and that the $\sigma$-conjugacy class of $\mu_{h_1}^{-1}(p)\in G(\Qp^{\nr})$ equals $\overline{b}$. 
Since $G(\Q)$ is dense in $G(\Qp)$ (\cite[Lem. 4.10]{Milne94}), there exists $g\in G(\Q)\cap K_p\cdot g_p^{-1}$. Then, the new special Shimura datum $(T,h):=\Int(g)(T_1,h_1)$ still enjoys the same properties and further satisfies that the unique hyperspecial subgroup of $T(\Qp)$ is contained in $K_p$. 
If moreover $(G,X)$ is of Hodge type, according to Lemma \ref{lem:isocrystal_of_special_pt}, for every $g_f\in G(\A_f)$, the special point $[h,g_f]_{K_p}\in \Sh_{K_p}(G,X)(\Qb)$ has good reduction at $\iota_p:\Qb\hra\Qpb$ and the $F$-isocrystal of its reduction in $\sS_{K_p}(G,X)(\Fpb)$ is $\overline{b}$.
$\qed$



\begin{cor} \label{cor:non-emptiness_criterion_of_ordinary_locus}
Let $(G,X)$ be a Shimura datum of Hodge type. Choose an embedding $\iota_p:\Qb\hra \Qpb$ and let $\wp$ be the prime of $E(G,X)$ induced by $\iota_p$. Suppose that $G_{\Q_p}$ is unramified and choose a hyperspecial subgroup $K_p$ of $G(\Qp)$. Then the reduction $\sS_{K_p}(G,X)\times\kappa(\wp)$ has non-empty ordinary locus if and only if $\wp$ has absolute height one (i.e. $E(G,X)_{\wp}=\Q_p$).
\end{cor}

\textsc{Proof.} The necessity of $E(G,X)_{\wp}=\Q_p$ was proved in \cite[Section3]{Bultel01}. 
For sufficiency, we note that it suffices to find a special Shimura sub-datum $(T,\{h\})$ with the following properties:

(i) $T_{\Qp}$ is unramified and the unique hyperspecial subgroup of $T(\Qp)$ is contained in $K_p$;

(ii) there exists a Borel subgroup $B$ over $\Qp$ of $G_{\Qp}$ containing $T_{\Qp}$ and such that $\mu_h\in X_{\ast}(T)$ lies in the closed Weyl chamber determined by $(T_{\Qp},B)$. 

Indeed, then the canonical Galois action of $\Gamma(p)$ on $X_{\ast}(T)$ coincides with the naive Galois action, hence one has that $E(G,X)_{\wp}=E(T,h)_{\frak{p}}$, where $\frak{p}$ is the prime of $E(T,h)$ induced by $\iota_p$. Given this, the corollary then follows from the fact (\cite[Lem.2.2]{Bultel01}) that for a CM-abelian variety $A$ over $\Qb(\subset\C)$ with CM-type $(F,\Phi)$, where $F$ is a CM-algebra and $\Phi$ is a subset of $\Hom(F,\C))$ with $\Phi\sqcup\iota\circ \Phi=\Hom(F,\C)$, the reduction of $A$ at $\iota_p:\Qb\hra\Qpb$ is an ordinary abelian variety if and only if the prime of the reflex field $E=E(F,\Phi)$ induced by $\iota_p$ is of absolute height one; when $A$ is defined by a special Shimura datum $(T,h)$, the reflex field $E(F,\Phi)$ equals the reflex field $E(T,h)$ defined earlier. 

Now, to prove the claim, let $K_p=G_{\Zp}(\Z_p)$ for a reductive $\Zp$-group scheme model $G_{\Zp}$ of $G_{\Qp}$ and 
$S$ a maximal $\Qp$-split torus of $G_{\Qp}$ such that in the Bruhat-Tits building $\mathcal{B}(G,\Qp)$, the hyperspecial point $x$ fixed by $K_p$ is contained in the apartment associated with $S$. 
For such $S$, let $S_{\Zp}$ be the reductive $\Zp$-group scheme with generic fiber $S$ and such that $S_{\Zp}(\Zp)$ is the (unique) hyperspecial subgroup of $S(\Qp)$.
Then the proof of Prop. 
\ref{prop:existence_of_Levi_and_Hodge_cocharacter}, (3) shows that
$S_{\Zp}$ and $T'_{\Zp}:=\mathrm{Cent}_{G_{\Zp}}(S_{\Zp})$ are smooth, closed subschemes of $G_{\Zp}$, and thus $T':=T'_{\Zp}\otimes\Qp$ is unramified.
Moreover, since $G_{\Qp}$ is quasi-split, $T'$ is a maximal torus of $G_{\Qp}$, and by construction the unique hyperspecial subgroup of $T'(\Qp)$, i.e. $T'_{\Zp}(\Zp)$ is contained in $K_p$.
Since $G_{\Qp}$ is quasi-split, there exists a Borel subgroup $B'$ defined over $\Qp$ which contains $T'$. Let $\mu'$ be the cocharacter in $c(G,X)$ factoring through $T'$ and lying in the closed Weyl chamber determined by $(T',B')$. Let $(T,\{h\})$ be a special Shimura sub-datum produced from $(T',\mu')$ as in the proof of Theorem \ref{thm:existence_of_special_point_with_good_reduction_and_prescribed_given_isocrystal} using Theorem \ref{thm:existence_of_special_SD_with_given_F-isocrystal} (so, in particular $T(\Z_p)\subset K_p$). Then, as $(T_{\Qp},\mu_h)$ is conjugate to $(T',\mu')$ under $G(\Qp)$, the claim is now clear.
$\qed$

\begin{rem} \label{rem:mu-ordinary_locus}
(1) In fact, for any special Shimura sub-datum $(T,h)$ with the properties (i), (ii) in the proof, every special point of $\Sh_{K_p}(G,X)$ defined by it has $\mu$-ordinary reduction at $\iota_p$. So, we get a more direct proof of non-emptiness of $\mu$-ordinary locus (without invoking a Levi subgroup as in Proposition \ref{prop:existence_of_Levi_and_Hodge_cocharacter}, etc).

(2) This corollary was proved in the PEL-type cases of Lie type $A$ or $C$ by Wedhorn \cite{Wedhorn99}. His argument was to first establish non-emptiness of the $\mu$-ordinary locus, and then to check that the $\mu$-ordinary locus equals the (usual) ordinary locus when $E(G,X)_{\wp}=\Qp$ (in the PEL-type cases of Lie type $A$ or $C$). But to check the second statement, he relied on a case-by-case analysis. In contrast, we observe that in view of (1), the equality of the $\mu$-ordinary Newton point and the ordinary Newton point when $E(G,X)_{\wp}=\Qp$ is a consequence of the existence of a special Shimura sub-datum $(T,h)$ with the properties (i), (ii) in the proof (since the reduction of such special point attains simultaneously the two Newton points when $E(G,X)_{\wp}=\Qp$).
\end{rem}


\section{Construction of Kottwitz triple with prescribed $F$-isocrystal}

In this section, assuming that $G^{\mathrm{der}}$ is simply-connected, we construct a Kottwitz triple $(\gamma_0;(\gamma)_{l\neq p},\delta)$ (see below for its definition) such that the $F$-isocrystal with $G_{\Qp}$-structure attached to $\delta\in G(L)$ equals the given one. According to the Langlands-Rapoport conjecture \cite{LR87}, which is now confirmed by Kisin \cite{Kisin13}, such a triple corresponds to an $\Fpb$-valued point in the reduction of $\sS_{K_p}(G,X)$, if it satisfies certain additional condition, and under this correspondence the $F$-isocrystal of the abelian variety corresponding to the point is identified with the $F$-isocrystal defined by the $\sigma$-conjugacy class of $\delta$ in the triple. The Kottwitz triple that we will construct in this section will satisfy the additional condition just mentioned. Therefore, combined with the result of Kisin, we will get another proof of the non-emptiness of Newton strata for Shimura varieties with simply connected derived group.

\subsection{Kottwitz triples} 
Our two main references for the material covered here are \cite{LR87}, pp. 182-183 and \cite{Kottwitz88}, $\S2$.   
Let $(G,X)$ be a Shimura datum. Throughout this section, we assume that $G^{\mathrm{der}}$ is simply connected.
Let $p$ be a rational prime. We fix a prime $\wp$ of $E(G,X)$ that is unramified above $p$. Recall that $\Qb$ is given as the algebraic closure of $\Q$ in $\C$ and $E(G,X)$ is by definition a subfield of $\Qb\subset\C$ (thus the embedding $\Q\hra\C$ determines an infinite place $E\hra\C$) . For every place $v$ of $\Q$, we fix an embedding $\Qb\hookrightarrow\Qvb$, such that for $v=p,\infty$ this should induce the above places of $E(G,X)$. For each $r\in\N$, let $L_r$ denote the unramified extension of degree $r$ of $\Qp$ in $\Qpb$ (i.e. $L(\F_{p^r})$ in our previous notation), and denote by $L=L(\Fpb)$ the completion of $\Qp^{\nr}$, the maximal unramified extension of $\Qp$ in $\Qpb$. Let $\sigma$ be the absolute Frobenius on $L$; so, $L_r$ is the subfield of $L$ of elements fixed under $\sigma^r$.

\subsubsection{} \label{subsubsec:pre-Kottwitz_triple}
We consider a triple $(\gamma_0;\gamma=(\gamma_l)_{l\neq p},\delta)$, where
\begin{itemize} \addtolength{\itemsep}{-4pt}
\item[(i)] $\gamma_0$ is a semi-simple element of $G(\Q)$ that is elliptic in $G(\R)$, defined up to conjugacy in $G(\overline{\Q})$;
\item[(ii)] for $l\neq p$, $\gamma_l$ is a semi-simple element in $G(\Q_l)$, defined up to conjugacy in $G(\Q_l)$, which is conjugate to $\gamma_0$ in $G(\overline{\Q}_l)$;
\item[(iii)] $\delta$ is an element of $G(L_n)$ (for some $n$), defined up to $\sigma$-conjugacy in $G(L)$, such that the norm $\Nm_n\delta$ of $\delta$  is conjugated to $\gamma_0$ under $G(\overline{L})$, where $\Nm_n\delta:=\delta\cdot\sigma(\delta)\cdots\sigma^{n-1}(\delta)\in G(\Qp)$.
\end{itemize}

There are some conditions to impose on such triple. Let $I_0$ be the centralizer $\Cent_G(\gamma_0)$ of $\gamma_0$ in $G$. For every place $v$ of $\Q$, we will define an algebraic $\Q_v$-group $I(v)$ and an inner twisting $\psi_v:I_0\rightarrow I(v)$ over $\Qb_v$. 
First, for each finite place $v\neq p$ of $\Q$, let $I(v)$ be the centralizer of $\gamma_v$ in $G_{\Q_v}$. If one chooses $g_v\in G(\overline{\Q}_v)$ such that $g_v\gamma_0g_v^{-1}=\gamma_v$, then $\Int (g_v)$, restricted to $I_0$, gives an inner twisting $\psi_v:I_0\rightarrow I(v)$, well defined up to inner automorphism of $I_0$. Next, for $v=p$, define an algebraic $\Qp$-group $I(p)$ by
\begin{equation*}I(p):=\{x\in \mathrm{Res}_{L_n/\Qp}(G_{L_n})\ |\ x^{-1}\delta\theta(x)=\delta\},\end{equation*}
Here, $\theta$ is the $\Qp$-automorphism of $\mathrm{Res}_{L_n/\Qp}(G_{L_n})$ induced by the restriction of $\sigma$ to $L_n$; for more details, we refer to \cite{Kottwitz82}, where $\theta$ and $I(p)$ are denoted by $s$ (on p. 801) and $I_{s\delta}$ (on p. 802), respectively. Then, Lemma 5.8 of loc.cit. gives us an inner twisting $\psi_p:I_0\rightarrow I(p)$, canonical up to inner automorphisms of $I_0$; since $\gamma_0$ and $\Nm_{L_n/\Qp}(\delta)$ are $G(\Qpb)$-conjugated, 
$\gamma_0\in G(\Qp)$ lies in the stable conjugacy class of $\delta$ (cf. loc.cit. $\S5$). 
Finally, at the infinite place, choose an elliptic maximal torus $T$ of $G_{\R}$ containing $\gamma_0$.
When we pick $h\in X$ such that $h$ factors through $T_{\R}$, by definition of Shimura datum, $\Int (h(i))$ induces a Cartan involution on $I_0/Z(G)$, which we use to twist $I_0$ over $\R$, obtaining an inner twisting $\psi_{\infty}:I_0\rightarrow I(\infty)$ with $I(\infty)/Z(G)$ being anisotropic over $\R$.

Then we consider the following condition on the triple $(\gamma_0;\gamma=(\gamma_l)_{l\neq p},\delta)$ as above; let $\psi_v:I_0\rightarrow I(v)$ be inner twistings as just described. 
\begin{itemize} \addtolength{\itemsep}{-4pt}
\item[(iv)] There exists a triple $(I,\psi,(j_v))$ consisting of a $\Q$-group $I$, an inner twisting $\psi:I_0\rightarrow I$ and for each place $v$ of $\Q$, an isomorphism $j_v:I\rightarrow I(v)$ over $\Q_v$, unramified almost everywhere, such that $j_v\circ\psi$ and $\psi_v$ differ by an inner automorphism of $I_0$ over $\Qb_v$.
\end{itemize}

\begin{defn}  \label{defn:Kottwitz_triple}
A Kottwitz triple is a triple $(\gamma_0;(\gamma_l)_{l\neq p},\delta)$ satisfying (i) - (iv), and the following two conditions ($\ast(\delta)$), ($\ast(\gamma_0)$):
\begin{itemize} \addtolength{\itemsep}{-4pt}
\item[($\ast(\delta)$)] the image of (the $\sigma$-conjugacy class of) $\delta$ under the map $\kappa_{G_{\Q_p}}:B(G_{\Q_p})\rightarrow \pi_1(G_{\Q_p})_{\Gamma(p)}=X^{\ast}(Z(\widehat{G})^{\Gamma(p)})$ (Thm. \ref{thm;RR96-Thm.1.15}) is equal to $\mu^{\natural}$ (where $\mu^{\natural}$ was defined in (\ref{subsubsec:mu_natural})).
\item[($\ast(\gamma_0)$)] Let $H$ be the centralizer in $G_{\Qp}$ of the maximal $\Qp$-split torus in the center of $(Z_{G}(\gamma_0))_{\Qp}$. Then,
there exists $\mu\in X_{\ast}(H)$ which is defined over $L_n$ and lies in $c(G,X)$ (via the given embeddings $\Qb\hookrightarrow\C$, $\Qb\hookrightarrow\Qpb$), and such that 
\begin{equation*}\lambda_H(\gamma_0)=\Nm_{L_n/\Qp}\mu,\end{equation*}
where $\lambda_H$ is the homomorphism $\lambda_{H}:H(\Q_p)\rightarrow \pi_1(H)^{\Gamma(p)}=X^{\ast}(Z(\widehat{H}))^{\Gamma(p)}$ defined (for quasi-split $\Qp$-groups) in \cite{Kottwitz84b}, $\S$3. 
\end{itemize}
\end{defn}

\begin{rem}
1) Following Kisin \cite{Kisin13}, we will say that a Kottwitz triple $(\gamma_0;(\gamma_l)_{l\neq p},\delta)$ is \textit{of level $n$} if $\delta\in G(L_n)$. 

2) Two triples $(\gamma_0;(\gamma_l)_{l\neq p},\delta)$, $(\gamma_0';(\gamma_l')_{l\neq p},\delta')$ as in (\ref{subsubsec:pre-Kottwitz_triple}) with $\delta,\delta'\in G(L_n)$ for (iii)
are said to be \textit{equivalent}, if $\gamma_0$ is stably conjugate to $\gamma_0'$, $\gamma_l$ is conjugate to $\gamma_l'$ in $G(\Ql)$ for each $l\neq p,\infty$, and $\delta$ is $\sigma$-conjugate to $\delta'$ in $G(L_n)$. Then, for two such equivalent triples, one of them is a Kottwitz triple of level $n$ if and only if the other one is so (cf. see \cite[Lemma 5.17]{LR87} for the condition ($\ast(\gamma_0)$)).

3) As noted in \cite{Kottwitz88}, p.172, the condition (iv) holds if the Kottwitz invariant $\alpha(\gamma_0;\gamma,\delta)$ (to be recalled below) is trivial.

4) The two conditions ($\ast(\delta)$), ($\ast(\gamma_0)$) are from \cite{LR87}, p.182-183. Some authors do not include these (especially, ($\ast(\gamma_0)$)) in their definition of Kottwitz triple. This will be partly a matter of taste. It is, however, necessary to point out \cite{LR87} that a triple $(\gamma_0;(\gamma_l)_{l\neq p},\delta)$ satisfying (i)-(iv) can correspond to an $\Fpb$-valued point of the Shimura variety only if it further satisfies these two conditions, so from the viewpoint of the question of finding an $\Fpb$-valued point of given Shimura variety using such triple (as in the current work), it is natural to include these conditions ($\ast(\delta)$), ($\ast(\gamma_0)$) in the definition of Kottwitz triple. 

5) In view of our sign convention of $c(G,X)$ in (\ref{sssec:Hodge_cocharacter:mu_h}), our formulation (especially the sign) of the condition ($\ast(\delta)$)) is the same as that of Kottwitz appearing in \cite[p.165]{Kottwitz88}. 
\end{rem}

\subsection{Kottwitz invariant} \label{subsec:Kottwitz_invariant}
We recall the definition of the Kottwitz invariant, cf. \cite{Kottwitz88}, $\S2$. Note that the centralizer $I_0$ of $\gamma_0$ in $G$ is connected and reductive since $\gamma_0$ is semisimple and the derived group of $G$ is simply connected. 
There is a natural embedding $Z(\widehat{G})\rightarrow Z(\widehat{I}_0)$, and the exact sequence 
\begin{equation*}1\rightarrow Z(\widehat{G})\rightarrow Z(\widehat{I}_0)\rightarrow Z(\widehat{I}_0)/Z(\widehat{G}) \rightarrow 1\end{equation*}
induces a homomorphism (\cite[Cor.2.3]{Kottwitz84a})
\begin{equation*}\pi_0((Z(\widehat{I}_0)/Z(\widehat{G}) )^{\Gamma})\rightarrow H^1(\Q,Z(\widehat{G})),\end{equation*}
where $\Gamma:=\Gal(\Qb/\Q)$. Let $\ker^1(\Q,Z(\widehat{G}))$ be the kernel of the map
$H^1(\Q,Z(\widehat{G}))\rightarrow \prod_v H^1(\Q_v,Z(\widehat{G}))$, and
let $\mathfrak{K}(I_0/\Q)$ denote the subgroup of $\pi_0((Z(\widehat{I}_0)/Z(\widehat{G}))^{\Gamma})$ consisting of elements whose image in $H^1(\Q,Z(\widehat{G}))$ lies in $\ker^1(\Q,Z(\widehat{G}))$. Since $\gamma_0$ is elliptic, we have an identification
\begin{equation*}\mathfrak{K}(I_0/\Q)=\left(\cap_vZ(\widehat{I}_0)^{\Gamma(v)}Z(\widehat{G})\right)/Z(\widehat{G}),\end{equation*}
and this is known to be a finite group. 
The Kottwitz invariant $\alpha(\gamma_0;\gamma,\delta)$ is then a character of $\mathfrak{K}(I_0/\Q)$ (i.e. an element of $\Hom(\mathfrak{K}(I_0/\Q),\C^{\times})$), defined as a product, over all places $v$ of $\Q$, of local components, where the $v$-component is the restriction to $\cap_vZ(\widehat{I}_0)^{\Gamma(v)}Z(\widehat{G})$ of a character $\beta_v(\gamma_0;\gamma,\delta)$ on $Z(\widehat{I}_0)^{\Gamma(v)}Z(\widehat{G})$: 
\begin{equation*}\alpha(\gamma_0;\gamma,\delta)=\prod_v \beta_v(\gamma_0;\gamma,\delta)|_{\cap_v Z(\widehat{I}_0)^{\Gamma(v)}Z(\widehat{G})}.\end{equation*}
This product is trivial on $Z(\widehat{G})$. 
The character $\beta_v(\gamma_0;\gamma,\delta)$ itself is the unique extension of another character $\alpha_v(\gamma_0;\gamma,\delta)$ on $Z(\widehat{I}_0)^{\Gamma(v)}$ with the restriction of $\beta_v$ to $Z(\widehat{G})$ being 
\begin{equation*}
\beta_v|_{Z(\widehat{G})}=
\begin{cases}
\quad \mu_1 &\text{ if }\quad  v=\infty \\
\ -\mu_1 &\text{ if }\quad  v=p \\
\text{ trivial } &\text{ if }\quad  v\neq p,\infty.
\end{cases}
\end{equation*}
Here, $\mu_1$ is the \emph{negative} of the restriction to $Z(\widehat{I_0})$ of the Weyl group orbit in $X^{\ast}(\widehat{T})$ dual to the Weyl group orbit $X_{\ast}(T)\cap c(G,X)$.
\footnote{This $\mu_1$ is identical to the character denoted by the same symbol in \cite[p.165]{Kottwitz88}.}
For the definition of the character $\alpha_v(\gamma_0;\gamma,\delta)$, here we will be just contented with directing readers to the reference \cite{Kottwitz88}, $\S2$.

\subsection{Second proof of non-emptiness of Newton strata}

\begin{thm} \label{thm:existence_of_Kottwitz_triple_for_any_b_in_B(G,X)}
Let $(G,X)$ be a Shimura datum such that $G^{\mathrm{der}}$ is simply connected and $G_{\Q_p}$ is unramified. Fix a hyperspecial subgroup $K_p$ of $G(\Qp)$.
Then, 
for any $\overline{b}\in B(G_{\Q_p},\mathcal{C})$, there exist 
\begin{itemize}
\item[$\bullet$] a maximal $\Q$-torus $T$ of $G$ which is elliptic at $\R$ and unramified over $\Qp$, and such that the unique hyperspecial subgroup of $T(\Qp)$ is contained in a $G(\Qp)$-conjugate of $K_p$, 
\item[$\bullet$] a Kottwitz triple $(\gamma_0;\gamma,\delta)\in T(\Q)\times G(\A_f^p)\times T(\Qp^{\nr})$ with $\delta\in\overline{b}$ and having trivial Kottwitz invariant, and 
\item[$\bullet$] a cocharacter $\mu\in X_{\ast}(T)\cap c(G,X)$ such that the condition $(\ast(\gamma_0))$ of Definition \ref{defn:Kottwitz_triple} holds for $(T,\mu)$ and that $\delta$ is $\sigma$-conjugate to $\mu(p)$ in $T(L)$. 
\end{itemize}
The Newton point $\overline{\nu}_{T}(\overline{\delta})\in \mathcal{N}(T_{\Qp})=X_{\ast}(T_{\Qp})^{\Gamma(p)}_{\Q}$ of $\overline{\delta}\in B(T_{\Qp})$ equals
\begin{equation*}\nu:=\frac{1}{r}\Nm_{L_r/\Qp}(\mu),\end{equation*}
 where $L_r$ is a splitting field of $T_{\Qp}$, and $T_{\Qp}$ is an elliptic maximal torus of the centralizer $J=\Cent_{G_{\Qp}}(\nu)$. Moreover, for any rational prime $l\neq p$, we can find such object $(T,\mu,(\gamma_0;\gamma,\delta))$ such that $T_{\Ql}$ is elliptic in $G_{\Ql}$. 
\end{thm}

Here, the statement that the condition $(\ast(\gamma_0))$ of Definition \ref{defn:Kottwitz_triple} holds for $(T,\mu)$ means that $\lambda_{T_{\Qp}}(\gamma_0)=\Nm_{L_r/\Qp}\mu \in X_{\ast}(T)^{\Gamma(p)}$ (note that $\mu$ must be defined over $L_r$). 

\textsc{Proof.}
Our construction of $T$, $\mu$, and $(\gamma_0;\delta,(\gamma_v)_{v\neq p})\in T(\Q\times L)\times G(\A_f)$ proceeds in three steps: we will find $(T,\mu)$, the pair $(\gamma_0;\delta)\in T(\Q\times L)$, and the elements $(\gamma_v)_{v\neq p}\in G(\A_f^p)$ in these orders.

\textbf{Step 1.} According to Proposition \ref{prop:existence_of_Levi_and_Hodge_cocharacter} and the first step of the proof of Theorem \ref{thm:existence_of_special_SD_with_given_F-isocrystal}, there exists a maximal $\Q$-torus $T$ of $G$ which is elliptic at $\R$ and at any given prime $l\neq p$ and is unramified over $\Qp$, and a triple $(M,b\in M(L),\mu\in X_{\ast}(T_{\Qp})\cap c(G,X))$ satisfying the condition $(\dagger)$ combined with $T_{\Qp}$ (particularly, $\kappa_M(\overline{b})=[\mu]_{M}$).

\textbf{Step 2.} Choose $r\in\N$ such that $T_{\Q_p}$ splits over $L_r=L(\F_{p^r})$. Pick any $\mu\in X_{\ast}(T)$ such that $\kappa_M(\overline{b})=[\mu]_{M}$. In this step, we construct $(\gamma_0,\delta)\in T(\Q\times L)$, but for this construction itself, we will not need the condition $\mu\in c(G,X)$, yet; that property will be invoked later when we construct the prime-to-$p$ components $(\gamma_v)_{v\neq p}$ (especially, the $l$-component) and verify that the associated Kottwitz invariant is zero.

We claim that there exists $v\in \mathcal{T}(\cO_{L_r})$ (the unique maximal compact open subgroup of $T(L_r)$) such that the element
\begin{equation*}\delta:=\mu(p) v\ \in T(L_r),\end{equation*}
has the following three properties:
\begin{itemize} \addtolength{\itemsep}{-4pt} 
\item[(a)] its norm $\Nm_r(\delta)\in T(\Qp)$ belongs to $T(\Q)$; 
\item[(b)]$\delta$ and $\mu(p)$ are $\sigma$-conjugate in $T(L)$;
\item[(c)] $\mu(b)$ and $b$ are $\sigma$-conjugate in $M(L)$. 
\end{itemize}

It follows from (b) that the Newton point $\overline{\nu}_{T_{\Qp}}(\overline{\delta})$ of $\overline{\delta}\in B(T_{\Qp})$ equals $\frac{1}{r}\Nm_{L_r/\Qp}(\mu)\in \mathcal{N}(T_{\Qp})=X_{\ast}(T)_{\Q}^{\Gamma(p)}$.

We will see that as long as $v\in\mathcal{T}(\cO_{L_r})$, the property (b) always holds and 
the property (c) is equivalent to that $\kappa_M(\overline{b})=[\mu]_{M}$. Then, it is the property (a) for which we have to make a restricted choice of $v$; interestingly, $\mu$ does not play any significant role in such choice of $v$.

\begin{lem} \label{lem:unramified_tori_over_local_fields}
Let $T$ be a torus over $\Qp$.

(1) For every $s\in\N$, the image of $\mathcal{T}(\cO_{L_s})$ under the norm map $\Nm_s:T(L_s)\rightarrow T(\Q_p):t\mapsto t\cdot \sigma(t)\cdots \sigma^{s-1}(t)$ is open.

(2) If $T$ is unramified, the $\sigma$-conjugacy class of any $u\in\mathcal{T}(\cO_L)$ is trivial.
\end{lem}
\textsc{Proof of lemma.} 
(1) The composition of the natural maps of algebraic tori over $\Q_p$
\begin{equation*}T_{\Q_p}\hookrightarrow \mathrm{Res}_{L_s/\Q_p}(T_{L_s})\stackrel{N_s}{\rightarrow}T_{\Q_p}\end{equation*} is an isogeny, so the image of the maximal compact subgroup of $T(\Q_p)$ is open in $T(\Q_p)$ 
(\cite[Ch.3, Cor.1]{PR94}).

(2) We recall (Theorem \ref{thm;RR96-Thm.1.15} (i)) the natural transformation $\kappa:B(\cdot)\rightarrow X_{\ast}(\cdot)_{\Gamma(p)}$ of set-valued functors on the category of connected reductive groups over $\Q_p$ which gives isomorphisms for tori. 
For an unramified $\Q_p$-torus $T$, this isomorphism $\kappa_{T}:B(T)\isom X_{\ast}(T)_{\Gamma(p)}$ is obtained as follows \cite[2.9]{Kottwitz85}: tensor the normalized valuation $L^{\times}\rightarrow\Z$ with $X_{\ast}(T)$ to get a canonical surjection $T(L)\rightarrow X_{\ast}(T)$. The isomorphism $\kappa_{T}$ is obtained by applying the functor $H^1(\langle\sigma\rangle,\cdot)$ to this surjection; on $\Gm$, we have $\kappa_{\Gm}(\overline{b})=\mathrm{ord}_L(b)$ for $b\in L^{\times}$.  Now, it is obvious that any $u\in \mathcal{T}(\cO_L)=\Hom(X^{\ast}(T),\cO_L^{\times})=X_{\ast}(T)\otimes_{\Z}\cO_L^{\times}$ is a linear combination of elements of $X_{\ast}(T)$ with coefficients in $\cO_L^{\times}$, hence maps to zero under the surjection $T(L)\rightarrow X_{\ast}(T)$. $\qed$

We now prove the claims. 
Since $T_{\Q_p}$ is unramified over $\Q_p$, $T(\Q)$ is dense in $T(\Q_p)$ (\cite[Lemma 4.10]{Milne92} or \cite[Theorem in $\S$11.5]{Voskresenskii98}). So, by Lemma \ref{lem:unramified_tori_over_local_fields} (1), there exists $v\in \mathcal{T}(\cO_{L_r})$ such that
\begin{equation*}\Nm_r(\mu(p))\cdot \Nm_r(v)\in T(\Q).\end{equation*} 
Namely, for such $v$, the element
\begin{equation*}\delta:=\mu(p)\cdot v\in T(L_r)\end{equation*} 
satisfies the property (a). Also, by Lemma \ref{lem:unramified_tori_over_local_fields} (2), 
$\delta$ is $\sigma$-conjugated to $\mu(p)$ in $T(L)$, which proves (b).

\textbf{Step 3.}
Finally, we assume that our $\mu\in X_{\ast}(T)$ further satisfies that $\mu\in c(G,X)$. Fix $h\in X$ which factors through $T_{\R}$ and let $\mu_h$ be the associated Hodge cocharacter. We will find an element $\gamma_l\in T(\Ql)$ which is conjugate to $\gamma_0$ in $G(\Qlb)$, characterized by certain cohomological property. In general, if an element $\gamma_l\in G(\Ql)$ is conjugate to $\gamma_0$ in $G(\Qlb)$, say, $\gamma_l=g_l\gamma_0g_l^{-1}\ (g_l\in G(\Qlb))$, one has that  $\tau(g_l\gamma_0g_l^{-1})=g_l\gamma_0g_l^{-1}$, i.e. $g_l^{-1}\tau(g_l)\in \Cent_{G_{\Ql}}(\gamma_0)(\Qlb)$, for all $\tau\in\Gal(\overline{\Q}_l/\Q_l)$.
In this way, one obtains a bijection between the set of $\Ql$-rational conjugacy classes of elements $\gamma_l\in G(\Ql)$ which are conjugate to $\gamma_0$ in $G(\Qlb)$ and the set
\begin{equation*}\mathrm{Ker}[H^1(\Q_l,I_0)\ra H^1(\Q_l,G)].\end{equation*}

On the other hand, we have a commutative diagram (for any finite place $v$) \cite[Prop. 6.4]{Kottwitz84a}
\begin{equation} \label{eq:alpha_l}
\xymatrix{
H^1(\Q_v,T) \ar[r] \ar[d]^{\cong} & H^1(\Q_v,I_0) \ar[r] \ar[d]^{\cong} & H^1(\Q_v,G) \ar[d]^{\cong} \\
\pi_0(Z(\widehat{T})^{\Gamma(v)})^D \ar[r] & \pi_0(Z(\widehat{I_0})^{\Gamma(v)})^D \ar[r] &\pi_0(Z(\widehat{G})^{\Gamma(v)})^D,}
\end{equation}
and for any connected reductive $\Q_v$-group $H$,
$\pi_0(Z(\widehat{H})^{\Gamma(v)})^D$ is the subgroup of $X^{\ast}(Z(\widehat{H})^{\Gamma(v)})=\Hom(Z(\widehat{H})^{\Gamma(v)},\C^{\times})$ of torsion elements (\cite[1.14]{RR96}).
Hence, if the image $\alpha_l$ of the cocharacter $-(\mu_h+\mu)\in X_{\ast}(T)$ in $X^{\ast}(Z(\widehat{I_0})^{\Gamma(l)})=\pi_1(I_0)_{\Gamma(l)}$ (via the map $X_{\ast}(T)\ra \pi_1(I_0)$) is a torsion element and furthermore maps to zero in $X^{\ast}(Z(\widehat{G})^{\Gamma(l)})=\pi_1(G)_{\Gamma(l)}$, it will belong to $\mathrm{Ker}[H^1(\Q_l,I_0)\ra H^1(\Q_l,G)]$ and accordingly we will get an element $\gamma_l\in G(\Ql)$ that is $\Qlb$-conjugate to $\gamma_0$. 
But, the cocharacter $-(\mu_h+\mu)\in X_{\ast}(T)$ in fact belongs to the subgroup $X_{\ast}(T^{\mathrm{der}})$, where $T^{\mathrm{der}}:=(T\cap G^{\mathrm{der}})^0(=T\cap G^{\mathrm{der}})$,  
and $X_{\ast}(T^{\mathrm{der}})_{\Gamma(l)}=X^{\ast}(\widehat{T^{\mathrm{der}}}^{\Gamma(l)})$ is a torsion group since $T^{\mathrm{der}}_{\Ql}$ is an anisotropic torus by our choice of $T$. Moreover, as $\mu$ is an image of $\mu_h^{-1}$ under the Weyl group action, $-(\mu_h+\mu)$ lies in the coroot lattice $\sum_{\alpha\in \Delta(T,G)}\alpha^{\vee}$, hence maps to zero in $\pi_1(G)$.

Now, we claim that when we put 
\begin{equation*}\gamma_0:=\Nm_r(\delta),\quad \gamma_v:=\gamma_0\ (v\neq p,l),\quad \gamma_l=\alpha_l\end{equation*}
the triple $(\gamma_0;(\gamma_v)_{v\neq p},\delta)$ forms a Kottwitz triple. 

First, we check the two conditions $(\ast(\delta))$, $(\ast(\gamma_0))$ of Definition \ref{defn:Kottwitz_triple}. 
For the condition $(\ast(\delta))$, we recall the identity (from (\ref{subsubsec:mu_natural}))
\begin{equation*}\mu^{\natural}=[\mu]_{G_{\Qp}}\ \in \pi_1(G)_{\Gamma(p)},\end{equation*}
as $\mu\in X_{\ast}(T)\cap c(G,X)$. 
As $\delta$ belongs to $T(L)$ and $T_{\Q_p}$ is also a maximal torus of $M$, the functoriality of the natural transformation $\kappa:B(\cdot)\rightarrow\pi_1(\cdot)_{\Gamma(p)}$ allows us to replace, in our verification, $G_{\Q_p}$ by $M$, in which case the condition $(\ast(\delta))$ is our assumption $\kappa_M(\overline{\delta})=[\mu]_{M}$.

For the condition $(\ast(\gamma_0))$, recall that we introduced the centralizer $H$ of the maximal split $\Qp$-subtorus of the center of $(Z_{G}(\gamma_0))_{\Qp}$; since $\gamma_0\in T(\Q)$ and $T_{\Q_p}$ is an unramified maximal torus of $G_{\Qp}$, $T_{\Qp}$ is also an \emph{unramified} maximal torus of $H$. Hence, by the functoriality of the natural transformation $\lambda_H:H(\Qp)\rightarrow \pi_1(H)^{\Gamma(p)}$ (\cite{Kottwitz84b}, Lemma 3.3), we may replace $H$ by $T_{\Qp}$, in which case $\lambda_T:T(\Qp)\rightarrow X_{\ast}(T)^{\Gamma(p)}=X^{\ast}(\widehat{T})^{\Gamma(p)}$ is obtained by tensoring the normalized valuation $(\Qp^{\nr})^{\times}\rightarrow\Z$ with $X_{\ast}(T)$ and taking invariants under $\Gamma(p)$, (\cite[$\S$3]{Kottwitz84b}). 
It follows from this definition that $\lambda(\mathcal{T}(\Z_p))=0$ and thus
\begin{equation*}\lambda_{T_{\Qp}}(\gamma_0)=\lambda_{T_{\Qp}}(\Nm_r(\delta))=\lambda_{T_{\Qp}}(\Nm_r(\mu)(p))
=\Nm_r(\mu),\end{equation*}
as required.

It remains to show the vanishing of the Kottwitz invariant $\alpha(\gamma_0;\gamma,\delta)$; recall from Subsection \ref{subsec:Kottwitz_invariant} that it is a product, over all places $v$ of $\Q$, of the restriction to $\cap_vZ(\widehat{I}_0)^{\Gamma(v)}Z(\widehat{G})$ of characters $\beta_v(\gamma_0;\gamma,\delta)$ on $Z(\widehat{I}_0)^{\Gamma(v)}Z(\widehat{G})$, which itself is a unique extension of another character $\alpha_v(\gamma_0;\gamma,\delta)$ on  $Z(\widehat{I}_0)^{\Gamma(v)}$. 
For $v\neq  l,p$, since $\gamma_0=\gamma_v$ it is obvious from the definition (cf. \cite[$\S2$]{Kottwitz88}) that $\alpha_v(\gamma_0;\gamma,\delta)$ and thus $\beta_v(\gamma_0;\gamma,\delta)$ as well are trivial. Also, by definition (\cite[$\S$2]{Kottwitz88}), $\alpha_l(\gamma_0;\gamma,\delta)$ is 
the cohomology class $[g_l^{-1}\tau(g_l)]$, regarded as an element of $\pi_0(Z(\widehat{I_0})^{\Gamma(v)})^D$ via the vertical isomorphism $H^1(\Q_v,I_0)\cong \pi_0(Z(\widehat{I_0})^{\Gamma(v)})^D$ of (\ref{eq:alpha_l}), where $\gamma_l=g_l\gamma_0 g_l^{-1}$ for $g_l\in G(\Qlb)$ with $g_l^{-1}\tau(g_l)\in I_0(\Qlb)$, hence is nothing but $-[\mu_h+\mu]_{{I_0}_{\Ql}}$ by the very definition of $\gamma_l$. 
Working out the definitions (first and second paragraphs on p.167 in $\S$2 of \cite{Kottwitz88}) for $v=p,\infty$ in a similar way, we get that 
$\alpha_v(\gamma_0;\gamma,\delta)\in X^{\ast}(Z(\widehat{I_0})^{\Gamma(v)})= \pi_1(I_0)_{\Gamma(v)}$ is given by 
\begin{equation*}\alpha_v(\gamma_0;\gamma,\delta)=
\begin{cases}
\kappa_{(I_0)_{\Qp}}(\overline{\delta}) & \text{ if }v=p \\
\quad [\mu_h]_{{I_0}_{\R}} & \text{ if }v=\infty,\\
-[\mu_h+\mu]_{{I_0}_{\Ql}} &\text{ if }v=l 
\end{cases}
\end{equation*} 
See the remark at the end of Subsection \ref{subsec:Kottwitz_invariant} for our sign convention. 

Recall that for $v=p$, we have (as a consequence of $\delta\in\mu(p)\cdot\mathcal{T}(\cO_{L_r})$) that $\kappa_{T_{\Qp}}(\overline{\delta})=[\mu]_{T_{\Qp}}$. 
Hence, when we fix embeddings $\widehat{T}\hra \widehat{I_0}\hra \widehat{G}$, we see that $\beta_{p}(\gamma_0;\gamma,\delta)$ 
is equal to the restriction of the character $\mu\in X_{\ast}(T)=X^{\ast}(\widehat{T})$ of $\widehat{T}$ to 
\begin{equation*}Z(\widehat{I}_0)^{\Gamma(p)} Z(\widehat{G})\hra\widehat{T}^{\Gamma(p)} Z(\widehat{G})\hra \widehat{T}\end{equation*} and similarly $\beta_{\infty}(\gamma_0;\gamma,\delta)$ is the restriction of $\mu_h^{-1}\in X_{\ast}(T)=X^{\ast}(\widehat{T})$ to $Z(\widehat{I}_0)^{\Gamma(\infty)}Z(\widehat{G})\hra \widehat{T}$. 
It follows that the product of the restrictions of $\beta_{\infty}(\gamma_0;\gamma,\delta)$ and $\beta_{p}(\gamma_0;\gamma,\delta)$ to 
$\bigcap_{v=\infty,p} Z(\widehat{T})^{\Gamma(v)} Z(\widehat{G})$ is the pull-back of the character $\mu+\mu_h\in X_{\ast}(T)=X^{\ast}(\widehat{T})$ to 
\begin{equation*}\bigcap_{v=\infty,p} Z(\widehat{I}_0)^{\Gamma(v)} Z(\widehat{G})\hra\bigcap_{v=\infty,p} \widehat{T}^{\Gamma(v)} Z(\widehat{G})\hra \widehat{T}.\end{equation*} 
As clearly $\beta_l(\gamma_0;\gamma,\delta)$ is the restriction of $-(\mu_h+\mu)\in X_{\ast}(T)=X^{\ast}(\widehat{T})$ to $Z(\widehat{I}_0)^{\Gamma(l)}Z(\widehat{G})\hra \widehat{T}$, we see that the Kottwitz invariant $\alpha(\gamma_0;\gamma,\delta)=\prod_v \beta_v(\gamma_0;\gamma,\delta)|_{\cap_v Z(\widehat{I}_0)^{\Gamma(v)}Z(\widehat{G})}$ vanishes.

It remains to show that $T_{\Qp}$ is elliptic in $J$. According to Proposition \ref{prop:existence_of_Levi_and_Hodge_cocharacter}, there exist a Levi $\Qp$-subgroup $M$ of $G_{\Qp}$ and a representative $b$ of $\overline{b}$ in $M(L)$ with $\nu_b\in \Hom_{\Qp}(\mathbb{D},G_{\Qp})$ such that $M=\Cent_{G_{\Qp}}(\nu_{b})$. Then, since $\overline{b}=\overline{\delta}$ in $B(G_{\Qp})$, it follows that $\nu_{b}=g \nu g^{-1}$ for some $g\in G(L)$ so that $M$ is an inner form of $J$ (indeed, $M=\Int(g)(J)$. On the other hand, since $\nu$, $\nu_{b}$ are defined over $\Qp$, $g^{-1}\sigma(g)\in \Cent_{G(L)}(\nu)=J(L)$). In particular, $\Int(g)$ induces a $\Qp$-isomorphism from $Z(J)$ to $Z(M)$. 
Now, let $T$ be as in Theorem \ref{thm:existence_of_special_SD_with_given_F-isocrystal} where the Levi subgroup $M'$ in the assumption $(\dagger)$ is $M$ just chosen, so that a $G(\Qp)$-conjugate of $T_{\Qp}$ is an elliptic maximal torus of $M$. But, as $T_{\Qp}\subset J$, $Z(J)\subset T_{\Qp}$. Then, $T_{\Qp}$ contains a subtorus $Z' \approx Z(M)$ such that $T_{\Qp}/Z'$ is anisotropic. But, since $T_{\Qp}\subset J$, $Z(J)\subset T_{\Qp}$. As $Z'\approx Z(J)$, thus $T_{\Qp}/ Z(J)$ is also anisotropic, i.e. $T_{\Qp}$ is also an elliptic maximal torus of $J$  (the maximal split subtorus $S$ of $T_{\Qp}$ is contained in $Z'$, so $Z(J)\approx Z'$ contains a split subtorus of $T_{\Qp}$ of same rank, which then must be $S$).  This finishes our proof of the theorem. $\qed$

Now, we can give another proof of the non-emptiness of Newton strata of Shimura varieties of Hodge type, using the resolution of the Langlands-Rapoport conjecture by Kisin \cite{Kisin13}.

\begin{cor} \label{for:second_proof_of_non-emptiness}
Let $(G,X)$ be a Shimura datum of Hodge type such that $G^{\mathrm{der}}$ is simply-connected. Let $p$ be a rational prime such that $G_{\Qp}$ is unramified. Fix a prime $\wp$ of $E(G,X)$ unramified over $p$ 
and a hyperspecial subgroup $K_p$ of $G(\Qp)$. Let $\Fpb$ be an algebraic closure of the residue field $\kappa(\wp)$ of $\wp$.

Then, for every $\overline{b}\in B(G_{\Qp},X)$, there exists a point $x\in \Sh_{K_p}(G,X)(\Fpb)$ whose associated $F$-isocrystal is $\overline{b}$.
\end{cor}

\textsc{Proof.}
This follows from Theorem \ref{thm:existence_of_Kottwitz_triple_for_any_b_in_B(G,X)}
and \cite{Kisin13}, Theorem 0.3. The statement of this theorem of Kisin is formulated in terms of the notion of Galois gerbs and admissible morphisms, but it is known (cf. \cite[$\S$5]{Milne92}, \cite{Milne09})
that this implies that for any Kottwitz triple with trivial Kottwitz invariant, there exists a point in $\Sh_{K_p}(G,X)(\Fpb)$ which has that Kottwitz triple as the associated Kottwitz triple. $\qed$

\appendix

\section{A lemma on algebraic groups over $p$-adic fields}

The goal of this appendix is to give a proof of the following fact, which is perhaps well-known (for example, it is mentioned without proof in the last paragraph on p. 172 of \cite{LR87}). But, since we could not find its proof in literatures, we present a proof here.

\begin{lem} \label{lem:existence_of_aniostropic_and_unramified_torus}
Let $H$ be an unramified semi-simple group defined over a $p$-adic field $F$ with ring of integers $\cO_F$.
Then $H$ admits a maximal $F$-torus which is anisotropic and unramified.
If $K$ is a hyperspecial subgroup of $H(F)$, there exists a such maximal $F$-torus $T'$ (i.e. aniostropic and unramified) such that the (unique) hyperspecial subgroup of $T'(F)$ is contained in $K$. 
\end{lem}

\textsc{Proof.} 
We first review some constructions in the general theory of reductive group schemes over $S=\mathrm{Spec}(\cO_F)$.

Let us use $H$ (by abuse of notation) again to denote the reductive group scheme over $S$ which is an integral model of the given $F$-group in question and such that $H(\cO_F)=K$. 
We fix a maximal $S$-torus $T$ of $H$; by definition (\cite[3.2.1]{Conrad11}), it is a $S$-torus $T\subset H$ whose every geometric fiber $T_{\overline{s}}$ is a maximal torus of $G_{\overline{s}}$ in the usual sense.
Let $N=N_H(T)$ be the normalizer scheme of $T$ in $H$ (\cite[2.1.1]{Conrad11}): it is a (finitely presented) closed subscheme of $G$, smooth over $S$ (\cite[2.1.2]{Conrad11}). Also, let $W:=N_H(T)/Z_H(T)$ be the Weyl group scheme of $T$, where $Z_{H}(T)$ is the centralizer scheme of $T$ in $H$ (\cite[2.2]{Conrad11}). It is a finite \'etale scheme over $\cO$ \cite[3.2.8]{Conrad11}). 
Let $X$ be the scheme of maximal tori of $H$: this is the scheme representing the functor on $S$-schemes
\begin{equation*}\underline{\mathrm{Tor}}_{H/S}: S'\mapsto \{\text{ maximal tori in }H_{S'}\}\end{equation*}
(which is a sheaf of sets in the fppf topology, \cite[3.2.6]{Conrad11}), and is a smooth affine scheme (\cite[XII, 5.4]{SGA3}).
Moreover, this also represents the quotient sheaf $H/N$ (\cite{Conrad11}, Thm. 2.3.1 and proof of Thm. 3.2.6). Now, to prove the lemma, it is enough to show the existence of an anisotropic maximal torus of the special fibre $H_{\kappa}=H\otimes_{\cO_F}\kappa$, where $\kappa$ is the residue field of $\cO_F$. Indeed, since $\underline{\mathrm{Tor}}_{H/S}$ is smooth over the henselian base $\cO_F$, then there exists a maximal torus $\mathcal{T}'$ over $\cO_F$ whose special fibre is anisotropic. But, since $\mathcal{T}'$ is a torus over $\cO_F$, the $\pi_1(S,\overline{\eta})$-module attached to the \'etale sheaf $X^{\ast}(\mathcal{T}')=\underline{\Hom}_{S_{\et}}(\mathcal{T}',\Gm)$ is unramified ($\overline{\eta}$ is a geometric generic point of $S$), 
hence the $F$-rank of $\mathcal{T}'_F$ equals the $\kappa$-rank of $\mathcal{T}'_{\kappa}$. Also, since $T'$ is a subscheme of $H$, the unique hyperspecial subgroup $\mathcal{T}'(\cO_F)$ of $\mathcal{T}(F)$ is contained in $K=H(\cO_F)$.

For the closed point $s$ of $S$ and a (fixed) algebraic closure $\overline{\kappa}=\overline{\kappa(s)}$ of its residue field $\kappa=\kappa(s)$, there exists a canonical diagram
\begin{equation*}\xymatrix{
X(\kappa) \ar[r]^{\delta} \ar[rd]_{\varphi:=j\circ \delta} & H^1(\kappa,N_{\kappa}) \ar[r] \ar[d]^{j} & H^1(\kappa,H_{\kappa}) \\
& H^1(\kappa,W_{\kappa}) \ar[r] & H^1(\kappa,\mathrm{Aut} T_{\kappa})
}\end{equation*}
Here, the upper sequence of pointed sets is exact and $\delta$ is defined as follows: for a maximal $\kappa$-torus $T'$ of $H_{\kappa}$, the image of $T'\in X(\kappa)$ under $\delta$ is the cohomology class of the cocycle \begin{equation*}\sigma\mapsto a_{\sigma}:=g^{-1}\cdot \sigma(g)\in N(\overline{\kappa}),\end{equation*} where $g\in H(\overline{\kappa})$ satisfies $T'=gTg^{-1}$. The lower horizontal map is induced from the injection $W\hra \mathrm{Aut} T$. 
On the other hand, the pointed set $H^1(\kappa,\mathrm{Aut} T)$ classifies the $\kappa$-isomorphism classes of $\kappa$-forms of $T_{\kappa}$, and for any $T'\in X(\kappa)$, its corresponding class in $H^1(\kappa,\mathrm{Aut} T_{\kappa})$ is the image of $[\varphi(T')]$ under the natural map
$H^1(\kappa,W_{\kappa})\rightarrow H^1(\kappa,\mathrm{Aut} T_{\kappa})$, namely, $T'$ is the twist of $T$ by (a cocycle in $W(\overline{\kappa})$ representing) $\varphi(T')$. 
Therefore, to prove the lemma, it will suffice to find a cohomology class $\xi\in H^1(\kappa,W_{\kappa})$ which lies in the image of $\varphi$ and such that the twist ${}_{\xi}T$ of $T_{\kappa}$ by $\xi$ is anisotropic. 
But, since $H^1(\kappa,H_{\kappa})=H^2(\kappa,T_{\kappa})=\{1\}$ (Lang's theorem),
$\varphi$ is always surjective, hence we only need to find a cohomology class $\xi\in H^1(\kappa,W_{\kappa})$ such that the twist ${}_{\xi}T_{\kappa}$ is anisotropic. 

This is also well-known. More explicitly, recalling the fact that every reductive group over a finite field is quasi-split, choose a Borel pair $(T_1,B_1)$ defined over $\kappa$, and let $\Delta$ be the associated set of simple roots of $(H_{\kappa},T_1)$, and $N(T_1)$ the normalizer of $T_1$; thus, the Frobenius automorphism $\sigma$ of $\kappa$ acts on $X^{\ast}(T_1)$, leaving $\Delta$ stable.
Then, what one needs is an element $w\in N(T_1)$ such that $w \sigma\in H(\overline{\kappa})\rtimes\Gal(\overline{\kappa}/\kappa)$ does not have a fixed vector in $X_{\ast}(T)_{\R}$. 
In this case, it is shown in \cite[Lemma 7.4]{Springer74} that if $\{\alpha_1,\cdots,\alpha_m\}$ is a set of representatives of the $\Gal(\overline{\kappa}/\kappa)$-orbits in $\Delta$, 
the product $\omega$ (in any order) of the corresponding simple reflections $r_i\ (i=1,\cdots,m)$, called ``twisted Coxeter element", has such property. 
$\qed$

Center for Geometry and Physics, Institute for Basic Science (IBS), 

Pohang 790-784, Republic of 
Korea 

\textit{Email:} dulee@ibs.re.kr


\begin{thebibliography}{99}

\bibitem[BO83]{BO83}
P. Berthelot, A. Ogus, 
\newblock F-isocrystals and de Rham cohomology. I. 
\newblock Invent. Math. 72 (1983), no. 2, 159-199.


\bibitem[Bla94]{Blasius94}
D. Blasius, 
\newblock A p-adic property of Hodge classes on abelian varieties. 
\newblock Motives (Seattle, WA, 1991), 293-308, Proc. Sympos. Pure Math., 55, Part 2, Amer. Math. Soc., Providence, RI, 1994.


\bibitem[BT84]{BT84}
F. Bruhat, J. Tits.
\newblock Groupes r\'eductifs sur un corps local. II. Sch\'emas en groupes. Existence d'une donn\'ee radicielle valu\'ee. 
\newblock Inst. Hautes \'Etudes Sci. Publ. Math. No. 60 (1984), 197-376.


\bibitem[Bou05]{Bourbaki05}
N. Bourbaki.
\newblock Lie groups and Lie algebras. Chapters 4-6, Chapters 7-9.
\newblock 
Elements of Mathematics (Berlin). Springer-Verlag, Berlin, 2002, 
2005.


\bibitem[Bul01]{Bultel01}
O. B\"ultel.
\newblock Density of the ordinary locus.
\newblock Bull. London Math. Soc.  33  (2001),  no. 2, 149-156.


\bibitem[Cha00]{Chai00}
C.-L. Chai. 
\newblock Newton polygons as lattice points. 
\newblock Amer. J. Math. 122 (2000), no. 5, 967-990.


\bibitem[Con11]{Conrad11}
B. Conrad.
\newblock Reductive group schemes.
\newblock Notes for ``SGA3 summer school", 2011, available at {\url{http://math.stanford.edu/~conrad/papers/luminysga3smf.pdf}}


\bibitem[Del71]{Deligne71}
P. Deligne.
\newblock Travaux de Shimura.
\newblock Seminaire Bourbaki, 23eme annee (1970/71), Exp. No. 389, 123-165.
Lecture Notes in Math., Vol. 244, Springer, Berlin, 1971.


\bibitem[Del77]{Deligne77}
P. Deligne.
\newblock Vari\'{e}t\'{e}s de Shimura: interpr\'{e}tation modulaire,
et techniques de construction de mod\`{e}les canoniques.
\newblock Automorphic forms, representations and $L$-functions
(Proc. Sympos. Pure Math., Oregon State Univ., Corvallis, Ore., 1977), Part 2, 247-289, Proc. Sympos. Pure Math., XXXIII, Amer. Math. Soc., Providence, R.I., 1979.


\bibitem[DMOS82]{DMOS82}
P. Deligne, J.S. Milne, A. Ogus, K-Y. Shih.
\newblock Hodge cycles, motives, and Shimura varieties. 
\newblock Lecture Notes in Math., 900. Springer-Verlag, Berlin-New York, 1982.


\bibitem[SGA3]{SGA3}
M. Demazure, A. Grothendieck. 
\newblock Sch\'emas en groupes I, II, III. 
\newblock Lecture Notes in Math., 151, 152, 153, Springer-Verlag, New York, 1970.


\bibitem[Far04]{Fargues04}
L. Fargues. 
\newblock Cohomologie des espaces de modules de groupes p-divisibles et correspondances de Langlands locales. 
\newblock Vari\'et\'es de Shimura, espaces de Rapoport-Zink et correspondances de Langlands locales. Ast\'erisque No. 291 (2004), 1-199.


\bibitem[FM87]{FM87}
J.-M. Fontaine, W. Messing. 
\newblock $p$-adic periods and $p$-adic \'etale cohomology. 
\newblock Current trends in arithmetical algebraic geometry (Arcata, Calif., 1985), 179-207, Contemp. Math., 67, Amer. Math. Soc., Providence, RI, 1987. 


\bibitem[Kis09]{Kisin09}
M. Kisin.
\newblock Integral canonical models of Shimura varieties. 
\newblock J. Th\'eor. Nombres Bordeaux 21 (2009), no. 2, 301-312.


\bibitem[Kis10]{Kisin10}
M. Kisin.
\newblock Integral models for Shimura varieties of abelian type.
\newblock  J. Amer. Math. Soc. 23 (2010), no. 4, 967-1012.


\bibitem[Kis13]{Kisin13}
M. Kisin.
\newblock Mod $p$-points on Shimura varieties of abelian type.
\newblock  preprint available on {\url{http://www.math.harvard.edu/~kisin/}}




\bibitem[Kos14]{Koskivirta14}
J.-S. Koskivirta
\newblock Sections of the Hodge bundle over Ekedahl-Oort strata of Shimura varieties of Hodge type.
\newblock {\url{http://arxiv.org/abs/1410.1317}}.


\bibitem[Kot82]{Kottwitz82}
R. Kottwitz. 
\newblock Rational conjugacy classes in reductive groups. 
\newblock Duke Math. J. 49 (1982), no. 4, 785-806. 


\bibitem[Kot84a]{Kottwitz84a}
R. Kottwitz. 
\newblock Stable Trace Formula: Cuspidal Tempered Terms. 
\newblock Duke Math. J. 51 (1984), 611-650.


\bibitem[Kot84b]{Kottwitz84b}
R. Kottwitz. 
\newblock Shimura Varieties and Twisted Orbital Integrals. 
\newblock Math. Ann. 269 (1984), 287-300.


\bibitem[Kot85]{Kottwitz85}
R. Kottwitz.
\newblock Isocrystals with additional structure.
\newblock Compositio Math. 56 (1985), no. 2, 201-220.


\bibitem[Kot86]{Kottwitz86}
R. Kottwitz. 
\newblock Stable Trace Formula: Elliptic Singular Terms.
\newblock Math. Ann. 275 (1986), 365-399.


\bibitem[Kot88]{Kottwitz88}
R. Kottwitz.
\newblock Shimura varieties and $\lambda$-adic representations. 
\newblock Automorphic forms, Shimura varieties, and L-functions, Vol. I (Ann Arbor, MI, 1988), 161-209, Perspect. Math., 10, Academic Press, Boston, MA, 1990.


\bibitem[Kot92]{Kottwitz92}
R. Kottwitz.
\newblock Points on some Shimura varieties over finite fields.
\newblock  J. Amer. Math. Soc. 5 (1992), no. 2, 373-444.


\bibitem[Kot97]{Kottwitz97}
R. Kottwitz.
\newblock Isocrystals with additional structure. II.
\newblock Compositio Math. 109 (1997), no. 3, 255-339.


\bibitem[KR03]{KR03}
R. Kottwitz, M. Rapoport. 
\newblock On the existence of F-crystals. 
\newblock Comment. Math. Helv. 78 (2003), no. 1, 153-184. 


\bibitem[Kot03]{Kottwitz03}
R. Kottwitz. 
\newblock On the Hodge-Newton decomposition for split groups. 
\newblock Int. Math. Res. Not. 2003, no. 26, 1433-1447.


\bibitem[Kre13]{Kret13}
A. Kret.
\newblock The trace formula and the existence of PEL type Abelian varieties modulo $p$.
\newblock {\url{http://arxiv.org/abs/1209.0264}}.


\bibitem[Lan83]{Langlands83}
R. P. Langlands. 
\newblock Les d\'ebuts d'une formule des traces stable. 
\newblock Publ. Math. de l'Univ. Paris VII 13, Paris 1983.


\bibitem[LR87]{LR87}
R. P. Langlands, M. Rapoport. 
\newblock Shimuravariet\"aten und Gerben. 
\newblock J. Reine Angew. Math. 378 (1987), 113-220.




\bibitem[Mil92]{Milne92}
J.S. Milne. 
\newblock The points on a Shimura variety modulo a prime of good reduction.  
\newblock The zeta functions of Picard modular surfaces, Univ. Mont\'eal, Montr\'eal, QC, 1992.


\bibitem[Mil94]{Milne94}
J.S. Milne.
\newblock Shimura varieties and motives.
\newblock Motives (Seattle, WA, 1991),  447--523, Proc. Sympos. Pure Math.,
55, Part 2, Amer. Math. Soc., Providence, RI, 1994.


\bibitem[Mil09]{Milne09}
J.S. Milne.
\newblock Points on Shimura varieties over finite fields: the conjecture of Langlands and Rapoport.
\newblock  {\url{http://arxiv.org/abs/0707.3173}}.


\bibitem[Moo98]{Moonen98}
B. Moonen. 
\newblock Models of Shimura varieties in mixed characteristics. 
\newblock Galois representations in arithmetic algebraic geometry (Durham, 1996), 267-350, London Math. Soc. Lecture Note Ser., 254, Cambridge Univ. Press, Cambridge, 1998. 


\bibitem[Noo96]{Noot96}
R. Noot.
\newblock Models of Shimura varieties in mixed characteristic.
\newblock J. Algebraic Geom.  5 (1996), no. 1, 187-207.


\bibitem[Oor13]{Oort13}
F. Oort. 
\newblock Moduli of abelian varieties in mixed and in positive characteristic. 
\newblock Handbook of moduli. Vol. III, 75-134, Adv. Lect. Math. (ALM), 26, Int. Press, Somerville, MA, 2013.


\bibitem[PR94]{PR94}
V. Platonov, A. Rapinchuk.
\newblock Algebraic groups and number theory. 
\newblock Pure and Applied Mathematics, 139, Academic Press, Inc., Boston, MA, 1994.


\bibitem[RR96]{RR96}
M. Rapoport, M. Richartz.
\newblock On the classification and specialization of $F$-isocrystals with additional structure.
\newblock Compositio Math. 103 (1996), no. 2, 153-181.


\bibitem[Rap03]{Rapoport03}
M. Rapoport. 
\newblock On the Newton stratification. 
\newblock S\'eminaire Bourbaki. Vol. 2001/2002. Ast\'erisque No. 290 (2003), Exp. No. 903, viii, 207-224.


\bibitem[Rap05]{Rapoport05}
M. Rapoport. 
\newblock A guide to the reduction modulo p of Shimura varieties. 
\newblock Automorphic forms. I. Ast\'erisque No. 298 (2005), 271-318.


\bibitem[SS13]{SS13}
P. Scholze, SW. Shin. 
\newblock On the cohomology of compact unitary group Shimura varieties at ramified split places. 
\newblock J. Amer. Math. Soc. 26 (2013), no. 1, 261-294.


\bibitem[Sh79]{Shelstad79}
D. Shelstad. 
\newblock Characters and inner forms of a quasisplit group over $\R$. 
\newblock Comp. Math. 39 (1979), 11-45.


\bibitem[Spr74]{Springer74}
T.A. Springer. 
\newblock Regular elements of finite reflection groups. 
\newblock Invent. Math. 25 (1974), 159-198. 




\bibitem[Spr98]{Springer98}
T.A. Springer.
\newblock Linear algebraic groups.
\newblock Progress in mathematics,  vol. 9, Birkh\"auser, Boston, Mass., 1998.


\bibitem[Sta97]{Stamm97}
H. Stamm.
\newblock On the reduction of the Hilbert-Blumenthal-moduli scheme with $\Gamma_0(p)$-level structure.
\newblock Forum Math. 9 (1997), no. 4, 405-455.


\bibitem[Tit79]{Tits79}
J. Tits. 
\newblock Reductive groups over local fields. 
\newblock Automorphic forms, representations and L-functions (Proc. Sympos. Pure Math., Oregon State Univ., Corvallis, Ore., 1977), Part 1, pp. 29-69, Proc. Sympos. Pure Math., XXXIII, Amer. Math. Soc., Providence, R.I., 1979.


\bibitem[Vas99]{Vasiu99}
A. Vasiu.
\newblock Integral canonical models of Shimura varieties of preabelian type.
\newblock Asian J. Math.  3  (1999), no. 2, 401-518.




\bibitem[Vas07]{Vasiu07}
A. Vasiu.
\newblock Good Reductions of Shimura Varieties of Hodge Type in Arbitrary Unramified Mixed Characteristic, Parts I and II.
\newblock {\url{http://arxiv.org/abs/0707.1668}} and {\url{http://arxiv.org/abs/0712.1572}}.


\bibitem[Vas08]{Vasiu08}
A. Vasiu.
\newblock Geometry of Shimura varieties of Hodge type over finite fields.
\newblock Proceedings of the NATO Advanced Study Institute on
``Higher dimensional geometry over finite fields", G\"ottingen,
Germany (June 25 - July 06, 2007), 197-243, IOS Press, 2008.




\bibitem[Vas11]{Vasiu11}
A. Vasiu. 
\newblock Manin problems for Shimura varieties of Hodge type. 
\newblock J. Ramanujan Math. Soc. 26 (2011), no. 1, 31-84.


\bibitem[Vas12]{Vasiu12}
A. Vasiu.
\newblock A motivic conjecture of Milne. 
\newblock J. Reine Angew. Math. 685 (2013), 181-247.


\bibitem[VW13]{VW13}
E. Viehmann, T. Wedhorn. 
\newblock Ekedahl-Oort and Newton strata for Shimura varieties of PEL-type. 
\newblock Math. Ann. 356 (2013), 1493-1550.



\bibitem[Vos98]{Voskresenskii98}
V. E. Voskresenskii. 
\newblock Algebraic groups and their birational invariants. 
\newblock Translations of Mathematical Monographs, 179. American Mathematical Society, Providence, RI, 1998. xiv+218.


\bibitem[Wed99]{Wedhorn99}
T. Wedhorn.
\newblock Ordinariness in good reductions of Shimura varieties of PEL-type.
\newblock Ann. Sci. \'Ecole Norm. Sup. (4)  32  (1999),  no. 5, 575-618.


\bibitem[Win05]{Wintenberger05}
J.-P. Wintenberger. 
\newblock Existence de F-cristaux avec structures suppl\'ementaires. 
\newblock Adv. Math. 190 (2005), no. 1, 196-224.


\bibitem[Wor13]{Wortmann13}
D. Wortmann.
\newblock The $\mu$-ordinary locus for Shimura varieties of Hodge type.
\newblock {\url{http://arxiv.org/abs/1310.6444}}.


\bibitem[Yu05]{Yu05}
C.-F. Yu. 
\newblock On the slope stratification of certain Shimura varieties. 
\newblock Math. Z. 251 (2005), no. 4, 859-873.


\end{thebibliography}
\end{document}